\documentclass[reqno,11pt]{article}
\usepackage[T1]{fontenc}
\usepackage[utf8]{inputenc}
\usepackage{comment}
\usepackage[a4paper,margin=1in]{geometry}
\usepackage{amsmath,amssymb,amsthm,mathtools,color}
\usepackage{mathrsfs}
\usepackage{bbm}
\usepackage{hyperref}

\usepackage{fancyhdr}
\usepackage{ifthen}

\newtheorem{theorem}{Theorem}
\newtheorem{lemma}[theorem]{Lemma}
\newtheorem{remark}[theorem]{Remark}
\theoremstyle{remark}
\newtheorem*{remark*}{Remark}
\theoremstyle{plain}
\newtheorem{corollary}[theorem]{Corollary}
\newtheorem{proposition}[theorem]{Proposition}
\newtheorem{step}{Step}

\newtheorem{example}{Example}

\newcommand{\R}{\mathbb{R}}
\newcommand{\supp}{\operatorname{supp}}

\newcommand{\pv}{\operatorname{p.v.}}

\numberwithin{equation}{section}

\hypersetup{
  pdftitle={A variational approach to an anisotropic nonlocal porous medium equation with drift},
  pdfauthor={Jose A. Carrillo; Ricardo Guimaraes; David Padilla-Garza; Ruiwen Shu}
}

\title{A variational approach to an anisotropic nonlocal porous medium equation with drift}
\author{
\textbf{Jos\'e A.\ Carrillo}\thanks{%
\protect\raggedright
Mathematical Institute, University of Oxford,
Oxford, UK.
Email: \texttt{jose.carrillo@maths.ox.ac.uk}.
}
\and
\textbf{Ricardo Guimar\~aes}\thanks{Universidade Estadual de Campinas, Campinas-SP, Brazil. Email: \texttt{rmmg1315@unicamp.br}}
\and
\textbf{David Padilla-Garza}\thanks{ \protect\raggedright Institute of Science and Technology Austria, Klosterneuburg, Austria. Email: \texttt{David.Padilla-Garza@ist.ac.at}}
\and
\textbf{Ruiwen Shu}\thanks{Department of Mathematics, University of Georgia, Athens, Georgia, USA. Email: \texttt{ruiwen.shu@uga.edu}}
}

\date{} 

\begin{document}
\maketitle

 \begin{abstract}
We study nonlocal porous-medium equations driven by anisotropic Riesz-type repulsive interactions and an external confining potential. The pressure is determined by convolution with an interaction potential of the form
\[
W(x)=|x|^{-s}\Omega\left(\frac{x}{|x|}\right),
\qquad d-2<s<d,
\]
where \(\Omega\) is a smooth, even, and strictly positive angular profile. For quadratic confinement {and
\(s\in(d-2,d)\cap(0,5]\)}, under suitable positivity assumptions on the anisotropic kernels, we prove that the unique energy minimizer is the only steady state among bounded continuous probability densities. We then establish the long-time convergence of global weak solutions to this equilibrium, both locally and globally in H\"older norms. To obtain these results, we derive uniform \(L^p\), \(L^\infty\), and H\"older estimates, adapting the regularity theory for fractional porous-medium equations to the anisotropic setting. We also construct global weak solutions through a minimizing-movement scheme in Wasserstein space, establish the corresponding energy-dissipation inequality, and prove a Talagrand-type transport inequality.
\end{abstract}

\section{Introduction}
We study a nonlocal porous-medium equation with external drift, in which the
pressure is determined by the density through convolution with an interaction
potential. More precisely, we consider
\begin{equation}\label{eq:our-evolution}
\begin{cases}
    \partial_t u
    =
    \nabla\cdot\Bigl(u\nabla\bigl((W*u)+U\bigr)\Bigr),\\[2mm]
    u(0)=u_0,
\end{cases}
\end{equation}
where \(W:\mathbb R^d\to\mathbb R\cup\{\infty\}\) is a repulsive interaction potential and
\(U:\mathbb R^d\to\mathbb R\) is a confining potential. Formally,
\eqref{eq:our-evolution} describes the mean-field evolution of particles moving
with velocity \(-\nabla(W*u+U)\), according to Darcy's law. It can also be
interpreted (see \cite{JKO98,CMV,CMV2,AGS05,villani2009optimal}) as the Wasserstein gradient flow of the interaction and confinement energy
\begin{equation}\label{eq:E_intro}
\mathcal E_U[\rho]
=
\frac12\int_{\mathbb R^d}(W*\rho)(x)d\rho(x)
+
\int_{\mathbb R^d}U(x)d\rho(x).
\end{equation}
Although several results are established for general confining potentials, our
main interest is the quadratic confinement \(U(x)=|x|^2\). Here, we consider anisotropic Riesz interaction potentials of the form
\begin{equation}\label{eq:W}
W(x)
=
|x|^{-s}\Omega\left(\frac{x}{|x|}\right),
\qquad x\neq0,
\end{equation}
where \(d\geq2\), \(d-2<s<d\), and
\(\Omega\in C^\infty(\mathbb S^{d-1})\) is even and strictly positive. The
angular profile \(\Omega\) allows the strength of the repulsive interaction to
depend on direction and therefore breaks the radial symmetry of the classical
Riesz interaction potential.

When \(U\equiv0\), equation \eqref{eq:our-evolution} reduces to an anisotropic
version of the fractional porous-medium equation with nonlocal pressure. In the isotropic case \(\Omega\equiv1\), the interaction potential \(W(x)=|x|^{-s}\) is,
up to normalization, the Riesz interaction potential associated with
\((-\Delta)^{-(d-s)/2}\). Thus, setting
\(\sigma:=(d-s)/2\in(0,1)\), the interaction potential generated by $u$, i.e. $W*u$, can be written as $(-\Delta)^{-\sigma}u$, up to a constant multiple. For the isotropic unconfined equation, \cite{caffarelli2011nonlinear} established the existence of weak
solutions, mass conservation, energy estimates, and finite propagation for
compactly supported data. A variational construction based on the Wasserstein
gradient-flow structure was later developed in \cite{Lisini_Mainini_Segatti_2018}. Boundedness and
\(L^p\)-decay estimates were obtained in
\cite{caffarelli2013regularity,Lisini_Mainini_Segatti_2018}, while Hölder
regularity was proved in \cite{caffarelli2013regularity} for
\(\sigma\neq1/2\), and in \cite{caffarelli2016regularity} for the critical
case \(\sigma=1/2\). This case is critical since the resulting diffusion has
the same order as the transport term, and the oscillation argument requires
successive geometric corrections to control the long-range part of the
velocity. 

The large-time behavior of the isotropic unconfined flow was studied in
\cite{caffarelli2010asymptotic}. In that setting, a self-similar change of
variables transforms the equation into a nonlocal Fokker--Planck equation, whose
stationary profiles are characterized through a fractional obstacle problem with quadratic confinement. Convergence to these profiles was then obtained by entropy methods. The regularity of minimizers of the energy \eqref{eq:E_intro} for the confined Riesz interaction potential with general confinement potentials, i.e. $\Omega=1$, was studied in \cite{CDM} via its connection to fractional obstacle problems under suitable assumptions on $U$.

\medskip

{\bf 1. Anisotropic potentials: steady states.-}
Anisotropic interaction potentials were considered in the literature due to their applications to two dimensional dislocations \cite{GPPS}. The authors in \cite{MRS,CMMRSV1} showed for the first time that anisotropic logarithmic interaction potentials in two dimensions lead to nonradial minimizers whose support is an ellipse or a degenerate measure supported on a vertical segment. Generalizations of these results were considered in \cite{CMMRSV2,MMRSV1}.

In the present work, we consider anisotropic potentials $W$ of the form \eqref{eq:W} with a quadratic confinement. The corresponding gradient flow \eqref{eq:our-evolution} has the same scaling as the quadratically confined Riesz potential, but we will show that the anisotropy replaces the radial equilibrium profile by an ellipsoidal one that requires a separate uniqueness theory for steady states. Under assumption \eqref{hw} and the assumptions on \(U\) in
Theorem\ref{thm:equilibrium-measure}, the energy \(\mathcal E_U\)
admits a unique minimizer on \(\mathcal P_2(\mathbb R^d)\).
Uniqueness follows from the conditional strict positive
definiteness of \(W\), which implies strict convexity of the
interaction energy along linear interpolations of probability
measures with finite energy; see \cite{ST},
\cite[Lemmas 1.7 and 1.9]{shu2025existence} and the references therein.
We call this minimizer the equilibrium measure and denote it
by \(\overline{\mu}\). The equilibrium measure associated with the
quadratically confined energy and anisotropic Riesz interaction potentials has been the object of extensive analysis. We first notice that the minimizers for anisotropic Riesz potentials are generically given by integrable densities $\overline{\mu} := \overline{\rho} dx$ under suitable assumptions. More precisely, in the isotropic case \(\Omega\equiv1\), the
unique minimizer over \(\mathcal P(\mathbb R^d)\), for
\(s\in(\max\{d-4,0\},d)\), is the Barenblatt-type density
\begin{equation*}
\rho_{\mathrm{iso}}(x)
=
C_{\mathrm{iso}}
\left(
1-\frac{|x|^2}{R_{\mathrm{iso}}^2}
\right)_+^{\frac{s+2-d}{2}},
\end{equation*}
where \(C_{\mathrm{iso}},R_{\mathrm{iso}}>0\) depend on \(d\), \(s\), and the
normalization of the energy, and are chosen so that
\(\rho_{\mathrm{iso}}\) has unit mass, see
\cite{caffarelli2010asymptotic,Carrillo_Shu_2023,Frank_22}. In our range \(d-2<s<d\), the
exponent \((s+2-d)/2\) lies in \((0,1)\), so
\(\rho_{\mathrm{iso}}\) is continuous, strictly positive in
\(B_{R_{\mathrm{iso}}}(0)\), and supported on
\(\overline{B_{R_{\mathrm{iso}}}(0)}\).

In the anisotropic setting, the ball supporting the isotropic minimizer is
deformed into an ellipsoid. A broad class of interaction energies was studied in \cite{Carrillo_Shu_2023} and they proved that the nonnegativity of the Fourier transform of the interaction potential is equivalent to linear interpolation convexity (LIC) leading to the uniqueness of minimizers. In two dimensions, it was shown in \cite{CS24cpam} that a large family of anisotropic Riesz-potentials have their minimizers supported on an ellipse and they are explicitly  given by a linear deformation of the isotropic profile. Related explicit characterizations were later obtained in three dimensions, first under additional symmetry assumptions and subsequently, without such assumptions, for the Riesz and Coulomb interactions; see \cite{Carrillo_Shu_3D,Mateu_Mora_Rondi_Scardia_Verdera_2023,frank2026explicit}. Related uniqueness results beyond a global LIC condition have been obtained for isotropic attractive--repulsive interaction energies. In \cite{shu2025extended}, the LIC radius quantifies strict convexity along linear interpolations of centered probability measures supported in a prescribed ball. This yields uniqueness up to translation when the LIC radius is large enough.

In the ranges $s\in [0,2)$, $d=2$ and \(s\in[d-3,d)\cap(0,5]\), $d\geq 3$, {where we interpret $s =0$ to correspond to a logarithmic singularity}, covered by
\cite{CS24cpam,Carrillo_Shu_3D,Mateu_Mora_Rondi_Scardia_Verdera_2023,frank2026explicit}, positivity of both \(W\) and its Fourier transform
implies that there exists a nonsingular matrix
\(M\in\mathbb R^{d\times d}\) such that
\begin{equation*}
\overline{\rho}(x)
=
|\det M|\rho_{\mathrm{iso}}(Mx)
=
|\det M|C_{\mathrm{iso}}
\left(
1-\frac{|Mx|^2}{R_{\mathrm{iso}}^2}
\right)_+^{\frac{s+2-d}{2}},
\end{equation*}
where \(|\det M|\) preserves the total mass, while
\(\operatorname{supp}\overline{\rho}\) is an ellipsoid centered at the origin. The
density \(\overline{\rho}\) satisfies the Euler--Lagrange conditions: there exists
\(C_0\in\mathbb R\) such that
\begin{equation*}
(W*\overline{\rho})(x)+|x|^2=C_0
\quad\text{on }\operatorname{supp}\overline{\rho},
\qquad
(W*\overline{\rho})(x)+|x|^2\geq C_0
\quad\text{in }\mathbb R^d.
\end{equation*}

The first main result of our paper, Theorem \ref{thm:uniqueness}, shows that steady states of \eqref{eq:our-evolution} with quadratic confinement under mild regularity assumptions are unique and characterized by the equilibrium measure $\overline{\mu} := \overline{\rho} dx$. This is achieved under certain positivity conditions on the anisotropic potential $W$, see \eqref{hl} and the definition of the operator $L$ in next section. More precisely, we prove that \(\overline{\rho}\) is the unique steady-state density among bounded continuous probability densities. The proof combines a maximum principle for the operator \(L\), with a sliding argument for the corresponding effective potentials.

\medskip

{\bf 2. Anisotropic potentials: Long-time asymptotics.-}
In order to analyse the evolution problem, we introduce some operators associated with the anisotropic interaction potential $W$. In particular, we obtain singular-integral representations for the operator \(L\), defined as the inverse of the operator $f \mapsto W \ast f$, and for \(A:=-\Delta L^{-1}\), and introduce the three positivity conditions \eqref{hw}, \eqref{hl}, and \eqref{ha} which will be necessary for our main theorems. They play distinct roles in the analysis: \eqref{hw} ensures that the interaction energy given by equation \eqref{eq:E_intro} is strictly positive-definite, \eqref{hl} yields the maximum principle underlying the uniqueness theory for steady states, and \eqref{ha} ensures that the singular kernel associated with \(A\) is uniformly comparable to \(|h|^{-s-2}\), which is used in the \(L^p\), \(L^\infty\), and H\"older estimates for the evolution problem. Moreover, the only general implications among these conditions established here are \eqref{ha} implies $\tilde{\Omega}>0$ and \(\eqref{hl}\) implies \(\eqref{hw}\), while no further implication holds in general; see Section \ref{sec:preliminaries}. We also provide explicit classes of admissible anisotropies satisfying these assumptions, including sufficiently small smooth perturbations of the isotropic interaction.

The second main result of this work, Theorem \ref{thm:asymptotic-convergence}, shows that weak solutions to \eqref{eq:our-evolution} with bounded probability measures as initial data converge to the unique anisotropic steady state $\overline{\rho}$ for large times obtained in \cite{frank2026explicit}. This is achieved by assuming the positivity conditions \eqref{hl} and \eqref{ha} on the interaction potential $W$. In order to prove this result, we need to provide a well-posedness theory in suitable spaces, show propagation of \(L^\infty\) and H\"older bounds uniform in time, and finally compactness for large times of the solutions to \eqref{eq:our-evolution}. Let us elaborate more on this proof strategy.

In fact, Theorem \ref{thm:global-wellposedness} constructs global weak solutions of \eqref{eq:our-evolution} through a minimizing-movement scheme in the quadratic Wasserstein distance. The result covers \(C^2\) confining potentials with quadratic growth and Laplacian bounded above. Adapting the variational framework of \cite{foghem2025gradient}, we use flow interchange to obtain Sobolev estimates for the discrete solutions. These estimates provide the compactness needed to pass to the limit in the discrete Euler--Lagrange equations. A De Giorgi variational interpolation yields the energy-dissipation inequality.

To use the characterization of steady states, Theorem \ref{thm:uniqueness}, in the analysis of the long-time behavior, we first need to derive uniform-in-time regularity estimates for the global solutions, which are of independent interest. Under \eqref{ha}, and under additional regularity assumptions on the confining potential $U$ we establish uniform \(L^p\) bounds for \(1<p<\infty\). If \(U(x)=|x|^2\), \(L^\infty\) and H\"older bounds are first proved for the corresponding unconfined flow and then transferred to the confined equation through a self-similar change of variables. By adapting the truncation argument and De Giorgi iteration in
\cite[Theorem 7.2]{Lisini_Mainini_Segatti_2018} to the anisotropic Dirichlet
form associated with \(A\), we show that solutions \(\rho\) of the unconfined equation 
$$
\partial_\tau\rho
=
\nabla_y\cdot\bigl(\rho\nabla_y(W*\rho)\bigr)
\quad
\mbox{satisfy}
\quad
\|\rho(\tau)\|_{L^\infty(\mathbb R^d)}
\leq
C_\infty \tau^{-\frac{d}{s+2}},
\qquad \tau>0.
$$
The exponent \(d/(s+2)\) is dictated by the natural scaling of the unconfined flow and coincides with the sharp decay rate in the isotropic setting. The self-similar transformation then gives an explicit smoothing bound for the confined solution. In particular, its \(L^\infty\)-norm remains uniformly bounded for large times by a constant independent of the initial \(L^\infty\)-norm.

We use the same strategy to establish H\"older regularity. For the unconfined
flow, we adapt the oscillation method in \cite{caffarelli2013regularity} when
\(\alpha{:=s-d+2}\neq1\), and the critical iteration in
\cite{caffarelli2016regularity} when \(\alpha=1\), to the anisotropic operator
\(A\). This yields local \(C^\eta\)-regularity for bounded solutions at every
positive time, for some \(\eta\in(0,1)\).
The original oscillation theory provides estimates only away from the initial
time. In our setting, the \(L^\infty\)-decay and the scaling invariance of the
unconfined flow allow us to upgrade these local estimates to global spatial
H\"older bounds. A self-similar change of variables then transfers these bounds
to the confined equation. The resulting regularity theory covers the full range
\(\alpha\in(0,2)\) and provides estimates that are uniform for large times. This is postponed to Appendix \ref{section_appendix} since the arguments follow the blueprint of \cite{caffarelli2013regularity}.

The proof of Theorem \ref{thm:asymptotic-convergence} can be finally obtained as follows.
The uniform \(L^\infty\) and H\"older estimates provide the compactness needed
to study the long-time behavior. Combined with the energy-dissipation
inequality, they allow us to consider time-translated solutions whose
asymptotic limits have zero dissipation and are therefore bounded continuous
steady states. The uniqueness theorem identifies every such limit with the
equilibrium density \(\overline{\rho}\), which in turn yields convergence of
the entire trajectory.

We first obtain strong H\"older convergence locally on every compact subset
of \(\mathbb R^d\), for every exponent below the H\"older exponent provided
by the regularity theory. Using the global H\"older bounds together with the
regularity of the explicit equilibrium profile, we then upgrade this result
to strong H\"older convergence on the whole space, with the admissible
exponents limited by both the regularity of the solution and the
\(C^{\alpha/2}\)-regularity of the equilibrium density. In particular, the
convergence is uniform in \(\mathbb R^d\) giving the final desired convergence results in Theorem \ref{thm:asymptotic-convergence}.

\medskip

{\bf 3. Organization of the paper.-} In Section \ref{sec:main-results}, we
introduce the standing assumptions and state the main results. The following
section collects preliminary properties of the anisotropic operators, together
with the existence and uniqueness of the equilibrium measure and the maximum
principle. We then prove the uniqueness of bounded continuous steady states.
The regularity analysis is divided into the \(L^p\) estimates of
Subsection \ref{subsec:Lp-estimates} and the \(L^\infty\) and H\"older estimates
of Subsection \ref{subsec:Linfty-holder-estimates}. These estimates are then
combined with the energy-dissipation inequality to establish convergence to
the equilibrium density as \(t\to\infty\). Subsequently, we develop the
minimizing-movement construction and prove the gradient-flow theorem, including
the convergence of the discrete solutions, the weak formulation, and the
energy-dissipation inequality. Finally, we prove the Talagrand-type cost
inequality. Appendix \ref{section_appendix} contains the adaptation of the
Caffarelli--Soria--V\'azquez oscillation method to the anisotropic operator
\(A=-\Delta L^{-1}\).

\section{Main results}\label{sec:main-results}

Before stating the main results, we introduce the notation and the
anisotropic operators used throughout the paper. We fix \(d\geq2\) and
\(d-2<s<d\).

The interaction potential \(W\) is an anisotropic Riesz-type potential of the form
\begin{equation*}
W(x)
=
|x|^{-s}
\Omega\left(\frac{x}{|x|}\right),
\qquad x\neq0.
\end{equation*}
Here, \(\Omega\in \mathcal{C}^\infty(\mathbb S^{d-1})\) is an even function,
and the Fourier transform of \(W\) has the form
\begin{equation*}
\widehat W(\xi)
=
|\xi|^{-d+s}
\widetilde\Omega\left(\frac{\xi}{|\xi|}\right),
\qquad \xi\neq0,
\end{equation*}
where \(\widetilde\Omega\in C^\infty(\mathbb S^{d-1})\).

Given a confining potential \(U:\mathbb R^d\to\mathbb R\), we define the energy
functional
\begin{equation*}
\mathcal E_U[\rho]
=
\frac12\int_{\mathbb R^d}(W*\rho)(x)d\rho(x)
+
\int_{\mathbb R^d}U(x)d\rho(x),
\end{equation*}
and the associated effective potential
\begin{equation*}
V_U[\rho](x)
:=
(W*\rho)(x)+U(x).
\end{equation*}
In the quadratic case \(U(x)=|x|^2\), we simply write
\[
\mathcal E[\rho]:=\mathcal E_{|\cdot|^2}[\rho],
\qquad
V[\rho]:=V_{|\cdot|^2}[\rho].
\]

We use two operators associated with \(W\).
{Whenever \(\widetilde\Omega\) is nowhere zero on \(\mathbb S^{d-1}\),
we define the Fourier multiplier \(L\) by}
\begin{equation}\label{defi_L_multiplier_main}
\widehat{L[f]}(\xi)
=
\psi(\xi)\widehat f(\xi)
:=
|\xi|^{d-s}
\widetilde\Omega\left(\frac{\xi}{|\xi|}\right)^{-1}
\widehat f(\xi),
\qquad \xi\neq0,
\end{equation}
{The positivity assumptions imposed below ensure this nonvanishing condition
in all applications. By construction,} \(L[W*f]=f\) in \(\mathcal S'(\mathbb R^d)\), whenever the convolution is
well defined. We also use the notation \(v=L^{-1}u\) to mean the Fourier
multiplier relation
\[
\widehat v(\xi)
=
\psi(\xi)^{-1}\widehat u(\xi),
\]
whenever the right-hand side is well defined.

In Lemma \ref{lemma_L_integral}, we prove that the operator \(L\) admits a
symmetric singular-integral representation
\begin{equation}\label{eq:L-kernel-main}
L[f](x)
=
\operatorname{p.v.}
\int_{\mathbb R^d}
\bigl(f(x)-f(y)\bigr)K_L(x-y)dy,
\end{equation}
where
\begin{equation*}
K_L(h)
=
|h|^{s-2d}
\omega_L\left(\frac{h}{|h|}\right),
\qquad h\neq0,
\end{equation*}
for an even function
\(\omega_L\in C^\infty(\mathbb S^{d-1})\).

{The second operator is defined directly by
\[
A[f]:=-\Delta[W*f],
\]
whenever the convolution is well defined. When \(L\) is defined,
this is equivalently \(A=-\Delta L^{-1}\).} Using the convention
\(\widehat{-\Delta f}(\xi)=|\xi|^2\widehat f(\xi)\) and setting
\(\alpha:=s-d+2\in(0,2)\), the Fourier symbol of \(A\) is
\begin{equation*}
\widehat{A[f]}(\xi)
=
\widetilde\psi(\xi)\widehat f(\xi)
:=
|\xi|^\alpha
\widetilde\Omega\left(\frac{\xi}{|\xi|}\right)
\widehat f(\xi),
\qquad \xi\neq0.
\end{equation*}
{Whenever \(L\) is defined, the symbols satisfy}
\[
\widetilde\psi(\xi)
=
|\xi|^2\psi(\xi)^{-1},
\qquad \xi\neq0.
\]

Analogously to \(L\) (see Lemma \ref{lemma_L_integral}), the operator \(A\)
admits a symmetric singular-integral representation
\begin{equation*}
A[f](x)
=
\operatorname{p.v.}
\int_{\mathbb R^d}
\bigl(f(x)-f(y)\bigr)K_A(x-y)dy,
\end{equation*}
where
\begin{equation*}
K_A(h)
=
|h|^{-d-\alpha}
\omega_A\left(\frac{h}{|h|}\right),
\qquad h\neq0,
\end{equation*}
for an even function
\(\omega_A\in C^\infty(\mathbb S^{d-1})\).

We now introduce the positivity conditions used throughout the paper. They
concern different objects associated with the interaction potential and play
distinct roles in the analysis.

\begin{equation}
\Omega(\theta)>0,
\qquad
\widetilde\Omega(\theta)>0
\qquad
\text{for every }\theta\in\mathbb S^{d-1}.
\tag{\textbf{H\(_W\)}}
\label{hw}
\end{equation}

\begin{equation}
\omega_L(\theta)>0
\qquad
\text{for every }\theta\in\mathbb S^{d-1}.
\tag{\textbf{H\(_L\)}}
\label{hl}
\end{equation}

\begin{equation}
\omega_A(\theta)>0
\qquad
\text{for every }\theta\in\mathbb S^{d-1}.
\tag{\textbf{H\(_A\)}}
\label{ha}
\end{equation}

Assumption \eqref{hw} provides the positivity properties needed for the
energy, \eqref{hl} yields the maximum principle used in the uniqueness theory
for steady states, and \eqref{ha} provides the ellipticity required for the
\(L^p\), \(L^\infty\), and H\"older estimates for the evolution problem.

The only general implication between these positivity conditions established
here is that \eqref{hl} implies \eqref{hw}; see
Lemma \ref{lem:HL-implies-HW}. Apart from this implication, the conditions
are in general independent, as illustrated by the counterexamples discussed
in Section \ref{sec:preliminaries}. Nevertheless, the class of interaction
potentials satisfying \eqref{hl} and \eqref{ha} is nonempty. In particular,
sufficiently small smooth perturbations of the isotropic interaction satisfy
all three conditions; see Proposition \ref{prop:small-perturbations}.

Since \(\omega_A\) is continuous and strictly positive on the compact sphere
\(\mathbb S^{d-1}\), there exist constants
\(0<\lambda_A\leq\Lambda_A<\infty\) such that
\begin{equation*}
0<\lambda_A
\leq
\omega_A(\theta)
\leq
\Lambda_A<\infty
\qquad
\text{for every }\theta\in\mathbb S^{d-1}.
\end{equation*}

With these structural assumptions in place, we now turn to the main results.
We first recall the functional framework used in the minimizing-movement
construction.

For a symmetric continuous symbol
\(\phi:\R^d\setminus\{0\}\to(0,\infty)\), we define
\[
\dot H_\phi(\R^d)
:=
\left\{
u\in\mathcal S'(\R^d):\ \widehat u\phi^{1/2}\in L^2(\R^d)
\right\},
\]
endowed with the homogeneous seminorm
\begin{equation*}
\|u\|_{\dot H_\phi(\R^d)}^2
:=
\frac{1}{(2\pi)^d}\int_{\R^d}|\widehat u(\xi)|^2\phi(\xi)d\xi.
\end{equation*}
We also define
\[
H_\phi(\R^d):=\dot H_\phi(\R^d)\cap L^2(\R^d),
\]
equipped with the norm
\[
\|u\|_{H_\phi(\R^d)}^2
:=
\|u\|_{L^2(\R^d)}^2+\|u\|_{\dot H_\phi(\R^d)}^2.
\]

We denote by \(\mathcal{P}(\mathbb{R}^d)\) the set of probability measures on
\(\mathbb{R}^d\), defined as
\[
\mathcal{P}(\mathbb{R}^d)
:=
\left\{
\mu \in \mathcal{M}_+(\mathbb{R}^d)
\ ; \
\mu(\mathbb{R}^d) = 1
\right\}.
\]
Furthermore, we consider the subset
\(\mathcal{P}_2(\mathbb{R}^d)\subset\mathcal{P}(\mathbb{R}^d)\) of probability
measures with finite second moment, given by
\[
\mathcal{P}_2(\mathbb{R}^d)
:=
\left\{
\mu \in \mathcal{P}(\mathbb{R}^d)
\ ; \
\int_{\mathbb{R}^d}|x|^2d\mu(x)<\infty
\right\}.
\]

For \(\mu,\nu\in\mathcal P_2(\R^d)\), we denote by \(W_2(\mu,\nu)\) the
quadratic Wasserstein distance, defined by
\[
W_2^2(\mu,\nu)
:=
\inf_{\gamma\in\Gamma(\mu,\nu)}
\int_{\R^d\times\R^d}|x-y|^2d\gamma(x,y),
\]
where \(\Gamma(\mu,\nu)\) is the set of transport plans with marginals
\(\mu\) and \(\nu\).  

We denote by
\[
AC^2(I,(\mathcal P_2(\R^d),W_2))
\]
the space of absolutely continuous curves
\(u:I\subset\R^+\to\mathcal P_2(\R^d)\) such that there exists
\(m\in L^2(I)\) satisfying
\[
W_2(u(t),u(s))
\leq
\int_s^t m(r)dr
\qquad
\text{for all }s,t\in I,\ s\leq t.
\]
Moreover, we denote by \(H\) the Boltzmann entropy
\begin{equation}\label{Boltzmann_defi}
H(u):=\int_{\R^d}u\log udx,
\end{equation}
with domain
\[
D(H):=\{u\in\mathcal P_2(\R^d):\ H(u)<\infty\}.
\]

Our first result provides a minimizing-movement construction of global weak
solutions to \eqref{eq:our-evolution} in the Wasserstein space. The proof,
given later in Section \ref{sec:proof-gradient-flow}, follows the strategy
developed for nonlocal pressure equations in \cite{foghem2025gradient}, which
in turn extends the gradient-flow framework of
\cite{Lisini_Mainini_Segatti_2018}.

\begin{theorem}[Gradient-flow structure]
\label{thm:global-wellposedness}
Assume the standing hypotheses $\eqref{hl}$ and $\eqref{ha}$.
Let \(U\in C^2(\R^d)\), and suppose that there exist positive
constants
\[
a_U,\ A_U,\ c_U,\ C_U,\ \Lambda_U
\]
such that
\begin{equation}
\label{eq:assumptions-U-gradient-flow}
c_U|x|^2-a_U
\le U(x)
\le C_U|x|^2+A_U
\qquad\text{for every }x\in\R^d,
\end{equation}
and
\begin{equation}
\label{eq:laplacian-upper-bound-U}
\Delta U(x)\le \Lambda_U
\qquad\text{for every }x\in\R^d.
\end{equation}
Let
\[
u_0\in
\dot H_{\psi^{-1}}(\R^d)\cap\mathcal P_2(\R^d).
\]
Then the following statements hold.

\begin{enumerate}
\item[(i)] \textbf{Existence of discrete solutions.}
For every \(\tau>0\), there exists a unique sequence
\((u_\tau^k)_{k\ge0}\), with \(u_\tau^0=u_0\), such that
\[
u_\tau^k\in
\dot H_{\psi^{-1}}(\R^d)\cap\mathcal P_2(\R^d)
\qquad\text{for every }k\ge0,
\]
and, for every \(k\ge1\),
\[
u_\tau^k\in
\operatorname*{argmin}_{\nu\in\mathcal P_2(\R^d)}
\left\{
\frac{1}{2\tau}W_2^2(\nu,u_\tau^{k-1})
+\mathcal E_U(\nu)
\right\}.
\]
Moreover, for every \(k\ge1\),
\[
u_\tau^k\in H_{\widetilde\psi}(\R^d)
\]
and
\begin{equation*}
\|u_\tau^k\|_{\dot H_{\widetilde\psi}(\R^d)}^2
\le
\frac{d}{\tau}+\Lambda_U.
\end{equation*}

\item[(ii)] \textbf{Convergence and regularity.}
Define the piecewise constant interpolation\footnote{Here
\(\lceil t/\tau\rceil\) denotes the ceiling function; hence
\(u_\tau(t)=u_\tau^k\) for
\(t\in((k-1)\tau,k\tau]\).}
by 
\[
u_\tau(t):=u_\tau^{\lceil t/\tau\rceil},
\qquad
u_\tau(0):=u_\tau^0.
\]
Then there exist a curve
\[
u\in AC^2\bigl([0,\infty),
(\mathcal P_2(\R^d),W_2)\bigr)
\]
and a subsequence \((\tau_n)_{n\ge1}\), with
\(\tau_n\to0\), such that
\[
u_{\tau_n}(t)\rightharpoonup u(t)
\qquad\text{narrowly for every }t\ge0.
\]
If, moreover, \(u_0\in D(H)\), then, for every \(T>0\),
\[
u_{\tau_n}\rightharpoonup u
\qquad\text{weakly in }
L^2(0,T;H_{\widetilde\psi}(\R^d)),
\]
and
\[
u_{\tau_n}\to u
\qquad\text{strongly in }
L^2(T_0,T;L^2_{\mathrm{loc}}(\R^d))
\]
for every \(0<T_0<T\).

\item[(iii)] \textbf{Convergence of the velocity.}
Define
\[
v_\tau:=L^{-1}u_\tau,
\qquad
v:=L^{-1}u.
\]
If \(u_0\in D(H)\), then, for every \(T>0\),
\[
\nabla v\in L^2(0,T;L^2(\R^d))
\]
and
\[
\nabla v_{\tau_n}\rightharpoonup\nabla v
\qquad\text{weakly in }
L^2(0,T;L^2(\R^d)).
\]

\item[(iv)] \textbf{Weak solution of the equation.}
If \(u_0\in D(H)\), then the limiting curve \(u\) is a weak
solution of
\[
\partial_tu
=
\nabla\cdot\bigl(u\nabla(v+U)\bigr),
\qquad
Lv=u,
\qquad
u(0)=u_0.
\]
More precisely, for every
\(\varphi\in C_c^\infty(\R^d\times(0,\infty))\),
\[
\int_0^\infty\int_{\R^d}
u
\Bigl(
\partial_t\varphi
-\nabla\varphi\cdot\nabla v
-\nabla U\cdot\nabla\varphi
\Bigr)dxdt
=0.
\]

\item[(v)] \textbf{Energy-dissipation inequality.}
If \(u\) is the weak solution from item \textup{(iv)}, then,
for every \(T>0\),
\[
\mathcal E_U(u(T))
+
\int_0^T\int_{\R^d}
u(t,x)
\left|
\nabla\bigl(v(t,x)+U(x)\bigr)
\right|^2
dxdt
\le
\mathcal E_U(u_0).
\]
\end{enumerate}
\end{theorem}

The next result is a Talagrand-type cost inequality. It bounds the Wasserstein distance to equilibrium in terms of the energy gap and therefore gives quantitative information on the geometry of the energy landscape near the minimizer.

\begin{proposition}[Talagrand-type cost inequality]
\label{prop:transport}
{Assume \eqref{hw}, and let \(U\in C(\R^d)\) have at least
quadratic growth: there exists \(c_U>0\) such that}
\begin{equation}
\label{eq:quadratic-growth}
U(x)\ge c_U(|x|^2-1)
\qquad\text{for every }x\in\R^d.
\end{equation}
{Let \(\overline{\mu}\) be the unique equilibrium measure
provided by Theorem~\ref{thm:equilibrium-measure}, and let
\(\nu\) be a probability measure on \(\R^d\).}
Then, there exists a positive constant \(C\) such that
\begin{equation}
\label{eq:transport}
W_2^4(\nu,\bar{\mu})
\le
C\left(
\mathcal E_U(\nu)-\mathcal E_U(\overline{\mu})
+
\left(\mathcal E_U(\nu)-\mathcal E_U(\overline{\mu})\right)^2
\right).
\end{equation}
\end{proposition}

We next characterize steady states. {For \(\nu=\rho dx\), where \(\rho\) is a continuous
probability density, set \(V_U[\rho]:=W*\rho+U\) and define
\begin{equation}
\label{eq:defV}
\mathcal V_U(\nu):=
\begin{cases}
\displaystyle
\int_{\{\rho>0\}}\rho(x)
\left|\nabla V_U[\rho](x)\right|^2dx,
&
V_U[\rho]\in H^1_{\mathrm{loc}}(\{\rho>0\}),\\[1ex]
+\infty, & \text{otherwise}.
\end{cases}
\end{equation}
Here \(\nabla V_U[\rho]\) denotes the distributional gradient
on the open set \(\{\rho>0\}\).
We call \(\nu\) a steady state if \(\mathcal V_U(\nu)=0\).}
This definition is natural because, for smooth solutions of \eqref{eq:our-evolution},
\[
\frac{d}{dt}\mathcal E_U(u(t))
=
-\mathcal V_U(u(t)).
\]
Thus steady states are precisely states with zero instantaneous energy dissipation. 

The equilibrium measure \(\overline{\mu}\) is clearly a steady state. To prove
that it is the only steady state among bounded continuous densities, we use the
explicit representation of its density. This representation will also be used
in Theorem \ref{thm:asymptotic-convergence} to identify the large-time limit and
to obtain global H\"older convergence.

For this purpose, we introduce the parameter range 
\begin{equation}\label{eq:explicit-minimizer-range}
\mathcal I_d :=(d-2,d)\cap(0,5].
\end{equation}
For \(s\in\mathcal I_d\), the equilibrium measure
\(\overline{\mu}=\overline{\rho}dx\) admits the explicit ellipsoidal
representation
\begin{equation}\label{eq:explicit-equilibrium-profile}
\overline{\rho}(x)
=
C\left(1-\frac{|Ax|^2}{R^2}\right)_+^{\alpha/2},
\qquad
\alpha=s-d+2\in(0,2),
\end{equation}
for some invertible matrix \(A\in\mathbb R^{d\times d}\), radius \(R>0\),
and constant \(C>0\); see
\cite{CS24cpam,Carrillo_Shu_3D,
Mateu_Mora_Rondi_Scardia_Verdera_2023,frank2026explicit}.
The following theorem shows that \(\overline{\mu}\) is the unique steady state
among bounded continuous probability densities.

\begin{theorem}[Uniqueness of bounded continuous steady states]
\label{thm:uniqueness}
Fix \(d\geq2\) and \(s\in\mathcal I_d\), where \(\mathcal I_d\) is defined in
\eqref{eq:explicit-minimizer-range}, and assume the standing hypothesis
\eqref{hl}. Let \(U(x)=|x|^2\), and let
\(\mathcal E_U\) and \(\mathcal V_U\) be defined by
\eqref{eq:E_intro} and \eqref{eq:defV}, respectively. If
\(\nu=\rho dx\) is a steady state in the sense that
\[
\mathcal V_U(\nu)=0,
\]
with
\[
\rho\in C(\mathbb R^d)\cap L^\infty(\mathbb R^d),
\qquad
\rho\geq0,
\qquad
\int_{\mathbb R^d}\rho(x)dx=1,
\]
then \(\rho=\overline{\rho}\) in \(\mathbb R^d\). Therefore,
\(\nu=\overline{\mu}=\overline{\rho}dx\).
\end{theorem}

We also obtain regularity estimates for solutions of \eqref{eq:our-evolution}. The next theorem gives uniform-in-time bounds in the \(L^p\), \(L^\infty\), and H\"older spaces.

\begin{theorem}[Regularity]
\label{thm:Lp_global_confined}
Fix \(d\geq2\) and \(d-2<s<d\). Assume
 \(\mathbf{(H_L)}\) and \(\mathbf{(H_A)}\). Let
\[
u_0\in
\dot H_{\psi^{-1}}(\mathbb R^d)
\cap\mathcal P_2(\mathbb R^d)
\cap D(H),
\]
and let \(u\) be the global weak solution given by
Theorem \ref{thm:global-wellposedness}\textup{(iv)}. Let
\(1<p<\infty\). Then the following statements hold.

\begin{enumerate}
\item[(i)] \textbf{\(L^p\) bounds.} If \(u_0\in L^p(\mathbb R^d)\), and \(U\in C^2(\mathbb R^d)\) is as in Theorem \ref{thm:global-wellposedness}, then
\[
\sup_{t\ge0}\|u(t)\|_{L^p(\mathbb R^d)}
\le
\max\left\{\|u_0\|_{L^p(\mathbb R^d)},C_p\right\},
\]
where \(C_p>0\) depends only on \(p\), on \(a_U,A_U,c_U,C_U,\Lambda_U\), and on the ellipticity constants in \eqref{ha} .

\item[(ii)] \textbf{\(L^\infty\) bounds.} If \(U(x)=|x|^2\), then for every \(t>0\),
\[
\|u(t)\|_{L^\infty(\mathbb R^d)}
\le
C_\infty e^{2dt}
\left(
\frac{e^{2(s+2)t}-1}{2(s+2)}
\right)^{-d/(s+2)}.
\]
In particular, for any $t > t_0$ with $t_0 > 0$, $\|u(t)\|_{L^\infty(\mathbb{R}^d)}$ is bounded by a constant independent of $\|u_0\|_{L^\infty(\mathbb{R}^d)}$.

\item[(iii)] \textbf{H\"older bounds.}
If \(U(x)=|x|^2\), then there exists
\(\eta=\eta(d,s,\Omega)\in(0,1)\) such that \(u\) is locally
H\"older continuous in space-time in
\(\mathbb R^d\times(0,\infty)\). More precisely, for every
compact set \(K\subset\mathbb R^d\) and every
\(0<t_0<t_1<\infty\), there exists a constant
\(C=C(d,s,\Omega,K,t_0,t_1)>0\) such that
\[
|u(t,x)-u(s,y)|
\leq
C\left(
|x-y|^\eta+|t-s|^\eta
\right)
\]
for every \(x,y\in K\) and \(t,s\in[t_0,t_1]\). Moreover, there exists a constant
\(C_\eta=C_\eta(d,s,\Omega)>0\), independent of time,
such that, for every \(t>0\),
\[
[u(t)]_{C^\eta(\mathbb R^d)}
\leq
C_\eta e^{2(d+\eta)t}
\left(
\frac{e^{2(s+2)t}-1}{2(s+2)}
\right)^{-\frac{d+\eta}{s+2}}.
\]
Hence, for every \(t_0>0\),
\[
\sup_{t\geq t_0}
\|u(t)\|_{C^\eta(\mathbb R^d)}
<\infty,
\]
with a bound depending only on \(d,s,\Omega,t_0\),
and independent of the initial \(L^\infty\)-norm.

\end{enumerate}
\end{theorem}

Finally, we combine the previous results to identify the asymptotic behavior of solutions in the case of a quadratic confining potential. The solution converges to the equilibrium density in a strong H\"older norm and, in particular, uniformly.

\begin{theorem}[Asymptotic convergence] \label{thm:asymptotic-convergence}
Fix \(d\geq2\) and \(s\in\mathcal I_d\), where \(\mathcal I_d\) is the interval defined in \eqref{eq:explicit-minimizer-range}, and assume that \eqref{hl} and \eqref{ha} hold. Set \(\alpha=s-d+2\in(0,2)\), and let \(U(x)=|x|^2\). Assume that
\[ 
u_0\in \dot H_{\psi^{-1}}(\mathbb R^d) \cap\mathcal P_2(\mathbb R^d) \cap D(H), 
\]
and denote by \(u\) the corresponding global weak solution to \eqref{eq:our-evolution} constructed in Theorem~\ref{thm:global-wellposedness}\textup{(iv)}. Further, denote the equilibrium measure by \(\overline{\mu}=\overline{\rho}dx\), and let \(\eta\in(0,1)\) be the H\"older exponent from Theorem \ref{thm:Lp_global_confined}\textup{(iii)}. Then, the following convergence results hold: for every compact set \(K\subset\mathbb R^d\) and every \(0<\beta<\eta\),
\[ 
u(t)\longrightarrow\overline{\rho} \qquad \text{strongly in }C^\beta(K) \quad\text{as }t\to\infty. 
\]
In addition, for every \(0<\beta<\eta_*:=\min\{\eta,\alpha/2\}\), the convergence holds globally:
\[ 
u(t)\longrightarrow\overline{\rho} \qquad \text{strongly in }C^\beta(\mathbb R^d) \quad\text{as }t\to\infty, 
\]
\end{theorem}

Before turning to the proofs of the main results, we collect in next Section several preliminary properties of the anisotropic operators and prove the auxiliary results needed throughout the paper. We then begin with the proof of the uniqueness theorem in Section \ref{sec:proof-uniqueness}, followed by the regularity and asymptotic convergence results. Finally, we prove the gradient-flow theorem in Section \ref{sec:proof-gradient-flow} and the Talagrand-type transport inequality in Section \ref{sec:proof-transport}.

\section{Preliminaries}\label{sec:preliminaries}

We begin with some foundational results that will be necessary for the proofs of the main results. 

\begin{theorem}[Existence and uniqueness of the equilibrium measure]
\label{thm:equilibrium-measure}
Assume \eqref{hw}, and let
\(U\in C(\mathbb R^d)\) satisfy
\[
U(x)\geq c_U|x|^2-a_U
\qquad\text{for every }x\in\mathbb R^d,
\]
for some constants \(c_U>0\) and \(a_U>0\). Then there exists a
unique measure \(\overline{\mu}\in\mathcal P_2(\mathbb R^d)\) such that
\[
\mathcal E_U(\overline{\mu})
=
\inf_{\mu\in\mathcal P(\mathbb R^d)}
\mathcal E_U(\mu).
\]
We call \(\overline{\mu}\) the equilibrium measure associated with \(U\).
\end{theorem}

\begin{proof}
Existence follows from the direct method. Uniqueness follows from the
linear-interpolation convexity argument in
\cite[Section 2]{Carrillo_Shu_2023} (especially (2.1) in this reference and the discussion immediately
following it): under \eqref{hw}, the strict
positivity of \(\widehat W\) makes the interaction term strictly convex
along every nontrivial linear interpolation of probability measures,
whereas the confinement term is linear.
\end{proof}

\begin{lemma}
\label{lemma_L_integral}
Let $\omega \in C^{\infty}(\mathbb{S}^{d-1})$ be an even function. Let $\gamma \in (0,2)$.
\begin{itemize}
    \item[1.] Define the kernel $K$ as
    \begin{equation}
    \label{eq:K}
        K(x)=|x|^{-(d+\gamma)}\omega\left(\frac{x}{|x|}\right).
    \end{equation}
    Given a Schwartz class function $u$, define the operator $G$ as
    \begin{equation}
    \label{eq:L-kernel}
        G[u](x)=\operatorname{p.v.}\int_{\mathbb R^d}\bigl(u(x)-u(y)\bigr)K(x-y)dy.
    \end{equation}
    Then
    \begin{equation*}
        \widehat{G[u]}(\xi)=\mathcal{RF}(K)(\xi)\widehat u(\xi),
    \end{equation*}
    where the renormalized Fourier transform $\mathcal{RF}$ is defined for an even input as
    \begin{equation*}
        \mathcal{RF}(K)(\xi)=\int_{\mathbb R^d}\bigl(1-\cos(z\cdot\xi)\bigr)K(z)dz.
    \end{equation*}

    \item[2.] The renormalized Fourier transform satisfies
    \begin{equation}
    \label{eq:firstform}
        \mathcal{RF}(K)(\xi)=C_{d,\gamma}|\xi|^\gamma
        \widetilde\omega\left(\frac{\xi}{|\xi|}\right),
    \end{equation}
    where
    \begin{equation*}
        C_{d,\gamma}:=\int_0^\infty(1-\cos u)u^{-1-\gamma}du,
    \end{equation*}
    and $\widetilde\omega\in C^\infty(\mathbb S^{d-1})$ is given by
    \begin{equation*}
        \widetilde\omega(\eta)=T_\gamma[\omega](\eta),\qquad \eta\in\mathbb S^{d-1},
    \end{equation*}
    where the operator $T_\gamma$ is defined as
    \begin{equation}
    \label{tgamma}
        T_\gamma[f](\eta):=\int_{\mathbb S^{d-1}}|\eta\cdot\theta|^\gamma f(\theta)d\theta.
    \end{equation}
    Equivalently, $\widetilde\omega$ may be written in terms of spherical harmonics: let
    \begin{equation}
    \label{eq:harmonics}
        \omega(\theta)=\sum_{l,m}a_{2l,m}Y_{2l,m}(\theta),
    \end{equation}
    so that only even spherical harmonics appear. Then
    \begin{equation*}
        \widetilde\omega(\theta)
        =\sum_{l,m}\lambda_{\gamma,2l}a_{2l,m}Y_{2l,m}(\theta),
    \end{equation*}
   where \(\lambda_{\gamma,2l}\neq0\) for every integer \(l\geq0\), and
    \begin{equation*}
        |\lambda_{\gamma,2l}|\asymp(1+l)^{-\frac d2-\gamma}.
    \end{equation*}

    \item[3.] Let $\widetilde\omega\in C^\infty(\mathbb S^{d-1})$ be an even function. Write $\widetilde\omega$ in terms of spherical harmonics as
    \begin{equation*}
        \widetilde\omega(\theta)=\sum_{l,m}\widetilde a_{2l,m}Y_{2l,m}(\theta).
    \end{equation*}
    Define $\omega\in C^\infty(\mathbb S^{d-1})$ as
    \begin{equation*}
        \omega(\theta)=\frac{1}{C_{d,\gamma}}
        \sum_{l,m}\frac{1}{\lambda_{\gamma,2l}}
        \widetilde a_{2l,m}Y_{2l,m}(\theta),
    \end{equation*}
    or equivalently
    \begin{equation*}
        \omega=\frac{1}{C_{d,\gamma}}T_\gamma^{-1}[\widetilde\omega].
    \end{equation*}
    Define $K$ as in equation \eqref{eq:K}, and the operator $G$ as in equation \eqref{eq:L-kernel}. Then
    \begin{equation}
    \label{eq:target}
        \widehat{G[u]}(\xi)=|\xi|^\gamma
        \widetilde\omega\left(\frac{\xi}{|\xi|}\right)\widehat u(\xi).
    \end{equation}
\end{itemize}
\end{lemma}

\begin{proof}
We prove each item separately.

\medskip

\noindent
\textbf{Proof of (1).}

Using the change of variables $z=y-x$, we obtain
\[
G[u](x)
=
\mathrm{p.v.}\int_{\mathbb R^d}(u(x)-u(x+z))K(z)dz.
\]
Since $u\in\mathcal S(\mathbb R^d)$ and $K$ is even, we may symmetrize the principal-value integral and write
\[
G[u](x)
=
\frac12\int_{\mathbb R^d}
\bigl(2u(x)-u(x+z)-u(x-z)\bigr)K(z)dz.
\]
By Taylor's theorem,
\[
\bigl\|2u(\cdot)-u(\cdot+z)-u(\cdot-z)\bigr\|_{L^1(\mathbb R^d)}
\leq C|z|^2
\qquad\text{for }|z|\leq1,
\]
while this $L^1$-norm is bounded by $4\|u\|_{L^1}$ for $|z|>1$. Since $|K(z)|\leq C|z|^{-d-\gamma}$, with $\gamma\in(0,2)$, it follows that
\[
\int_{\mathbb R^d}|K(z)|
\int_{\mathbb R^d}
\bigl|2u(x)-u(x+z)-u(x-z)\bigr|dxdz<\infty.
\]
Thus Fubini's theorem applies. Taking Fourier transforms gives
\begin{align*}
\widehat{G[u]}(\xi)
&=
\frac12\int_{\mathbb R^d}
\left(
\int_{\mathbb R^d}
 e^{-ix\cdot\xi}
 \bigl(2u(x)-u(x+z)-u(x-z)\bigr)dx
\right)K(z)dz \\
&=
\frac12\int_{\mathbb R^d}
\bigl(2-e^{iz\cdot\xi}-e^{-iz\cdot\xi}\bigr)K(z)dz\widehat u(\xi) \\
&=
\left(
\int_{\mathbb R^d}
(1-\cos(z\cdot\xi))K(z)dz
\right)\widehat u(\xi).
\end{align*}
This proves item \textup{(1)}.

\medskip

\noindent
\textbf{Proof of (2).}
By definition,
\[
\mathcal{RF}(K)(\xi)
=
\int_{\mathbb R^d}
(1-\cos(z\cdot\xi))
|z|^{-d-\gamma}
\omega\!\left(\frac{z}{|z|}\right)dz.
\]
Passing to polar coordinates $z=r\theta$, with $r>0$ and $\theta\in S^{d-1}$, we obtain
\[
\mathcal{RF}(K)(\xi)
=
\int_{S^{d-1}}
\omega(\theta)
\int_0^\infty
(1-\cos(r\xi\cdot\theta))
r^{-1-\gamma}drd\theta.
\]
For fixed $\theta$, set $a=|\xi\cdot\theta|$. Using the change of variables $u=ar$,
\[
\int_0^\infty
(1-\cos(r\xi\cdot\theta))
r^{-1-\gamma}dr
=
a^\gamma
\int_0^\infty
(1-\cos u)u^{-1-\gamma}du.
\]
Define
\[
C_{d,\gamma}
=
\int_0^\infty
(1-\cos u)u^{-1-\gamma}du.
\]
Then
\[
\mathcal{RF}(K)(\xi)
=
C_{d,\gamma}
\int_{S^{d-1}}
|\xi\cdot\theta|^\gamma
\omega(\theta)d\theta
=
C_{d,\gamma}
|\xi|^\gamma
\int_{S^{d-1}}
\left|
\frac{\xi}{|\xi|}\cdot\theta
\right|^\gamma
\omega(\theta)d\theta.
\]
This proves equation \eqref{eq:firstform}.

We now prove the spherical harmonic representation. Since the kernel defining $T_\gamma$ depends only on $\eta\cdot\theta$, the Funk--Hecke theorem \cite[Theorem 2.22]{atkinson2012spherical} implies that every spherical harmonic $Y_{l,m}$ is an eigenfunction of $T_\gamma$:
\[
T_\gamma[Y_{l,m}]
=
\lambda_{\gamma,l}Y_{l,m},
\]
where
\[
\lambda_{\gamma,l}
=
|S^{d-2}|
\int_{-1}^1
|t|^\gamma
P_l^{(d)}(t)
(1-t^2)^{\frac{d-3}{2}}dt,
\]
and $P_l^{(d)}$ denotes the Gegenbauer polynomial of degree $l$ appearing in the Funk--Hecke formula. Therefore, from equation \eqref{eq:harmonics},
\[
\widetilde\omega(\theta)
=
\sum_{l,m}
\lambda_{\gamma,2l}a_{2l,m}Y_{2l,m}(\theta).
\]
For $l\ge0$, the preceding integral can be evaluated explicitly as
\[
\lambda_{\gamma,2l}
=
2\pi^{\frac{d-1}{2}}
\frac{\Gamma\!\left(\frac{1+\gamma}{2}\right)}
{\Gamma\!\left(-\frac{\gamma}{2}\right)}
(-1)^l
\frac{\Gamma\!\left(l-\frac{\gamma}{2}\right)}
{\Gamma\!\left(l+\frac{d+\gamma}{2}\right)}.
\]
Since $\gamma\in(0,2)$, all the factors in this expression are finite and nonzero, and hence
\[
\lambda_{\gamma,2l}\neq0
\qquad\text{for every }l\ge0.
\]
Moreover, Stirling's formula gives
\[
|\lambda_{\gamma,2l}|
\asymp
(1+l)^{-\frac d2-\gamma}.
\]
Finally, since $\omega\in C^\infty(S^{d-1})$, its spherical harmonic coefficients decay rapidly. The preceding estimate shows that multiplication by $\lambda_{\gamma,2l}$ preserves rapid decay, and therefore
\[
\widetilde\omega\in C^\infty(S^{d-1}).
\]

\medskip

\noindent
\textbf{Proof of (3).}
Since $\widetilde\omega\in C^\infty(S^{d-1})$, the coefficients $\widetilde a_{2l,m}$ decay rapidly, while the estimate above shows that $\lambda_{\gamma,2l}^{-1}$ grows at most polynomially. Hence the series defining $\omega$ converges in $C^\infty(S^{d-1})$.

By the definition of $\omega$,
\[
T_\gamma[\omega]
=
\frac{1}{C_{d,\gamma}}\widetilde\omega.
\]
Applying Part \textup{(2)}, we obtain
\[
\mathcal{RF}(K)(\xi)
=
C_{d,\gamma}|\xi|^\gamma
T_\gamma[\omega]\!\left(\frac{\xi}{|\xi|}\right)
=
|\xi|^\gamma
\widetilde\omega\!\left(\frac{\xi}{|\xi|}\right).
\]
Using Part \textup{(1)},
\[
\widehat{G[u]}(\xi)
=
\mathcal{RF}(K)(\xi)\widehat u(\xi),
\]
and equation \eqref{eq:target} follows. This completes the proof.
\end{proof}

We now clarify the class of admissible anisotropic potentials satisfying the hypotheses $\eqref{hw}$, $\eqref{hl}$, and $\eqref{ha}$.

\begin{proposition}[Positivity for small anisotropic perturbations]
\label{prop:small-perturbations}
Let \(d\geq2\), fix \(\gamma\in(0,2)\), and set \(\delta:=2-\gamma\). 
Choose
\[
N>d-\frac12+\max\{\gamma,\delta\}.
\]
Then there exists \(\varepsilon_*>0\) such that, if \(\widetilde\Omega\in C^\infty(\mathbb S^{d-1})\) is even and
\[
\|\widetilde\Omega-1\|_{H^N(\mathbb S^{d-1})}<\varepsilon_*,
\]
then, denoting by \(\Omega\) the corresponding angular profile under the inverse homogeneous Fourier transform,
\[
\Omega>0,\qquad
\widetilde\Omega>0,\qquad
T_\gamma^{-1}[\widetilde\Omega]>0,\qquad
T_\delta^{-1}[\widetilde\Omega^{-1}]>0
\]
pointwise on \(\mathbb S^{d-1}\), where the inverse operators are understood on the subspace of even functions. Moreover, \(T_\gamma^{-1}[\widetilde\Omega]\) and \(T_\delta^{-1}[\widetilde\Omega^{-1}]\) belong to \(C^\infty(\mathbb S^{d-1})\). As a consequence, the hypotheses $\eqref{hw}$, $\eqref{hl}$, and $\eqref{ha}$ are satisfied. 
\end{proposition}

\begin{proof}
\textbf{Step 1.} We first claim that
\[
T_\kappa^{-1}:H^{q+\frac d2+\kappa}_{\mathrm{even}}(\mathbb S^{d-1})
\longrightarrow H^q_{\mathrm{even}}(\mathbb S^{d-1})
\]
is continuous for every \(q\in\mathbb R\). Indeed, for every \(\kappa\in(0,2)\), the Funk--Hecke eigenvalues of \(T_\kappa\) on spherical harmonics of even degree satisfy
\begin{equation}\label{ineq_eigenvalues}
|\lambda_{\kappa,2j}|\asymp(1+j)^{-\frac d2-\kappa}.
\end{equation}
 Let \(\mu_\ell=\ell(\ell+d-2)\) denote the eigenvalues of \(-\Delta_{\mathbb S^{d-1}}\). If \(f\) is even, then
\[
f=\sum_{j,m}f_{2j,m}Y_{2j,m},
\qquad
T_\kappa^{-1}f=\sum_{j,m}\lambda_{\kappa,2j}^{-1}f_{2j,m}Y_{2j,m}.
\]
By \eqref{ineq_eigenvalues} and the relation \(1+\mu_{2j}\asymp(1+j)^2\), we have
\[
|\lambda_{\kappa,2j}|^{-2}
\leq C(1+\mu_{2j})^{\frac d2+\kappa}.
\]
Therefore, we get
\[
\begin{aligned}
\|T_\kappa^{-1}f\|_{H^q}^2
&=
\sum_{j,m}(1+\mu_{2j})^q|\lambda_{\kappa,2j}|^{-2}|f_{2j,m}|^2\\
&\leq
C\sum_{j,m}(1+\mu_{2j})^{q+\frac d2+\kappa}|f_{2j,m}|^2
=
C\|f\|_{H^{q+\frac d2+\kappa}}^2.
\end{aligned}
\]
By the choice of \(N\), we may choose \(r>(d-1)/2\) such that $N>r+\frac d2+\max\{\gamma,\delta\}.$
Taking \(q=r\) and using the Sobolev embedding \(H^r(\mathbb S^{d-1})\hookrightarrow C^0(\mathbb S^{d-1})\), we obtain
\begin{equation}\label{eq:T-kappa-inverse-Sobolev-C0}
\|T_\kappa^{-1}f\|_{L^\infty(\mathbb S^{d-1})}
\leq C_\kappa\|f\|_{H^{r+\frac d2+\kappa}(\mathbb S^{d-1})}.
\end{equation}
Moreover, for \(\kappa\in\{\gamma,\delta\}\), the choice of \(N\) gives $H^N(\mathbb S^{d-1})
\hookrightarrow
H^{r+\frac d2+\kappa}(\mathbb S^{d-1}),$ and hence
\begin{equation}\label{eq:T-kappa-inverse-HN-C0}
\|T_\kappa^{-1}f\|_{L^\infty(\mathbb S^{d-1})}
\leq C\|f\|_{H^N(\mathbb S^{d-1})}.
\end{equation}

We also control the spatial angular profile. Let \(\mathcal F_\gamma\) denote the angular map induced by the homogeneous Fourier transform, so that \(\widetilde\Omega=\mathcal F_\gamma\Omega\), and let \(\Omega_{\mathrm{iso}}:=\mathcal F_\gamma^{-1}[1]>0\). The inverse multipliers of \(\mathcal F_\gamma\) on spherical harmonics of degree \(2j\) satisfy
\[
|b_{\gamma,2j}^{-1}|\asymp(1+j)^{\frac d2-2+\gamma}.
\]
Hence, choosing \(r>(d-1)/2\) as above and using Sobolev embedding,
\[
\|\Omega-\Omega_{\mathrm{iso}}\|_{L^\infty}
\leq C\|\widetilde\Omega-1\|_{H^N}.
\]
After decreasing \(\varepsilon_*\), this gives \(\Omega\geq \Omega_{\mathrm{iso}}/2>0\).

\textbf{Step 2.} We now claim that
\begin{equation}\label{eq:inverse-Omega-small-perturbation}
\|\widetilde\Omega^{-1}-1\|_{H^N}
\leq
C\|\widetilde\Omega-1\|_{H^N},
\end{equation}
for some constant $C$. 

\textbf{Substep 2.1} We first prove that, for \(N>(d-1)/2\), the space \(H^N(\mathbb S^{d-1})\) is a Banach algebra, that is, $||fg||_{H^N} \leq C_N ||f||_{H^N} ||g||_{H^N}$. 
Indeed, by the fractional Leibniz rule on the sphere,
\[
\|fg\|_{H^N(\mathbb S^{d-1})}
\leq C_N\Bigl(
\|f\|_{L^\infty(\mathbb S^{d-1})}
\|g\|_{H^N(\mathbb S^{d-1})}
+
\|g\|_{L^\infty(\mathbb S^{d-1})}
\|f\|_{H^N(\mathbb S^{d-1})}
\Bigr).
\]
Since \(N>(d-1)/2\), the Sobolev embedding
\[
H^N(\mathbb S^{d-1})\hookrightarrow
L^\infty(\mathbb S^{d-1})
\]
holds. Consequently,
\[
\|fg\|_{H^N(\mathbb S^{d-1})}
\leq C_N
\|f\|_{H^N(\mathbb S^{d-1})}
\|g\|_{H^N(\mathbb S^{d-1})},
\]
so \(H^N(\mathbb S^{d-1})\) is a Banach algebra.

\textbf{Substep 2.2} Estimating the inverse of \(\widetilde\Omega\).

By the Sobolev embedding and the assumption on \(\widetilde\Omega\), after possibly decreasing \(\varepsilon_*\), we have
\[
\widetilde\Omega\geq\frac12
\qquad\text{on }\mathbb S^{d-1}.
\]
Set \(h:=\widetilde\Omega-1\). 
By assumption, \(\|h\|_{H^N}<\varepsilon_*\). Decreasing \(\varepsilon_*\) if necessary, we may assume that for some constant $C_N$ we have
\[
C_N\|h\|_{H^N}<C_N\varepsilon_*<1. 
\]
Thus, the Neumann series $\widetilde\Omega^{-1}
=(1+h)^{-1}=\sum_{k=0}^\infty(-h)^k$ converges in \(H^N(\mathbb S^{d-1})\). Thus, we have $\widetilde\Omega^{-1}-1=\sum_{k=1}^\infty(-h)^k.$ By iterating the algebra estimate, we have
\[
\|h^k\|_{H^N}
\leq
C_N^{k-1}\|h\|_{H^N}^k
\qquad\text{for every }k\geq1.
\]
Hence,
\begin{align*}
\|\widetilde\Omega^{-1}-1\|_{H^N}
&\leq
\sum_{k=1}^\infty\|h^k\|_{H^N}
\leq
\|h\|_{H^N}
\sum_{k=1}^\infty
C_N^{k-1}\|h\|_{H^N}^{k-1}\\
&\leq
\|h\|_{H^N}
\sum_{k=1}^\infty
(C_N\varepsilon_*)^{k-1}
=
\frac{1}{1-C_N\varepsilon_*}
\|h\|_{H^N}.
\end{align*}
Hence,
\begin{equation*}
\|\widetilde\Omega^{-1}-1\|_{H^N}
\leq
C\|\widetilde\Omega-1\|_{H^N},
\qquad
C=(1-C_N\varepsilon_*)^{-1}.
\end{equation*}

\textbf{Step 3.} Conclusion.

Since the constant function is a spherical harmonic of degree zero,
\[
T_\kappa^{-1}[1]=\lambda_{\kappa,0}^{-1}>0.
\]
By the linearity of \(T_\gamma^{-1}\), the estimate \eqref{eq:T-kappa-inverse-HN-C0}, and the assumption on \(\widetilde\Omega\), we obtain
\begin{align}
\left\|T_\gamma^{-1}[\widetilde\Omega]
-\lambda_{\gamma,0}^{-1}\right\|_{L^\infty}
&=
\left\|T_\gamma^{-1}[\widetilde\Omega-1]\right\|_{L^\infty} \\
&\leq
C\|\widetilde\Omega-1\|_{H^N}
<
C\varepsilon_*.
\end{align}
Similarly, using \eqref{eq:inverse-Omega-small-perturbation}, we have
\begin{align*}
\left\|T_\delta^{-1}[\widetilde\Omega^{-1}]
-\lambda_{\delta,0}^{-1}\right\|_{L^\infty}
&=
\left\|T_\delta^{-1}[\widetilde\Omega^{-1}-1]\right\|_{L^\infty} \\
&\leq
C\|\widetilde\Omega^{-1}-1\|_{H^N}
\leq
C\|\widetilde\Omega-1\|_{H^N}
<
C\varepsilon_*.
\end{align*}
Therefore, choosing \(\varepsilon_*>0\) sufficiently small so that
\[
C\varepsilon_*
<
\frac12
\min\left\{
\lambda_{\gamma,0}^{-1},
\lambda_{\delta,0}^{-1}
\right\},
\]
the reverse triangle inequality gives
\[
T_\gamma^{-1}[\widetilde\Omega](\theta)
\geq
\lambda_{\gamma,0}^{-1}
-
\left\|T_\gamma^{-1}[\widetilde\Omega]
-\lambda_{\gamma,0}^{-1}\right\|_{L^\infty}
\geq
\frac{1}{2\lambda_{\gamma,0}}>0.
\]
The same argument gives
\[
T_\delta^{-1}[\widetilde\Omega^{-1}]
\geq
\frac{1}{2\lambda_{\delta,0}}>0
\qquad\text{on }\mathbb S^{d-1}.
\]
\end{proof}

\begin{remark}\label{remark_inverse_T}
As mentioned in the previous section, we apply
Lemma \ref{lemma_L_integral} to the operators \(L\) and \(A\). In the notation of
Proposition \ref{prop:small-perturbations}, this corresponds to
\(\gamma=\alpha=s-d+2\) and \(\delta=2-\gamma=d-s\). Accordingly, the
anisotropic kernels \(\omega_L\) and \(\omega_A\) are defined by
\begin{equation}\label{omegaAT}
\omega_L:=\frac{1}{C_{d,d-s}}T_{d-s}^{-1}[\widetilde{\Omega}^{-1}]
=\frac{1}{C_{d,\delta}}T_\delta^{-1}[\widetilde{\Omega}^{-1}],
\qquad
\omega_A:=\frac{1}{C_{d,\alpha}}T_{\alpha}^{-1}[\widetilde{\Omega}].
\end{equation}
Since the constants \(C_{d,d-s}\) and \(C_{d,\alpha}\) are positive, these normalization factors do not affect the positivity assumptions \eqref{hl} and \eqref{ha}.
Moreover, Lemma \ref{lemma_L_integral} shows that \(K_A=\Delta W\), since their renormalized Fourier transforms coincide.

The second equation in \eqref{omegaAT} shows that \eqref{ha} implies $\tilde{\Omega}>0$ because the operator $T_\alpha$ has a positive kernel, see \eqref{tgamma}. In fact, the implication from $\Delta W>0$ to the positive definiteness of $W$ is true for more general interaction potentials $W$, see \cite[Lemma 1.14]{shu2025existence}.
\end{remark}

We next clarify the relations between the positivity assumptions \eqref{hw}, \eqref{hl}, and \eqref{ha}. We first prove that Assumption \eqref{hl} implies Assumption \eqref{hw}. We then consider the two-dimensional setting and provide explicit counterexamples showing that, in general, the remaining implications between these assumptions do not necessarily hold.
\begin{lemma}[Positivity of the inverse kernel]
\label{lem:HL-implies-HW}
Let $d\ge2$ and $d-2<s<d$, and set $\delta=d-s\in(0,2)$. Then \eqref{hl} implies \eqref{hw}.
\end{lemma}

\begin{proof}

Assume \eqref{hl}. We first prove the positivity of $\tilde{\Omega}$. By Lemma \ref{lemma_L_integral}, we can express $\psi$ in terms of $\omega_L$ as
\begin{equation}
    \psi(\xi)
= \mathcal{RF}[|\cdot|^{-d-\delta}\omega_L](\xi) = 
C_{d,\delta}|\xi|^\delta T_\delta[\omega_L]\left(\frac{\xi}{|\xi|}\right),\qquad \xi\ne 0
\end{equation}
where we notice that $C_{d,\delta}>0$ and the integral kernel of $T_\delta$ defined in \eqref{tgamma} is nonnegative. Since $\omega_L$ has a positive minimum and a finite maximum,
there exist constants $c,C>0$ such that
\begin{equation}
\label{eq:HL-symbol-bounds}
c|\xi|^\delta\le\psi(\xi)\le C|\xi|^\delta
\qquad\text{for all }\xi\in\mathbb R^d.
\end{equation}
In particular,
\[
\widetilde\Omega(\theta)=\psi(\theta)^{-1}>0
\qquad\text{for every }\theta\in\mathbb S^{d-1}.
\]

Next we prove the positivity of $\Omega$. We claim that 
\begin{equation}\label{claimpt}
    \text{For every $t>0$, $p_t := \mathcal{F}^{-1}[e^{-t\psi}]$ is a strictly positive function.}
\end{equation}
The positivity of $\Omega$ follows from this claim. Indeed, For $\varphi\in\mathcal S(\mathbb R^d)$, we have
\[
\begin{aligned}
\int_0^\infty\int_{\mathbb R^d}
e^{-t\psi(\xi)}|\widehat\varphi(-\xi)|d\xi dt
&=
\int_{\mathbb R^d}
\frac{|\widehat\varphi(-\xi)|}{\psi(\xi)}d\xi\\
&\le
c^{-1}\int_{\mathbb R^d}
|\xi|^{-\delta}|\widehat\varphi(-\xi)|d\xi
<\infty.
\end{aligned}
\]
The last integral converges near the origin because $\delta<d$,
and at infinity because $\widehat\varphi$ is rapidly decreasing.
Thus Fourier inversion and Fubini's theorem give
\begin{equation}
\label{eq:HL-potential-identity}
\begin{aligned}
\int_0^\infty \langle p_t,\varphi\rangle dt
&=
(2\pi)^{-d}
\int_{\mathbb R^d}
\left(\int_0^\infty e^{-t\psi(\xi)}dt\right)
\widehat\varphi(-\xi)d\xi\\
&=
(2\pi)^{-d}
\int_{\mathbb R^d}
\frac{\widehat\varphi(-\xi)}{\psi(\xi)}d\xi
=
\langle W,\varphi\rangle.
\end{aligned}
\end{equation}
Hence 
\begin{equation}
    W=\int_0^\infty p_tdt
\end{equation}
in the distributional sense. To pass from this distributional identity to pointwise
strict positivity, define
\[
g(x):=\int_1^2 p_t(x)dt.
\]
The bound \eqref{eq:HL-symbol-bounds} implies
\[
e^{-t\psi(\xi)}\le e^{-c|\xi|^\delta},
\qquad 1\le t\le2.
\]
Since the right-hand side is integrable, Fourier
inversion and dominated convergence show that
$(t,x)\mapsto p_t(x)$ is continuous on
$[1,2]\times\mathbb R^d$. Hence $g$ is continuous,
and \eqref{claimpt} gives $g(x)>0$ for every $x$.

For every nonnegative
$\varphi\in C_c^\infty(\mathbb R^d\setminus\{0\})$,
\eqref{eq:HL-potential-identity} and $p_t\ge0$ give
\[
\langle W-g,\varphi\rangle
=
\int_{(0,1)\cup(2,\infty)}
\langle p_t,\varphi\rangle dt
\ge0.
\]
Thus $W\ge g$ distributionally on
$\mathbb R^d\setminus\{0\}$.
Since both functions are continuous there, this
inequality holds pointwise. Therefore
\[
W(x)\ge g(x)>0
\qquad (x\ne0),
\]
and in particular
$\Omega(\theta)=W(\theta)>0$
for every $\theta\in\mathbb S^{d-1}$.

Finally we prove the claim \eqref{claimpt} using the L\'evy--Khintchine theorem from probability theory, following \cite{schilling2016introduction}. We first decompose the kernel
\[
|x|^{-d-\delta}\omega_L(x/|x|)
=
\frac m2|x|^{-d-\delta}+K_1(x),\quad m:=\min_{\mathbb S^{d-1}}\omega_L>0
\]
where
\[
K_1(x)
=
|x|^{-d-\delta}
\left(\omega_L(x/|x|)-\frac m2\right)\ge0.
\]
Consequently, on the Fourier side we have the decomposition
\[
\psi(\xi)=a|\xi|^\delta+\psi_1(\xi),
\qquad a>0,
\]
where $\psi_1 = \mathcal{RF}[K_1]$. Then we write 
\begin{equation}
    e^{-t\psi(\xi)} = e^{-t\psi_1(\xi)}\cdot e^{-t a|\xi|^\delta}
\end{equation}
and analyze the inverse Fourier transform of each factor separately.

In the sense of \cite[Definition 6.10]{schilling2016introduction}, the measure $\nu(dx):=K_1(x)dx$ is a symmetric L\'evy
measure, since $K_1\ge 0$ and
\[
\int_{\mathbb R^d}(1\wedge |x|^2)K_1(x)dx<\infty.
\]
Apply \cite[Theorem 7.5]{schilling2016introduction} to the L\'evy triplet
$(0,0,\nu)$. By symmetry, the characteristic exponent in
\cite[equation (7.1)]{schilling2016introduction} reduces to $\psi_1 = \mathcal{RF}[K_1]$. Hence there exists a L\'evy process $X_t$ satisfying
$\mathbb E e^{i\xi\cdot X_t}=e^{-t\psi_1(\xi)}$.
Its probability law $\nu_t:=\operatorname{Law}(X_t)$ therefore
satisfies, using the evenness of $\psi_1$,
\[
\widehat{\nu_t}(\xi)
=
\mathbb E e^{-i\xi\cdot X_t}
=
e^{-t\psi_1(\xi)}.
\]

We next show that
$q_t:=\mathcal F^{-1}[e^{-ta|\cdot|^\delta}]$
is a strictly positive probability density.
Since $f(\lambda)=a\lambda^{\delta/2}$ is a Bernstein
function, \cite[Theorem 5.2]{SchillingSongVondracek2012}
provides a nonnegative measure $\eta_t$ on $[0,\infty)$
satisfying
\[
\int_{[0,\infty)}e^{-\lambda r}\eta_t(dr)
=e^{-ta\lambda^{\delta/2}}.
\]
Letting $\lambda\downarrow0$ and $\lambda\to\infty$,
respectively, gives
\[
\eta_t([0,\infty))=1,
\qquad \eta_t(\{0\})=0.
\]
Taking $\lambda=|\xi|^2$, we obtain
\[
e^{-ta|\xi|^\delta}
=\int_{(0,\infty)}e^{-r|\xi|^2}\eta_t(dr).
\]
Moreover, Tonelli's theorem yields
\[
\int_{(0,\infty)}
 \int_{\mathbb R^d}e^{-r|\xi|^2}d\xi\eta_t(dr)
=
\int_{\mathbb R^d}e^{-ta|\xi|^\delta}d\xi
<\infty.
\]
Thus Fourier inversion and Fubini's theorem give, for
every $x\in\mathbb R^d$,
\[
q_t(x)
=
\int_{(0,\infty)}
(4\pi r)^{-d/2}
\exp\left(-\frac{|x|^2}{4r}\right)\eta_t(dr)>0.
\]
The strict inequality follows because the Gaussian
integrand is strictly positive and
$\eta_t((0,\infty))=1$.
Integrating this identity in $x$ and applying Tonelli's
theorem shows that $\int_{\mathbb R^d}q_t(x)dx=1$.
Also, since $e^{-ta|\cdot|^\delta}\in L^1(\mathbb R^d)$,
the function $q_t$ is bounded and continuous.

Consequently,
\[
p_t(x)
=(q_t*\nu_t)(x)
=\int_{\mathbb R^d}q_t(x-y)\nu_t(dy)>0
\qquad\text{for every }x\in\mathbb R^d.
\]
Indeed, the convolution is a bounded continuous
probability density with Fourier transform
$e^{-t\psi}$, and it is strictly positive because
$q_t$ is strictly positive and $\nu_t$ is a probability
measure. This proves \eqref{claimpt} and completes the proof.
\end{proof}

We now turn to counterexamples showing that the remaining implications between the positivity assumptions fail in general. For this purpose, we restrict ourselves to the two-dimensional setting, where explicit examples can be constructed. Let \(d=2\). Then \(\alpha=s\), and writing \(W(r,\theta)=r^{-s}\Omega(\theta)\), we have
\[
\Delta W(r,\theta)
=
r^{-s-2}\bigl(\Omega''(\theta)+s^2\Omega(\theta)\bigr).
\]
Since \(K_A=\Delta W\) by Remark \ref{remark_inverse_T}, it follows that
\begin{equation}\label{eq:omega-A-two-dimensional}
\omega_A(\theta)
=
\Omega''(\theta)+s^2\Omega(\theta).
\end{equation}

\begin{example}[Assumption \eqref{hw} does not imply \eqref{ha} {or \eqref{hl}} in general]
Let \(d=2\) and \(0<s<2\). We consider the fourth angular mode, which is convenient since \(\partial_{\theta\theta}\cos(4\theta)=-16\cos(4\theta)\), while the Fourier transform preserves the same angular mode up to a multiplicative factor. Set
\begin{equation}\label{eq:Omega-cos4-counterexample}
\Omega(\theta):=1+\frac{s}{4}\cos(4\theta).
\end{equation}
Since \(0<s<2\), we have \(\Omega(\theta)\geq1-s/4>1/2\), and hence \(\Omega>0\).

We next verify the positivity of \(\widetilde\Omega\). Writing \(\xi=\rho(\cos\varphi,\sin\varphi)\), by the linearity of the Fourier transform and the formulas in \cite[Lemma 2.1 and Sections 9.1--9.2]{CS24cpam}, we have
\[
\mathcal F\bigl[|x|^{-s}\bigr](\xi)=c_s|\xi|^{s-2},\qquad
\mathcal F\bigl[|x|^{-s}\cos(4\theta)\bigr](\xi)
=
c_s|\xi|^{s-2}\frac{(4-s)(2-s)}{s(s+2)}\cos(4\varphi),
\qquad c_s>0.
\]
Indeed, the second identity follows from the explicit computations for the modes \(\cos^2\theta\) and \(\cos^4\theta\), together with \(\cos(4\theta)=8\cos^4\theta-8\cos^2\theta+1\). Therefore, using \eqref{eq:Omega-cos4-counterexample}, we obtain
\[
\widehat W(\xi)
=
c_s|\xi|^{s-2}\left(1+\frac{(4-s)(2-s)}{4(s+2)}\cos(4\varphi)\right).
\]
Since \(\widehat W(\xi)=|\xi|^{s-2}\widetilde\Omega(\varphi)\), we conclude
\begin{equation}\label{eq:tilde-Omega-cos4-counterexample}
\widetilde\Omega(\varphi)
=
c_s\left(1+\frac{(4-s)(2-s)}{4(s+2)}\cos(4\varphi)\right).
\end{equation}
Moreover,
\[
1-\frac{(4-s)(2-s)}{4(s+2)}
=
\frac{s(10-s)}{4(s+2)}>0,
\]
then \eqref{eq:tilde-Omega-cos4-counterexample} yields \(\widetilde\Omega>0\). Hence Assumption \eqref{hw} holds. On the other hand, by \eqref{eq:omega-A-two-dimensional} and \eqref{eq:Omega-cos4-counterexample}, we have that
\[
\omega_A(0)
=
\frac{s}{4}\bigl(s^2+4s-16\bigr)<0,
\]
where we use that \(0<s<2\). Thus \(\omega_A\) is not positive on \(\mathbb S^1\), and Assumption \eqref{ha} fails. Therefore, we conclude
\[\eqref{hw}\not\Longrightarrow\eqref{ha}\qquad\text{for every }0<s<2.
\]

{Next, we focus on the case $s=1$, for which we have verified \eqref{hw} for the example \eqref{eq:Omega-cos4-counterexample}. In this case we have $\tilde{\Omega}(\varphi)=c_1 (1+\frac{1}{4}\cos(4\varphi))$, and thus by Remark \ref{remark_inverse_T}, $\omega_L$ is a positive constant multiple of $T_1^{-1}[\tilde{g}]$, where $\tilde{g}(\varphi)=1/(1+\frac{1}{4}\cos(4\varphi))$. Applying \cite{CS24cpam} in the case $s=1$, we see that its corresponding angular profile in the physical space is $g(\theta)=\tilde{g}(\theta+\frac{\pi}{2}) = 1/(1+\frac{1}{4}\cos(4\theta))$ up to a positive constant multiple. Remark \ref{remark_inverse_T} with 
 \eqref{eq:omega-A-two-dimensional} show that $T_1^{-1}[\tilde{g}]$ can be expressed in terms of $g$, up to a positive constant multiple, as $g''(\theta)+g(\theta)$. We verify that
\[
g''\big(\frac{\pi}{4}\big)+g\big(\frac{\pi}{4}\big) = -\frac{52}{9}<0,
\]
which implies that \eqref{hl} fails. Therefore we conclude that 
\[\eqref{hw}\not\Longrightarrow\eqref{hl}\qquad\text{for }s=1.
\]
}
\end{example}

\begin{example}[Assumption \eqref{ha} does not imply \eqref{hw} in general]
Let \(d=2\), \(\sqrt{2}<s<2\), and consider
\[
\Omega(\theta)=1+\varepsilon\cos(2\theta), \qquad \varepsilon>1.
\]
Then \(\Omega(\pi/2)=1-\varepsilon<0\), so the physical-space positivity condition in \eqref{hw} fails. On the other hand, by \eqref{eq:omega-A-two-dimensional}, we have
\[
\omega_A(\theta)
=\Omega''(\theta)+s^2\Omega(\theta)
=s^2+\varepsilon(s^2-4)\cos(2\theta).
\]
Since \(\sqrt{2}<s<2\) and taking $1 < \varepsilon < \frac{s^2}{4-s^2}$, we obtain
\[
\omega_A(\theta)
\geq s^2-\varepsilon(4-s^2) > 0.
\]
Therefore, we conclude
\[
\eqref{ha}\not\Longrightarrow\eqref{hw}.
\]
\end{example}

\begin{remark}
The previous example also shows that Assumption \eqref{ha} does not imply Assumption \eqref{hl} in general. Indeed, since $\eqref{hl}\Longrightarrow\eqref{hw},$ the implication \(\eqref{ha}\Longrightarrow\eqref{hl}\) would imply \(\eqref{ha}\Longrightarrow\eqref{hw}\), contradicting the previous example. Hence,
\[
\eqref{ha}\not\Longrightarrow\eqref{hl}
\qquad\text{in general}.
\]
\end{remark}

\begin{example}[Assumption \eqref{hl} does not imply \eqref{ha} in general]
Let \(d=2\), \(s=1\), and consider
\[
\Omega(\theta):=\frac{1}{\sqrt{\cos^2\theta+3\sin^2\theta}}.
\]
Then, setting \(x=r(\cos\theta,\sin\theta)\), we have
\[
W(r,\theta)
=
r^{-1}\Omega(\theta)
=
\frac{1}{r\sqrt{\cos^2\theta+3\sin^2\theta}}
=
\frac{1}{\sqrt{x_1^2+3x_2^2}}.
\]
Now, using the change of variables \(y_1=x_1\), \(y_2=\sqrt{3}x_2\), we obtain \(x_1^2+3x_2^2=|y|^2\), and hence
\[
\widehat W(\xi)
=
\frac{1}{\sqrt{3}}
\mathcal F\bigl(|y|^{-1}\bigr)
\left(\xi_1,\frac{\xi_2}{\sqrt{3}}\right)
=
c\left(\xi_1^2+\frac{\xi_2^2}{3}\right)^{-1/2},
\qquad c>0,
\]
where we used the classical Riesz-potential identity
\[
\mathcal F\bigl(|y|^{-1}\bigr)(\eta)
=
c|\eta|^{-1},
\qquad c>0.
\]
Hence, from \eqref{defi_L_multiplier_main}, the symbol of \(L\) is
\[
\psi(\xi)
=
c^{-1}\left(\xi_1^2+\frac{\xi_2^2}{3}\right)^{1/2}.
\]
Thus, its singular kernel is, up to a positive constant, given by
\[
K_L(h)
=
c_L\left(h_1^2+3h_2^2\right)^{-3/2},
\qquad c_L>0.
\]
Since \(K_L>0\), Assumption \eqref{hl} holds. On the other hand,  since  \(\Omega(0)=1\) and \(\Omega''(0)=-2\), it follows from \eqref{eq:omega-A-two-dimensional} that \(\omega_A(0)<0\). Hence, assumption \eqref{ha} fails.
\end{example}

\section{Proof of Theorem \ref{thm:uniqueness}}
\label{sec:proof-uniqueness}

This section is devoted to the proof of Theorem \ref{thm:uniqueness}, namely
the uniqueness of steady states among bounded continuous densities. Throughout
the section, we work with the quadratic confinement \(U(x)=|x|^2\). We first
recall the notation for measures with fixed mass. We denote 
\[
\mathcal P(\mathbb R^d)
:=
\left\{
\mu\in\mathcal M_+(\mathbb R^d)
\;:\;
\mu(\mathbb R^d)=1
\right\}.
\]
For \(m>0\), we set
\[
\mathcal P_{2,m}(\mathbb R^d)
:=
\left\{
\mu\in\mathcal M_+(\mathbb R^d):
\mu(\mathbb R^d)=m,\ 
\int_{\mathbb R^d}|x|^2d\mu(x)<\infty
\right\}.
\]
In particular, \(\mathcal P_{2,1}(\mathbb R^d)=\mathcal P_2(\mathbb R^d)\).
When \(\mu\in\mathcal P(\mathbb R^d)\) is absolutely continuous with respect to
Lebesgue measure, we write \(d\mu(x)=\rho(x)dx\) and identify \(\mu\) with its
density \(\rho\).

We now prove the maximum principle for \(L\). By Lemma \ref{lemma_L_integral},
\(L\) admits a representation as a difference-quotient integral. Thus,
if \(\omega_L>0\), \(L\) satisfies a monotonicity property: when two ordered
functions touch at a point, their values at that point satisfy the reverse
inequality.

\begin{lemma}[Maximum principle]\label{lemma_maximum_principle}
Assume that $\eqref{hl}$ holds. Let
\(p_0,p_1\in\mathcal S(\mathbb R^d)\) satisfy \(p_0\le p_1\) in
\(\mathbb R^d\), and suppose that there exists \(x_0\in\mathbb R^d\) such that
\(p_0(x_0)=p_1(x_0)\). Then
\begin{equation}\label{eq:mp}
L[p_0](x_0)\ge L[p_1](x_0).
\end{equation}
Moreover, equality holds if and only if \(p_0\equiv p_1\) in \(\mathbb R^d\).
\end{lemma}

\begin{proof}
By \eqref{eq:L-kernel}, for \(i=0,1\),
\[
L[p_i](x_0)
=
\mathrm{p.v.}
\int_{\mathbb R^d}
\big(p_i(x_0)-p_i(y)\big)
|x_0-y|^{s-2d}
\omega_L\left(\frac{x_0-y}{|x_0-y|}\right)dy.
\]
Subtracting the identities for \(i=0\) and \(i=1\), and using
\(p_0(x_0)=p_1(x_0)\), we obtain
\[
L[p_0](x_0)-L[p_1](x_0)
=
\mathrm{p.v.}
\int_{\mathbb R^d}
\big(p_1(y)-p_0(y)\big)
|x_0-y|^{s-2d}
\omega_L\left(\frac{x_0-y}{|x_0-y|}\right)dy.
\]
Since \(p_1-p_0\ge0\) and \(\omega_L>0\), the integrand is nonnegative. Hence
\(L[p_0](x_0)-L[p_1](x_0)\ge0\), which proves \eqref{eq:mp}. Now, if equality holds, then the nonnegative integrand above must vanish a.e. in
\(\mathbb R^d\). Since the kernel is strictly positive away from \(x_0\), we get
\(p_1=p_0\) a.e. in \(\mathbb R^d\). By continuity, \(p_0\equiv p_1\).
\end{proof}

\begin{corollary}\label{corollary_principle_maximum}
Assume that $\eqref{hl}$ holds. Let
\(f\in\mathcal S(\mathbb R^d)\) have a fixed sign in \(\mathbb R^d\). If there
exists \(x_*\in\mathbb R^d\) such that $f(x_*)=0$ and $L[f](x_*)=0$, then \(f\equiv0\) in \(\mathbb R^d\).
\end{corollary}

\begin{proof}
It is enough to consider the case \(f\ge0\). 
Since \(f(x_*)=0\), the representation
\eqref{eq:L-kernel-main} gives
\[
L[f](x_*)
=
-\mathrm{p.v.}
\int_{\mathbb R^d}f(y)K_L(x_*-y)dy.
\]
By $\eqref{hl}$, \(K_L\) is strictly positive away from the
origin. Hence the integrand \(f(y)K_L(x_*-y)\) is nonnegative, and the identity
\(L[f](x_*)=0\) implies
\[
\int_{|y-x_*|>\varepsilon} f(y)K_L(x_*-y)dy=0
\qquad\text{for every }\varepsilon>0.
\]
Therefore \(f(y)K_L(x_*-y)=0\) for a.e. \(y\neq x_*\), and so \(f=0\) a.e. in
\(\mathbb R^d\). Since \(f\) is continuous, \(f\equiv0\) in \(\mathbb R^d\).
\end{proof}

\begin{remark}
Although Lemma \ref{lemma_maximum_principle} and Corollary \ref{corollary_principle_maximum}
are stated for Schwartz functions, the same conclusions apply to the potentials used below.
Indeed, if \(q\in C(\mathbb R^d)\cap L^\infty(\mathbb R^d)\) satisfies
\[
Lq=g
\]
in the sense of distributions, with \(g\in C(\mathbb R^d)\), then \(q\) is also a viscosity
solution of this equation. Hence, if \(q\) attains a global minimum at \(x_0\), the maximum
principle gives \(g(x_0)\leq 0\).

Moreover, if
\[
q\geq 0,\qquad q(x_0)=0,\qquad g(x_0)=0,
\]
then the strict positivity of \(K_L\) away from the origin yields the strong maximum
principle \(q\equiv 0\). We apply this below to differences of potentials of the form
\(W*f\), for which \(L(W*f)=f\) distributionally. Thus Lemma \ref{lemma_maximum_principle} and Corollary \ref{corollary_principle_maximum} remain applicable even though these
potentials are not Schwartz functions.
\end{remark}

\begin{lemma}\label{lem:farfield-decay}
Let \(\rho\in \mathcal{P}(\mathbb R^d)\cap L^\infty(\mathbb R^d)\) and
\(p:=W*\rho\). Then \(p(x)\to0\) as \(|x|\to\infty\).
\end{lemma}

\begin{proof}
Fix \(\varepsilon>0\).
Since \(s<d\), we have \(\int_{|z|<L}|z|^{-s}dz
=
\frac{|\mathbb S^{d-1}|}{d-s}L^{d-s}.
\) Choose \(L>0\) such that $\|\Omega\|_{L^\infty}\|\rho\|_{L^\infty}
\frac{|\mathbb S^{d-1}|}{d-s}L^{d-s}
<\varepsilon,$ it follows that
\[
|p(x)|
\le
\|\Omega\|_{L^\infty}
\int_{|x-y|<L}|x-y|^{-s}\rho(y)dy
+
\|\Omega\|_{L^\infty}
\int_{|x-y|\ge L}|x-y|^{-s}\rho(y)dy.
\]
The first term is bounded by \(\varepsilon\). For the second term, choose \(R>0\) so large that \(\int_{|y|>R}\rho(y)dy<\varepsilon L^s.\) Then, we obtain
\[
\int_{|x-y|\ge L}|x-y|^{-s}\rho(y)dy
\le
\int_{\substack{|y|\le R\\ |x-y|\ge L}}
|x-y|^{-s}\rho(y)dy
+
L^{-s}\int_{|y|>R}\rho(y)dy.
\]
For fixed \(R\), the first term tends to \(0\) as \(|x|\to\infty\), by
dominated convergence. The second term is bounded by \(\varepsilon\). Hence
\[
\limsup_{|x|\to\infty}|p(x)|
\le
(1+\|\Omega\|_{L^\infty})\varepsilon.
\]
Since \(\varepsilon>0\) is arbitrary, \(p(x)\to0\) as \(|x|\to\infty\).
\end{proof}

The next lemma derives a geometric consequence of the Euler--Lagrange condition.
Since the effective potential is constant on the support of the density and lies
above this value elsewhere, the maximum principle for \(L\) implies a convexity
property for $V[\rho] = W*\rho+|\cdot|^2$ outside the support. This convexity will be used in the sliding argument in the
proof of Theorem \ref{thm:uniqueness}, where one has to control possible contact points between competing potentials.

\begin{lemma}[Euler--Lagrange condition]\label{lemma_convex}
Assume that $\eqref{hw}$ and $\eqref{hl} $ hold. Let \(\rho \in C_c(\mathbb{R}^d)\) be a nonnegative density with fixed positive total mass satisfying the Euler--Lagrange equation for some constant \(C_0 \in \mathbb{R}\), so that
\begin{equation}\label{eq:EL_rho}
V(x) = C_0 \quad \text{for } x \in \operatorname{supp}\rho, \qquad \text{and} \qquad V(x) \geq C_0 \quad \text{for } x \in \mathbb{R}^d  \setminus \supp \rho.
\end{equation}
Then, \( V[\rho] = P + |\cdot|^2 \) is convex on $\mathbb{R}^d$. Moreover, denoting \(D := \mathbb{R}^d \setminus \operatorname{supp}\rho\), we have
\begin{equation*} 
V(x) > C_0 \quad \textrm{and} \quad \nabla V(x) \neq 0 \qquad \text{for every } x \in D,
\end{equation*}
where \(P := W * \rho\).
\end{lemma}

\begin{proof}
For simplicity, write \(V(x):=V[\rho](x)\). Fix a unit vector
\(e\in\mathbb R^d\) and \(h>0\). We define the second difference
\[
\Delta^2_{h,e}V(x):=V(x+he)+V(x-he)-2V(x).
\]
Since \(U(x)=|x|^2\), we have
\begin{equation}\label{eq:convex_U_second_diff}
\Delta^2_{h,e}U(x)
=
|x+he|^2+|x-he|^2-2|x|^2
=
2h^2.
\end{equation}
Moreover, for \(x\in D\), since \(\rho(x)=0\), we have
\begin{equation}\label{eq:Lu_nonneg_on_D}
L[\Delta^2_{h,e}V](x)
=
\rho(x+he)+\rho(x-he)-2\rho(x)
=
\rho(x+he)+\rho(x-he)
\ge0.
\end{equation}

{\smallskip
\noindent\emph{\textbf{Step 1:}}}
We now prove that \(\Delta^2_{h,e}V\ge0\) in \(\mathbb R^d\). Suppose, by
contradiction, that
\[
\inf_{\mathbb R^d}\Delta^2_{h,e}V<0.
\]
By Lemma \ref{lem:farfield-decay}, \(P(x)\to0\) as \(|x|\to\infty\), which
implies
\[
\Delta^2_{h,e}P(x)
=
P(x+he)+P(x-he)-2P(x)
\longrightarrow0
\qquad\text{as } |x|\to\infty.
\]
Combining this with \eqref{eq:convex_U_second_diff}, we obtain
\begin{equation}\label{eq:u_limit_infty}
\Delta^2_{h,e}V(x)
=
2h^2+\Delta^2_{h,e}P(x)
\longrightarrow 2h^2>0
\qquad\text{as }|x|\to\infty.
\end{equation}
Therefore, the negative infimum of \(\Delta^2_{h,e}V\) is attained at some point
\(x_*\in\mathbb R^d\), namely
\[
\Delta^2_{h,e}V(x_*)
=
\inf_{\mathbb R^d}\Delta^2_{h,e}V
<0.
\]

We first show that \(x_*\in D\). Indeed, if
\(x_*\in\operatorname{supp}\rho\), then \eqref{eq:EL_rho} implies
\(V(x_*)=C_0\), while \(V\ge C_0\) in \(\mathbb R^d\). Hence
\[
\Delta^2_{h,e}V(x_*)
=
V(x_*+he)+V(x_*-he)-2V(x_*)
\ge
C_0+C_0-2C_0
=
0,
\]
which contradicts \(\Delta^2_{h,e}V(x_*)<0\). Therefore
\(x_*\notin\operatorname{supp}\rho\), and hence \(x_*\in D\).

Since \(x_*\in D\), \eqref{eq:Lu_nonneg_on_D} gives
\[
L[\Delta^2_{h,e}V](x_*)\ge0.
\]
On the other hand, since \(x_*\) is a global minimum point of
\(\Delta^2_{h,e}V\), the singular-integral representation of \(L\) gives
\[
L[\Delta^2_{h,e}V](x_*)
=
\mathrm{p.v.}
\int_{\mathbb R^d}
\big(\Delta^2_{h,e}V(x_*)-\Delta^2_{h,e}V(y)\big)
K_L(x_*-y)dy
\le0.
\]
Combining both inequalities, we obtain
$
L[\Delta^2_{h,e}V](x_*)=0.
$
Set
$
\Phi(x):=\Delta^2_{h,e}V(x)-\Delta^2_{h,e}V(x_*).
$
Then \(\Phi\ge0\) in \(\mathbb R^d\), \(\Phi(x_*)=0\), and, since \(L\)
annihilates constants,
$
L[\Phi](x_*)=0.
$
By Corollary \ref{corollary_principle_maximum}, we conclude that
\(\Phi\equiv0\). Hence \(\Delta^2_{h,e}V\) is constant in \(\mathbb R^d\).
Passing to the limit \(|x|\to\infty\) in \eqref{eq:u_limit_infty}, this
constant must be \(2h^2>0\), contradicting
\(\Delta^2_{h,e}V(x_*)<0\).

As a consequence, we have shown that
\begin{equation}\label{eq:second-difference-nonnegative-global}
\Delta^2_{h,e}V(x)\ge0
\qquad
\text{for every }x\in\mathbb R^d,\ h>0,\ e\in\mathbb S^{d-1}.
\end{equation}
Since \(V\) is continuous, \eqref{eq:second-difference-nonnegative-global}
implies that \(V\) is convex. In particular, \(V\) is convex on every convex
subset of \(D\). This proves \eqref{eq:second-difference-nonnegative-global}.

\smallskip
\noindent\emph{\textbf{Step 2:}}
We now prove that
\begin{equation}\label{eq:second-difference-positive-D}
\Delta^2_{h,e}V(x)>0
\qquad\text{for every }x\in D.
\end{equation}
Assume, by contradiction, that there exists \(x_0\in D\) such that
\(\Delta^2_{h,e}V(x_0)=0\). By
\eqref{eq:second-difference-nonnegative-global}, the function
\(\Delta^2_{h,e}V\) is nonnegative in \(\mathbb R^d\) and vanishes at \(x_0\).
Moreover, since \(x_0\in D\), \eqref{eq:Lu_nonneg_on_D} gives
\(L[\Delta^2_{h,e}V](x_0)\ge0\). On the other hand, the representation
\eqref{eq:L-kernel-main} and the nonnegativity of \(\Delta^2_{h,e}V\) imply
\(L[\Delta^2_{h,e}V](x_0)\le0\). Hence
\(L[\Delta^2_{h,e}V](x_0)=0\). From
Corollary \ref{corollary_principle_maximum}, we get
\(\Delta^2_{h,e}V\equiv0\), contradicting \eqref{eq:u_limit_infty}. Therefore
\eqref{eq:second-difference-positive-D} holds.

\smallskip
\noindent\emph{\textbf{Step 3:}}
We claim that
\begin{equation}\label{eq:V-strictly-above-C0-on-D}
V(x)>C_0
\qquad\text{for every }x\in D.
\end{equation}
Suppose, by contradiction, that there exists \(x_0\in D\) such that
\(V(x_0)=C_0\). Let \(x_\partial\in\partial D\) satisfy
\(|x_0-x_\partial|=\operatorname{dist}(x_0,\partial D)\). Then
\(V(x_\partial)=C_0\), and the segment joining \(x_0\) to \(x_\partial\),
except for the endpoint \(x_\partial\), is contained in \(D\). Define
\(g(t):=V((1-t)x_0+t x_\partial)\) for \(t\in[0,1]\). Since \(V\geq C_0\),
we have \(g(t)\geq C_0\), while \(g(0)=g(1)=C_0\). On the other hand, by
convexity,
\(g(t)\leq(1-t)g(0)+tg(1)=C_0\). Hence, \(g\equiv C_0\) on \([0,1]\).

Now fix \(t_0\in(0,1)\) and set
\(x_{t_0}:=(1-t_0)x_0+t_0x_\partial\),
\(r:=|x_\partial-x_0|\), and \(e:=(x_\partial-x_0)/r\). For \(h>0\)
sufficiently small, the points \(x_{t_0}\pm he\) remain on the segment
between \(x_0\) and \(x_\partial\), and hence belong to \(D\). Since
\(V\equiv C_0\) along this segment, we obtain
\[
\Delta^2_{h,e}V(x_{t_0})
=
V(x_{t_0}+he)+V(x_{t_0}-he)-2V(x_{t_0})
=
0,
\]
which contradicts \eqref{eq:second-difference-positive-D}. Therefore,
\eqref{eq:V-strictly-above-C0-on-D} holds.

\smallskip
\noindent\emph{\textbf{Step 4:}}
Finally, we prove that \(\nabla V(x)\neq0\) for every \(x\in D\). Since
\(P=W*\rho\) and \(\operatorname{dist}(x,\operatorname{supp}\rho)>0\) for
\(x\in D\), the kernel \(W(x-\cdot)\) is smooth on
\(\operatorname{supp}\rho\). Hence \(P\), and therefore
\(V=P+|\cdot|^2\), is smooth in \(D\).

Fix \(x\in D\), and choose \(x_\partial\in\partial D\) such that
\(|x-x_\partial|=\operatorname{dist}(x,\partial D)\). Consider \(g\) as in
\emph{Step 3}, with \(x_0\) replaced by \(x\). By
\eqref{eq:V-strictly-above-C0-on-D}, \(g(0)=V(x)>C_0\), while
\(g(1)=V(x_\partial)=C_0\). Since \(g\) is convex,
\[
g'_+(0)\leq g(1)-g(0)<0.
\]
By the smoothness of \(V\) in \(D\),
\[
g'_+(0)
=
\nabla V(x)\cdot(x_\partial-x).
\]
Therefore,
\[
\nabla V(x)\cdot(x_\partial-x)<0,
\]
which implies that \(\nabla V(x)\neq0\). Since \(x\in D\) was arbitrary, the
proof is complete.
\end{proof}

\begin{lemma}\label{lem:rho_m}
Assume $\eqref{hw}$ holds. Suppose that \(\overline{\rho}\in\mathcal P(\mathbb R^d)\) minimizes \(\mathcal E\) with Euler--Lagrange level \(C_0\). For \(m>0\), define $a:=m^{-\frac{1}{s+2}}$ and $\overline{\rho}_m(x):=m a^d\overline{\rho}(ax).$ Then \(\overline{\rho}_m\in\mathcal P_{2,m}(\mathbb R^d)\) minimizes \(\mathcal E\) over \(\mathcal P_{2,m}(\mathbb R^d)\). Moreover, the constant of the Euler-Lagrange is given by \(C_{0,m}=m^{\frac{2}{s+2}}C_0.\)
\end{lemma}

\begin{proof}
{For simplicity of notation, we present the proof for measures admitting densities.
The general case follows along the same lines, with rescalings interpreted as pushforwards under dilations.}

Clearly, \(\overline{\rho}_m\in\mathcal P_{2,m}(\mathbb R^d)\). Since \(W\) is homogeneous of degree \(-s\), we have
\[
(W*\overline{\rho}_m)(x) = m a^d\int_{\mathbb R^d}W(x-y)\overline{\rho}(ay)dy = m a^s(W*\overline{\rho})(ax).
\]
Since \(a^{s+2}=m^{-1}\), we have \(m a^s=a^{-2}=m^{2/(s+2)}\), which implies that
\[
V[\overline{\rho}_m](x) = (W*\overline{\rho}_m)(x)+|x|^2 = m^{\frac{2}{s+2}}V[\overline{\rho}](ax).
\]
The Euler--Lagrange condition for \(\overline{\rho}\) gives \(V[\overline{\rho}]\ge C_0\) on \(\mathbb R^d\) and \(V[\overline{\rho}]=C_0\) on \(\operatorname{supp}\overline{\rho}\). Hence, \(V[\overline{\rho}_m]\ge m^{2/(s+2)}C_0\) on \(\mathbb R^d\), with equality on \(\operatorname{supp}\overline{\rho}_m\), showing that \(C_{0,m}=m^{2/(s+2)}C_0\).

To prove minimality, given \(\rho\in\mathcal P_{2,m}(\mathbb R^d)\), we set \(\eta(x):=\frac{1}{m a^d}\rho\left(\frac{x}{a}\right),\) so that \(\rho(x)=m a^d\eta(ax)\) and \(\eta\in\mathcal P(\mathbb R^d)\). The scaling identities yield
\[
\frac12\int_{\mathbb R^d}(W*\rho)\rho dx = m^2a^s \frac12\int_{\mathbb R^d}(W*\eta)\eta dx \qquad \text{and} \qquad \int_{\mathbb R^d}|x|^2\rho dx = m a^{-2}\int_{\mathbb R^d}|x|^2\eta dx.
\]
Since \(m^2a^s=ma^{-2}=m^{(s+4)/(s+2)}\), it follows that \(\mathcal E[\rho] = m^{\frac{s+4}{s+2}}\mathcal E[\eta]\). By the minimality of \(\overline{\rho}\) in \(\mathcal P(\mathbb R^d)\), we conclude that
\[
\mathcal E[\rho] = m^{\frac{s+4}{s+2}}\mathcal E[\eta] \ge m^{\frac{s+4}{s+2}}\mathcal E[\overline{\rho}] = \mathcal E[\overline{\rho}_m].
\]
Thus, \(\overline{\rho}_m\) minimizes \(\mathcal E\) over \(\mathcal P_{2,m}(\mathbb R^d)\).
\end{proof}

We briefly describe the proof of Theorem \ref{thm:uniqueness}. The argument is divided into four steps. In Step 1, we collect preliminary facts about the steady-state density \(\mu\), its interaction potential \(p:=W*\mu\), and the connected components of its positivity set. We identify the component containing the origin and select a reference minimizer \(\overline{\rho}_{\lambda_*}\) whose Euler--Lagrange level agrees with the level of \(p+U\) on that component. In Step 2, we perform a sliding argument from above and prove that \(P_{\lambda_*}\geq p\) in \(\mathbb R^d\), where \(P_{\lambda_*}:=W*\overline{\rho}_{\lambda_*}\). In Step 3, we show that \(p=P_{\lambda_*}\) on \(\operatorname{supp}\overline{\rho}_{\lambda_*}\). The final step uses a sliding argument from below to obtain \(p\equiv P_{\lambda_*}\) in \(\mathbb R^d\). Applying \(L\) then gives \(\mu=\overline{\rho}_{\lambda_*}\), and comparison of the total masses yields \(\lambda_*=1\). Therefore, \(\mu=\overline{\rho}\).

\begin{proof}[\textbf{Proof of Theorem \ref{thm:uniqueness}}]

\vspace{0.2cm}
\noindent
\begin{step}[\textit{Preliminaries}]\label{step1}{\rm

Let \(\overline{\rho}\) denote the unique minimizer of the energy
\(\mathcal E\) among probability measures in \(\mathcal P(\mathbb R^d)\). Writing $\overline P:=W*\overline{\rho},$
the Euler--Lagrange conditions imply that there exists a constant \(C_0>0\)
such that
\begin{equation}\label{eq:EL-reference-unit-mass}
\overline P(x)+U(x)=C_0
\quad x\in\operatorname{supp}\overline{\rho},
\qquad\text{and}\qquad
\overline P(x)+U(x)\geq C_0
\quad x\in\mathbb R^d.
\end{equation}

Let \(\mu\in C(\mathbb R^d)\cap L^\infty(\mathbb R^d)\) be a nonnegative steady state, and set \(p:=W*\mu\geq0\). Since \(W\) is locally integrable, \(\mu\) is bounded, and \(\mu\) unit mass, the Dominated Convergence Theorem implies that \(p\) is continuous. Hence, the effective potential \(V[\mu]:=p+U\) is also continuous.

For each \(x\in\{\mu>0\}\), let \(\Gamma_x\) denote the connected component of \(\{\mu>0\}\) containing \(x\). {Since \(\mu\) is a steady state, \eqref{eq:defV} implies that
\(p+U\in H^1_{\mathrm{loc}}(\{\mu>0\})\) and its weak gradient
vanishes almost everywhere there. Hence \(p+U\) is constant
almost everywhere on each connected component of \(\{\mu>0\}\).
By continuity, \(p+U\) is constant on each connected component of \(\{\mu>0\}\).}
Hence, for every \(x\in\{\mu>0\}\), there exists \(C_{\Gamma_x}\in\mathbb R\) such that
\begin{equation}\label{condition_steady_state_mu}
p(y)+U(y)=C_{\Gamma_x}
\qquad\text{for }y\in\Gamma_x.
\end{equation}
By continuity, the same identity holds on \(\overline{\Gamma_x}\). The constant \(C_{\Gamma_x}\) depends only on the connected component. In particular, \(C_{\Gamma_x}=C_{\Gamma_y}\) whenever \(\Gamma_x=\Gamma_y\).

From Lemma \ref{lem:farfield-decay}, \(p(x)\to0\) as \(|x|\to\infty\). Since \(p\geq0\) and \(p\not\equiv0\), it attains a positive global maximum. We claim that no global maximum point of \(p\) can belong to \(\{\mu=0\}\). Indeed, if \(x_0\in\{\mu=0\}\) were such a point, since \(p\) is not constant and the kernel is strictly positive under the assumption \eqref{hl}, it follows by the singular-integral representation of \(L\) that \(L[p](x_0)>0\), contradicting \(L[p](x_0)=\mu(x_0)=0\).

Thus, \(p\) attains its global maximum at some \(\bar x\in\{\mu>0\}\). Let \(\Gamma_{\bar x}\) be the connected component containing \(\bar x\), and let \(C_{\Gamma_{\bar x}}\) be its corresponding constant in \eqref{condition_steady_state_mu}. Since \(\mu(\bar x)>0\), there exists \(r>0\) such that \(B_r(\bar x)\subset\Gamma_{\bar x}\). Using \(U(x)=|x|^2\), we obtain
\begin{equation}\label{ident_ball}
p(x)=C_{\Gamma_{\bar{x}}}-|x|^2,
\qquad x\in B_r(\bar x).
\end{equation}
Fix \(e\in\mathbb S^{d-1}\) and \(0<h<r\). Since \(\bar x\) is a global maximizer of \(p\), \eqref{ident_ball} gives \(C_{\Gamma_{\bar{x}}}-|\bar x|^2\geq C_{\Gamma_{\bar{x}}}-|\bar x\pm he|^2\), and hence \(|\bar x|^2\leq|\bar x\pm he|^2\). Expanding the two inequalities yields
\[
-h\leq2\bar x\cdot e\leq h.
\]
Taking \(h\to0\), we obtain \(\bar x\cdot e=0\). Since \(e\) is arbitrary, \(\bar x=0\). Therefore, \(p\) attains its maximum at the origin and \(\mu(0)>0\). We denote by $\Gamma_0$ the connected component of \(\{\mu>0\}\) containing the origin. Then, we have
\begin{equation}\label{eq:central-component-level}
p(x)+U(x)=C_{\Gamma_0}
\quad\text{for }x\in\Gamma_0,
\qquad
C_{\Gamma_0}=p(0).
\end{equation}

We now introduce the family of reference minimizers. For every \(\lambda>0\), define
\[
\overline{\rho}_\lambda(x):=\lambda a_\lambda^d\overline{\rho}(a_\lambda x),
\qquad
a_\lambda:=\lambda^{-\frac1{s+2}},
\]
and set \(P_\lambda:=W*\overline{\rho}_\lambda\). By Lemma \ref{lem:rho_m}, \(\overline{\rho}_\lambda\) is the minimizer of mass \(\lambda\), and its Euler--Lagrange level is given by $C_\lambda=\lambda^{\frac{2}{s+2}}C_0,$ where \(C_0\) denotes the Euler--Lagrange level associated with the mass-one minimizer \(\overline{\rho}\). Thus, we have
\begin{equation}\label{eq:EL-reference-mass-lambda}
P_\lambda(x)+U(x)=C_\lambda
\quad\text{on }\operatorname{supp}\overline{\rho}_\lambda,
\qquad
P_\lambda(x)+U(x)\geq C_\lambda
\quad\text{in }\mathbb R^d.
\end{equation}

Since \(C_{\Gamma_0}=p(0)>0\) and \(C_\lambda=\lambda^{\frac{2}{s+2}}C_0\), the map \(\lambda\mapsto C_\lambda\) is strictly increasing from \(0\) to \(+\infty\). Hence, there exists a unique \(\lambda_*>0\) such that
\begin{equation}\label{eq:choice-lambda-level}
C_{\lambda_*}=C_{\Gamma_0}.
\end{equation}
We fix this value and set \(P_{\lambda_*}:=W*\overline{\rho}_{\lambda_*}\). Taking \(\lambda=\lambda_*\) in \eqref{eq:EL-reference-mass-lambda}, we obtain the Euler--Lagrange conditions for \(\overline{\rho}_{\lambda_*}\):
\begin{equation}\label{eq:P0EL}
P_{\lambda_*}(x)+U(x)=C_{\lambda_*}
\quad\text{on }\operatorname{supp}\overline{\rho}_{\lambda_*},
\qquad
P_{\lambda_*}(x)+U(x)\geq C_{\lambda_*}
\quad\text{in }\mathbb R^d.
\end{equation}
Moreover, since \(0\in\operatorname{int}(\operatorname{supp}\overline{\rho}_{\lambda_*})\) and \(U(0)=0\), equations \eqref{eq:central-component-level}, \eqref{eq:P0EL}, and \eqref{eq:choice-lambda-level} yield
\begin{equation}\label{eq:potentials-agree-origin}
P_{\lambda_*}(0)=C_{\lambda_*}=C_{\Gamma_0}=p(0).
\end{equation}

Our goal is to prove that \(P_{\lambda_*}\equiv p\) in \(\mathbb R^d\). Applying \(L\) will then give \(\mu=\overline{\rho}_{\lambda_*}\). Since \(\mu\) has unit mass, comparison of the total masses yields \(\lambda_*=1\), and hence \(\mu=\overline{\rho}\). We begin with a sliding argument from above.

}
\end{step}

\vspace{0.2cm}
\noindent
\begin{step}[\textit{Sliding from above}]\label{step2}{\rm

We now compare \(P_{\lambda_*}\) and \(p\) by a sliding argument from above. Define
\begin{equation}\label{eq:P_t}
\mathcal{P}_t:=P_{\lambda_*}+t,\qquad t\geq0.    
\end{equation}

Since \(p\) is continuous and \(p(x)\to0\) as \(|x|\to\infty\) by Lemma \ref{lem:farfield-decay}, and since \(P_{\lambda_*}\geq0\), we can choose \(t\) sufficiently large so that \(\mathcal{P}_t>p\) in \(\mathbb R^d\). Hence, there exists a minimal value \(t_0\geq0\) such that \(\mathcal{P}_{t_0}\geq p\) in \(\mathbb R^d\). If \(t_0=0\), then \(P_{\lambda_*}\geq p\), and there is nothing to prove. We therefore assume that \(t_0>0\). Since \(p(x)-P_{\lambda_*}(x)\to0\) as \(|x|\to\infty\), the positive maximum
\[
t_0=\max_{\mathbb R^d}\bigl(p-P_{\lambda_*}\bigr)
\]
is attained at some point \(x_0\in\mathbb R^d\). Thus, we have
\begin{equation}\label{eq:touch}
\mathcal{P}_{t_0} (x) = P_{\lambda_*}(x) + t_0 \qquad \textrm{and} \quad \mathcal{P}_{t_0}(x_0)=p(x_0).
\end{equation}

\smallskip
\noindent\textit{Case 1.} Assume that \(\mu(x_0)=0\). Set \(h:=\mathcal{P}_{t_0}-p\geq0\), so that \(h(x_0)=0\). Since \(L\) annihilates constants,
\begin{equation}\label{eq1_case1}
L[h](x_0)=L[P_{\lambda_*}](x_0)-L[p](x_0)
=\overline{\rho}_{\lambda_*}(x_0)-\mu(x_0)
=\overline{\rho}_{\lambda_*}(x_0).
\end{equation}
If \(\overline{\rho}_{\lambda_*}(x_0)>0\), then \(L[h](x_0)>0\). On the other hand, applying the maximum principle, Lemma \ref{lemma_maximum_principle}, with \(p_0=0\) and \(p_1=h\) yields \(L[h](x_0)\le0\), a contradiction.
Therefore \(\overline{\rho}_{\lambda_*}(x_0)=0\), and \eqref{eq1_case1} gives \(L[h](x_0)=0\). Since \(h\ge0\) and \(h(x_0)=0\),
Corollary \ref{corollary_principle_maximum} implies \(h\equiv0\) in \(\mathbb R^d\), that is, \(\mathcal{P}_{t_0}\equiv p\).
Hence, from \eqref{eq:touch}, we have \(t_0=0\), contradicting the assumption \(t_0>0\). Hence this case cannot occur.

\smallskip
\noindent\textit{Case 2.} Assume that \(\mu(x_0)>0\). 

We first prove that \(x_0\in\operatorname{supp}\overline{\rho}_{\lambda_*}\). Indeed, by \eqref{condition_steady_state_mu}, we have $\nabla(p+U)(x_0)=0$. Since $p-P_{\lambda_*}=(p+U)-(P_{\lambda_*}+U)$ achieves its maximum at $x_0$, we deduce that $\nabla(P_{\lambda_*}+U)(x_0) =0$. Applying Lemma \ref{lemma_convex} to $\bar{\rho}_{\lambda_*}$, we see that $\nabla(P_{\lambda_*}+U)(x)\ne0$ for every $x\notin \operatorname{supp}\bar{\rho}_{\lambda_*}$, which implies \(x_0\in\operatorname{supp}\overline{\rho}_{\lambda_*}\).

By \eqref{condition_steady_state_mu}, \eqref{eq:P0EL}, and \eqref{eq:touch},
\begin{equation}\label{eq:contact-level-component}
C_{\Gamma_{x_0}}
=
p(x_0)+U(x_0)
=
\mathcal{P}_{t_0}(x_0)+U(x_0)
=
C_{\lambda_*}+t_0.
\end{equation}

Consider the segment \(y_\tau:=\tau x_0\), \(0\leq\tau\leq1\). Since
\(\operatorname{supp}\overline{\rho}_{\lambda_*}\) is an ellipsoid centered at
the origin and \(x_0\in\operatorname{supp}\overline{\rho}_{\lambda_*}\), we have $y_\tau\in\operatorname{supp}\overline{\rho}_{\lambda_*}$ for every $\tau\in[0,1].$ Moreover, for every \(0\leq\tau<1\), the point \(y_\tau\) belongs to the interior
of \(\operatorname{supp}\overline{\rho}_{\lambda_*}\), and then $\overline{\rho}_{\lambda_*}(y_\tau)>0.$ 
We claim that
\begin{equation}\label{eq:segment-positive-mu}
\mu(y_\tau)>0
\qquad\text{for every }0\leq\tau\leq1.
\end{equation}
Suppose, by contradiction, that \(\mu(y_\tau)=0\) for some \(\tau\in[0,1]\). Since the map \(\tau\mapsto\mu(y_\tau)=\mu(\tau x_0)\) is continuous, the set
\[
Z:=\{\tau\in[0,1]:\mu(y_\tau)=0\}
\]
is closed in \([0,1]\). Moreover, \(\mu(y_0)=\mu(0)>0\) and \(\mu(y_1)=\mu(x_0)>0\), so \(Z\subset(0,1)\). Hence \(Z\) is a nonempty compact subset of \([0,1]\), and there exists \(\tau_0:=\max Z\in(0,1)\). Then \(\mu(y_{\tau_0})=0\) and \(\mu(y_\tau)>0\) for every \(\tau\in(\tau_0,1]\). Since the path \(\{y_\tau:\tau\in(\tau_0,1]\}\) is connected, contained in \(\{\mu>0\}\), and contains \(x_0=y_1\), it follows that
\[
\{y_\tau:\tau\in(\tau_0,1]\}\subset\Gamma_{x_0}.
\]
Therefore, by \eqref{condition_steady_state_mu} and by continuity of \(p+U\), it follows that
\begin{equation}\label{eq_C_gamma_0}
p(y_\tau)+U(y_\tau)=C_{\Gamma_{x_0}}
\qquad\text{for every }\tau\in[\tau_0,1]
\end{equation}
On the other hand, since \(y_{\tau_0}\in\operatorname{supp}\overline{\rho}_{\lambda_*}\), it follows from \eqref{eq:P0EL}, \eqref{eq:contact-level-component}, and \eqref{eq_C_gamma_0}, we have
\[
\mathcal{P}_{t_0}(y_{\tau_0})+U(y_{\tau_0})
=
C_{\lambda_*}+t_0 = C_{\Gamma_{x_0}} = p(y_{\tau_0}) + U(y_{\tau_0}),
\]
which implies \(\mathcal{P}_{t_0}(y_{\tau_0})=p(y_{\tau_0})\). Thus, defining \(h:=\mathcal{P}_{t_0}-p\geq0\), we have \(h(y_{\tau_0})=0\) and
\[
L[h](y_{\tau_0})
=
\overline{\rho}_{\lambda_*}(y_{\tau_0})-\mu(y_{\tau_0})
=
\overline{\rho}_{\lambda_*}(y_{\tau_0})>0.
\]
This contradicts Lemma \ref{lemma_maximum_principle} (Maximum Principle), which gives \(L[h](y_{\tau_0})\leq0\). Therefore, \eqref{eq:segment-positive-mu} holds.

Since the segment \(\{y_\tau:0\leq\tau\leq1\}\) is connected, contained in \(\{\mu>0\}\), and joins \(0\) to \(x_0\), we have
\[
\Gamma_{x_0}=\Gamma_0.
\]
Hence, \(C_{\Gamma_{x_0}}=C_{\Gamma_0}\). Combining this identity with \eqref{eq:choice-lambda-level} and \eqref{eq:contact-level-component}, we obtain
\[
C_{\lambda_*}+t_0
=
C_{\Gamma_{x_0}}
=
C_{\Gamma_0}
=
C_{\lambda_*}.
\]
Thus \(t_0=0\), contradicting the assumption \(t_0>0\). Therefore, we conclude
\begin{equation}\label{eq:P0-above-p}
P_{\lambda_*}\geq p
\qquad\text{in }\mathbb R^d.
\end{equation}

} \end{step}

\vspace{0.2cm}
\noindent
\begin{step}[\textit{\(p=P_{\lambda_*}\) on \(\operatorname{supp}\overline{\rho}_{\lambda_*}\)}]\label{step3}{\rm

We claim that \(p=P_{\lambda_*}\) on \(\operatorname{supp}\overline{\rho}_{\lambda_*}\). Set
\[
S:=\operatorname{int}\bigl(\operatorname{supp}\overline{\rho}_{\lambda_*}\bigr),
\qquad
E:=\{x\in\mathbb R^d:p(x)=P_{\lambda_*}(x)\},
\qquad
A:=S\cap E.
\]
We first prove that \(A=S\). The set \(S\) is open and connected, since it is the interior of an ellipsoid. Since \(p\) and \(P_{\lambda_*}\) are continuous, \(E\) is closed. Moreover, \eqref{eq:potentials-agree-origin} and \(0\in S\) imply that \(0\in A\), so \(A\neq\varnothing\). Therefore, \(A\) is relatively closed in \(S\). It remains to prove that \(A\) is relatively open in \(S\).

\emph{i)} We first show that the relative boundary of \(A\) in \(S\) is empty. Since \(E\) is closed,
\begin{equation}\label{defi_A_border}
\partial_S A=\partial E\cap S
=\partial\{p=P_{\lambda_*}\}\cap\operatorname{int}\bigl(\operatorname{supp}\overline{\rho}_{\lambda_*}\bigr),
\end{equation}
where \(\partial_S\) denotes the boundary relative to \(S\). Suppose that \(x_1\in\partial_S A\). Since \(x_1\in S\), there exists \(r_1>0\) such that \(B_{r_1}(x_1)\subset S\), and hence, by \eqref{eq:P0EL},
\[
P_{\lambda_*}(x)+U(x)=C_{\lambda_*}
\qquad\text{for }x\in B_{r_1}(x_1).
\]

Assume first that \(\mu(x_1)>0\). Since \(\Gamma_{x_1}\) is open, there exists \(r_2>0\) such that \(B_{r_2}(x_1)\subset\Gamma_{x_1}\). By \eqref{condition_steady_state_mu},
\[
p(x)+U(x)=C_{\Gamma_{x_1}}
\qquad\text{for }x\in B_{r_2}(x_1).
\]
Since \(x_1\in\partial E\) and \(E\) is closed, we have \(p(x_1)=P_{\lambda_*}(x_1)\). Evaluating the two identities at \(x_1\) gives \(C_{\Gamma_{x_1}}=C_{\lambda_*}\). Therefore, for \(r:=\min\{r_1,r_2\}\),
\[
p(x)+U(x)
=
C_{\Gamma_{x_1}}
=
C_{\lambda_*}
=
P_{\lambda_*}(x)+U(x)
\qquad\text{for }x\in B_r(x_1),
\]
and hence \(p=P_{\lambda_*}\) on \(B_r(x_1)\). This implies \(x_1\in\operatorname{int}(E)\), contradicting \(x_1\in\partial E\). Thus, \(\mu(x_1)=0\). By Step \ref{step2}, \(p\leq P_{\lambda_*}\) in \(\mathbb R^d\), while \(p(x_1)=P_{\lambda_*}(x_1)\). Applying Lemma \ref{lemma_maximum_principle} (Maximum Principle) at \(x_1\)  with \(p_0=p\) and \(p_1=P_{\lambda_*}\), we obtain
\[
L[p](x_1)\geq L[P_{\lambda_*}](x_1).
\]
Since \(L[p]=\mu\), \(L[P_{\lambda_*}]=\overline{\rho}_{\lambda_*}\), \(\mu(x_1)=0\), and \(x_1\in S\), it follows that
\[
0=\mu(x_1)\geq\overline{\rho}_{\lambda_*}(x_1)>0,
\]
which is a contradiction. Therefore, we have \(\partial_S A=\varnothing\).

\emph{ii)} We now prove that \(A\) is relatively open in \(S\). Let \(x\in A\). Since \(x\in E\) and \(\partial_S A=\varnothing\), the point \(x\) cannot belong to the relative boundary of \(A\) in \(S\). Since \(A\) is relatively closed in \(S\), it follows that \(x\) belongs to the relative interior of \(A\). Thus, \(A\) is relatively open in \(S\).

Therefore, \(A\) is nonempty, relatively open, and relatively closed in the connected set \(S\). Hence \(A=S\), and consequently \(p=P_{\lambda_*}\) in \(S\). By continuity, this equality extends to \(\overline S=\operatorname{supp}\overline{\rho}_{\lambda_*}\). Thus,
\begin{equation}\label{eq:p-equals-P0-support}
p=P_{\lambda_*}
\qquad\text{on }\operatorname{supp}\overline{\rho}_{\lambda_*}.
\end{equation}

}
\end{step}

\vspace{0.2cm}
\noindent
\begin{step}[\textit{\(p=P_{\lambda_*}\) in \(\mathbb R^d\)}]\label{step4}{\rm

From Steps \ref{step2} and \ref{step3}, we have
\[
p\leq P_{\lambda_*}
\quad\text{in }\mathbb R^d \quad\textrm{and}
\quad
p=P_{\lambda_*}
\quad\text{on }\operatorname{supp}\overline{\rho}_{\lambda_*}.
\]
We now perform a sliding argument from below. For \(t\geq0\), define
\[
Q_t:=P_{\lambda_*}-t.
\]
Since \(P_{\lambda_*}\) is bounded and \(p\geq0\), we have \(Q_t\leq p\) in \(\mathbb R^d\) for every sufficiently large \(t\). Set
\[
t_1:=\inf\bigl\{t\geq0:Q_t\leq p\text{ in }\mathbb R^d\bigr\}.
\]
Then \(Q_{t_1}\leq p\) in \(\mathbb R^d\). If \(t_1=0\), then \(P_{\lambda_*}\leq p\), and the conclusion follows immediately from \(p\leq P_{\lambda_*}\). We therefore assume by contradiction that \(t_1>0\). Since \(P_{\lambda_*}(x)-p(x)\to0\) as \(|x|\to\infty\), the positive maximum \(t_1=\max_{\mathbb R^d}(P_{\lambda_*}-p)\) is attained at some point \(y_0\in\mathbb R^d\). Hence,
\begin{equation}\label{eq:touch-below}
Q_{t_1}(y_0)=p(y_0).
\end{equation}
By \eqref{eq:p-equals-P0-support}, \(P_{\lambda_*}-p=0\) on \(\operatorname{supp}\overline{\rho}_{\lambda_*}\). Since \(t_1>0\), the touching point satisfies
\[
y_0\notin\operatorname{supp}\overline{\rho}_{\lambda_*}.
\]

Set \(w:=p-Q_{t_1}\geq0\). Then \(w(y_0)=0\). Since \(L\) annihilates constants, applying Lemma \ref{lemma_maximum_principle} at \(y_0\) gives
\[
L[p](y_0)\leq L[Q_{t_1}](y_0)=L[P_{\lambda_*}](y_0).
\]
Therefore,
\[
\mu(y_0)\leq\overline{\rho}_{\lambda_*}(y_0)=0.
\]
Since \(\mu\geq0\), we conclude that \(\mu(y_0)=0\). Consequently,
\[
L[w](y_0)=L[p](y_0)-L[Q_{t_1}](y_0)
=\mu(y_0)-\overline{\rho}_{\lambda_*}(y_0)=0.
\]
Since \(w\geq0\) and \(w(y_0)=0\), Corollary \ref{corollary_principle_maximum} implies that \(w\equiv0\) in \(\mathbb R^d\). Hence,
\[
p=Q_{t_1}=P_{\lambda_*}-t_1
\qquad\text{in }\mathbb R^d.
\]
On the other hand, \eqref{eq:p-equals-P0-support} gives \(p=P_{\lambda_*}\) on the nonempty set \(\operatorname{supp}\overline{\rho}_{\lambda_*}\). Therefore, \(t_1=0\), contradicting the assumption \(t_1>0\). We conclude that \(t_1=0\), and hence \(P_{\lambda_*}\leq p\) in \(\mathbb R^d\). Combining this with Step \ref{step2}, we obtain
\[
p\equiv P_{\lambda_*}
\qquad\text{in }\mathbb R^d.
\]

Applying \(L\) and recalling that \(P_{\lambda_*}=W*\overline{\rho}_{\lambda_*}\), we obtain
\[
\mu=L[p]=L[P_{\lambda_*}]=\overline{\rho}_{\lambda_*}.
\]
Since \(\mu\) has total mass \(1\), while \(\overline{\rho}_{\lambda_*}\) has total mass \(\lambda_*\), it follows that
\[
1=\int_{\mathbb R^d}\mu(x)dx
=\int_{\mathbb R^d}\overline{\rho}_{\lambda_*}(x)dx
=\lambda_*.
\]
Therefore, \(\lambda_*=1\) and \(\mu=\overline{\rho}\), completing the proof.

}
\end{step}
\end{proof}

\begin{corollary}
\label{cor:connected-support-steady-state}
Under the assumptions of Theorem \ref{thm:uniqueness}, if
\(\mu=\rho dx\) is a steady state, then
\[
\operatorname{supp}\rho
=
\operatorname{supp}\overline{\rho}.
\]
In particular, \(\operatorname{supp}\rho\) is connected.
\end{corollary}

\section{Proof of Theorem \ref{thm:Lp_global_confined}}

\subsection{\texorpdfstring{\(L^p\) estimates}{Lp estimates}}
\label{subsec:Lp-estimates}

We first prove the \(L^p\)-estimate in
Theorem \ref{thm:Lp_global_confined}\textup{(i)}. The argument is based on the
flow-interchange method applied to the entropy
\[
\mathcal G_p(u):=\frac{1}{p-1}\int_{\mathbb R^d}u(x)^pdx,
\qquad 1<p<\infty.
\]
We shall use the ellipticity of the kernel in $\eqref{ha}$. Namely, the
quadratic form associated with \(A=-\Delta L^{-1}\) controls the standard
fractional seminorm of order \(\alpha/2\). More precisely, there exists
\(c_A>0\), depending only on the ellipticity constants in
$\eqref{ha}$, such that
\begin{equation}\label{ine:norm_sobolev_A}
\int_{\mathbb R^d} f A[f]dx
\ge
c_A\|f\|_{\dot H^{\alpha/2}(\mathbb R^d)}^2, \qquad f \in \mathcal{S}(\R^d).
\end{equation}
The positivity of the kernel then allows
us to control the nonlocal part of the dissipation through a
Stroock--Varopoulos inequality. The confinement term is controlled by the upper
bound on \(\Delta U\). Together with a fractional interpolation inequality, this
yields a discrete logistic estimate for
\(\|u_\tau^k\|_{L^p(\mathbb R^d)}\). Iterating the inequality gives a bound
uniform in the time step, which then passes to the limiting curve obtained in
Theorem \ref{thm:global-wellposedness}.

\begin{proof}[Proof of Theorem \ref{thm:Lp_global_confined}\textup{(i)}]
Since
\[
u_0\in L^p(\mathbb R^d)\cap\mathcal P_2(\mathbb R^d),
\]
we also have $u_0\in D(H)$: the positive part of $u_0\log u_0$ is
controlled by the $L^p$-norm, while the finite second moment controls
its negative part. Let 
\[
u\in D(\mathcal G_p),
\qquad
\mathcal E_U(u)<\infty
\]
be fixed, and let $(S_h^{(p)})_{h\ge0}$ denote the $2$-Wasserstein
gradient flow of $\mathcal G_p$. Notice that, if
$w_h:=S_h^{(p)}u$, then
\begin{equation}\label{eq:PME_auxiliary_Lp}
\partial_h w_h=\Delta(w_h^p),
\qquad
w_h|_{h=0}=u.
\end{equation}
Since $u_\tau^0=u_0\in D(\mathcal G_p)$, the flow-interchange theorem,
applied inductively to the minimizing movement sequence $(u_\tau^k)_k$,
gives $u_\tau^k\in D(\mathcal G_p)$ and
\begin{equation}\label{eq:flow-interchange-Lp}
\mathcal G_p(u_\tau^k)
+
\tau D_{\mathcal G_p}[\mathcal E_U](u_\tau^k)
\le
\mathcal G_p(u_\tau^{k-1}),
\end{equation}
for every $k\ge1$, where
\[
D_{\mathcal G_p}[\mathcal E_U](u)
:=
\limsup_{h\to0^+}
\frac{\mathcal E_U(u)-\mathcal E_U(S_h^{(p)}u)}{h}.
\]

We estimate the dissipation of $\mathcal E_U$ along the auxiliary flow.
Recall that
\[
\mathcal E_U(u)=F(u)+\int_{\mathbb R^d}U(x)u(x)dx,
\qquad
F(u)=\frac12\int_{\mathbb R^d}(L^{-1}u)(x)u(x)dx.
\]
The following calculations may first be performed for smooth
approximations of the auxiliary flow and then passed to the limit. The
coercivity estimate \eqref{ine:norm_sobolev_A}, initially stated for
Schwartz functions, extends by density and lower semicontinuity to the
form domain of $A$.

Since $\delta F/\delta u=L^{-1}u$, setting $v_h:=L^{-1}w_h$, we obtain
\[
\frac{d}{dh}F(w_h)
=
\int_{\mathbb R^d}v_h\partial_h w_hdx
=
\int_{\mathbb R^d}v_h\Delta(w_h^p)dx
=
\int_{\mathbb R^d}w_h^p\Delta v_hdx.
\]
Since $A=-\Delta L^{-1}$, we have $A[w_h]=-\Delta v_h$, and hence
\begin{equation}\label{eq:dFdh-A-Lp}
\frac{d}{dh}F(w_h)
=
-\int_{\mathbb R^d}w_h^p A[w_h]dx := - D_p(w_h).
\end{equation}
Using the singular-integral representation in $\eqref{ha}$, the
symmetry of $K_A$, and the positivity of the kernel, we get
\begin{equation*} 
D_p(w_h)
=
\frac12
\iint_{\mathbb R^d\times\mathbb R^d}
\big(w_h(x)^p-w_h(y)^p\big)
\big(w_h(x)-w_h(y)\big)
K_A(x-y)dxdy.
\end{equation*}
Using that,
\[
\big(a^p-b^p\big)(a-b)
\ge
\frac{4p}{(p+1)^2}
\left(
a^{\frac{p+1}{2}}-b^{\frac{p+1}{2}}
\right)^2,
\qquad a,b\ge0,
\]
we obtain
\begin{equation}\label{eq:Dp-kernel-lower-Lp}
D_p(w_h)
\ge
\frac{2p}{(p+1)^2}
\iint_{\mathbb R^d\times\mathbb R^d}
\left|
w_h(x)^{\frac{p+1}{2}}
-
w_h(y)^{\frac{p+1}{2}}
\right|^2
K_A(x-y)dxdy.
\end{equation}
By the coercivity of the quadratic form recalled in
\eqref{ine:norm_sobolev_A}, it follows from
\eqref{eq:Dp-kernel-lower-Lp} that
\begin{equation}\label{eq:Dp-Halpha-lower-Lp}
D_p(w_h)
\ge
\frac{4p}{(p+1)^2}c_A
\left\|w_h^{\frac{p+1}{2}}\right\|_{\dot H^{\alpha/2}(\mathbb R^d)}^2.
\end{equation}

Let $S_{d,\alpha/2}$ denote the constant in the fractional Sobolev inequality, and set $q:= \frac{d(p+1)}{d-\alpha}$. Then
\[
\|\phi\|_{L^q(\mathbb R^d)}^{p+1}
=
\left\|\phi^{\frac{p+1}{2}}\right\|_{L^{\frac{2d}{d-\alpha}}(\mathbb R^d)}^2
\le
S_{d,\alpha/2}^2
\left\|\phi^{\frac{p+1}{2}}\right\|_{\dot H^{\alpha/2}(\mathbb R^d)}^2.
\]
Interpolating $L^p$ between $L^1$ and $L^q$, we obtain the fractional
interpolation inequality
\begin{equation}\label{eq:LMS-interpolation-Lp}
\|\phi\|_{L^p(\mathbb R^d)}^{p+1}
\le
S_{d,\alpha/2}^{2\theta_p}
\|\phi\|_{L^1(\mathbb R^d)}^{(1-\theta_p)(p+1)}
\left\|\phi^{\frac{p+1}{2}}\right\|_{\dot H^{\alpha/2}(\mathbb R^d)}^{2\theta_p},
\quad \textrm{where} \quad
\theta_p:=\frac{d(p^2-1)}{p(\alpha+dp)}\in(0,1).
\end{equation}
Applying \eqref{eq:LMS-interpolation-Lp} to $\phi=w_h$, and using
$\|w_h\|_{L^1(\mathbb R^d)}=1$, we get
\[
\left\|w_h^{\frac{p+1}{2}}\right\|_{\dot H^{\alpha/2}(\mathbb R^d)}^2
\ge
S_{d,\alpha/2}^{-2}
\|w_h\|_{L^p(\mathbb R^d)}^{\frac{p+1}{\theta_p}}.
\]
Combining this with \eqref{eq:Dp-Halpha-lower-Lp}, we obtain
\begin{equation}\label{eq:Dp-Lp-lower}
D_p(w_h)
\ge
B_p\|w_h\|_{L^p}^{\frac{p+1}{\theta_p}}
=
B_p\left(\|w_h\|_{L^p}^p\right)^{\beta_p},
\end{equation}
where
\begin{equation*} 
B_p
:=
\frac{4p}{(p+1)^2}c_A S_{d,\alpha/2}^{-2},
\qquad
\beta_p
:=
\frac{p+1}{p\theta_p}
=
1+\frac{d+\alpha}{d(p-1)}.
\end{equation*}

Finally, let us estimate the confinement contribution. Set
\[
\widetilde U:=U+a_U.
\]
Then $\widetilde U\ge c_U|x|^2\ge0$. Choose nondecreasing concave
functions $\chi_R\in C^2([0,\infty))$ such that
\[
0\le\chi_R'\le1,
\qquad
\chi_R''\le0,
\]
with $\chi_R(r)=r$ for $r\le R$, $\chi_R$ constant for $r\ge2R$, and
$\chi_R(r)\uparrow r$ as $R\to\infty$. Since $\widetilde U$ is
coercive, the derivatives of $\chi_R(\widetilde U)$ have compact support.
Moreover,
\[
\Delta\bigl(\chi_R(\widetilde U)\bigr)
=
\chi_R'(\widetilde U)\Delta U
+
\chi_R''(\widetilde U)|\nabla U|^2
\le
\Lambda_U.
\]
Testing \eqref{eq:PME_auxiliary_Lp} against
$\chi_R(\widetilde U)$, after subtracting its constant value at
infinity, and integrating in time, we obtain
\[
\int_{\mathbb R^d}\chi_R(\widetilde U)w_hdx
-
\int_{\mathbb R^d}\chi_R(\widetilde U)u dx
\le
\Lambda_U\int_0^h\|w_r\|_{L^p(\mathbb R^d)}^pdr.
\]
Letting $R\to\infty$ and using monotone convergence, together with
mass conservation, yields
\begin{equation}\label{eq:confinement-general-Lp}
\int_{\mathbb R^d}U(x)w_h(x)dx
-
\int_{\mathbb R^d}U(x)u(x)dx
\le
\Lambda_U\int_0^h\|w_r\|_{L^p(\mathbb R^d)}^pdr.
\end{equation}

Combining \eqref{eq:dFdh-A-Lp}, \eqref{eq:Dp-Lp-lower}, and
\eqref{eq:confinement-general-Lp}, we obtain
\[
\mathcal E_U(u)-\mathcal E_U(w_h)
\ge
\int_0^h
\left[
B_p\left(\|w_r\|_{L^p}^p\right)^{\beta_p}
-
\Lambda_U\|w_r\|_{L^p}^p
\right]dr.
\]
Dividing by $h$, taking the limit $h\to0^+$, and using that
$S_h^{(p)}u\to u$ strongly in $L^p(\mathbb R^d)$, it follows that
\begin{equation}\label{eq:D-Gp-E-lower-Lp}
D_{\mathcal G_p}[\mathcal E_U](u)
\ge
B_p\|u\|_{L^p}^{p\beta_p}
-
\Lambda_U\|u\|_{L^p}^p.
\end{equation}

Applying \eqref{eq:D-Gp-E-lower-Lp} to $u=u_\tau^k$, and using
\[
\mathcal G_p(u_\tau^k)
=
\frac{1}{p-1}\|u_\tau^k\|_{L^p(\mathbb R^d)}^p,
\]
we obtain from \eqref{eq:flow-interchange-Lp} that
\begin{equation}\label{eq:discrete-Lp-logistic}
\|u_\tau^k\|_{L^p}^p
+
\tau(p-1)B_p\|u_\tau^k\|_{L^p}^{p\beta_p}
-
\tau(p-1)\Lambda_U\|u_\tau^k\|_{L^p}^p
\le
\|u_\tau^{k-1}\|_{L^p}^p.
\end{equation}
Set
\begin{equation}\label{eq:Ap-Lp}
A_p:=(p-1)B_p
=
\frac{4p(p-1)}{(p+1)^2}c_A S_{d,\alpha/2}^{-2}.
\end{equation}
Then, we can rewrite \eqref{eq:discrete-Lp-logistic} as
\begin{equation}\label{eq:discrete-Lp-bracket}
\|u_\tau^k\|_{L^p}^p
\left[
1+\tau\left(
A_p\|u_\tau^k\|_{L^p}^{p(\beta_p-1)}
-
(p-1)\Lambda_U
\right)
\right]
\le
\|u_\tau^{k-1}\|_{L^p}^p.
\end{equation}

Define
\begin{equation*} 
C_p
:=
\left(
\frac{(p-1)\Lambda_U}{A_p}
\right)^{\frac{1}{p(\beta_p-1)}}.
\end{equation*}
We claim that
\begin{equation}\label{eq:Lp-discrete-bound}
\|u_\tau^k\|_{L^p(\mathbb R^d)}
\le
\max\left\{
\|u_\tau^0\|_{L^p(\mathbb R^d)},C_p
\right\}
\qquad\text{for every }k\ge0.
\end{equation}
The claim is trivial for $k=0$. Assume it holds up to $k-1$. If
$\|u_\tau^k\|_{L^p}\le C_p$, there is nothing to prove. If
$\|u_\tau^k\|_{L^p}>C_p$, then by the definition of $C_p$,
\[
A_p\|u_\tau^k\|_{L^p}^{p(\beta_p-1)}-(p-1)\Lambda_U>0.
\]
Hence the bracket in \eqref{eq:discrete-Lp-bracket} is strictly larger
than $1$, and so
\[
\|u_\tau^k\|_{L^p}
<
\|u_\tau^{k-1}\|_{L^p}
\le
\max\left\{
\|u_\tau^0\|_{L^p},C_p
\right\},
\]
where we used the induction hypothesis in the last inequality. This proves
\eqref{eq:Lp-discrete-bound}.

Since $u_\tau^0=u_0$, for the piecewise constant interpolation we obtain
\begin{equation}\label{eq:Lp-discrete-final}
\sup_{t\ge0}
\|u_\tau(t)\|_{L^p(\mathbb R^d)}
\le
\max\left\{
\|u_0\|_{L^p(\mathbb R^d)},C_p
\right\}.
\end{equation}
Using $\beta_p-1=\frac{d+\alpha}{d(p-1)}$ and
\eqref{eq:Ap-Lp}, the threshold can be written as
\begin{equation*} 
C_p
=
\left[
\frac{\Lambda_U(p+1)^2S_{d,\alpha/2}^{2}}
{4pc_A}
\right]^{\frac{d(p-1)}{p(d+\alpha)}}.
\end{equation*}

It remains to pass to the limit $\tau\to0$. Let $(\tau_n)_n$ be the
subsequence from Theorem \ref{thm:global-wellposedness}. Fix $t\ge0$. From
\eqref{eq:Lp-discrete-final}, the sequence $(u_{\tau_n}(t))_n$ is bounded in
$L^p(\mathbb R^d)$. Since $1<p<\infty$, every subsequence has a further
subsequence converging weakly in $L^p(\mathbb R^d)$. On the other hand, by
the convergence of the JKO scheme,
\[
u_{\tau_n}(t)\rightharpoonup u(t)
\qquad\text{narrowly in }\mathcal P_2(\mathbb R^d).
\]
Every weak $L^p$ cluster point therefore agrees with $u(t)$ when tested
against $C_c^\infty(\mathbb R^d)$. Hence $u(t)$ has an $L^p$-density and
\[
u_{\tau_n}(t)\rightharpoonup u(t)
\qquad\text{weakly in }L^p(\mathbb R^d).
\]
By the lower semicontinuity of the $L^p$-norm,
\[
\|u(t)\|_{L^p(\mathbb R^d)}
\le
\liminf_{n\to\infty}
\|u_{\tau_n}(t)\|_{L^p(\mathbb R^d)}
\le
\max\left\{
\|u_0\|_{L^p(\mathbb R^d)},C_p
\right\}.
\]
Therefore, we conclude the proof of item \textup{(i)} in
Theorem \ref{thm:Lp_global_confined}.
\end{proof}

\subsection{\texorpdfstring{\(L^\infty\) and H\"older estimates}{Linfty and Holder estimates}}
\label{subsec:Linfty-holder-estimates}

We prove the \(L^\infty\) and H\"older estimates in
Theorem \ref{thm:Lp_global_confined}\textup{(ii)--(iii)}. Throughout this
section, we restrict to the quadratic confinement \(U(x)=|x|^2\). Let \(u\)
be the unit-mass limiting curve obtained in
Theorem \ref{thm:global-wellposedness}. Using the self-similar scaling, we
associate with \(u\) an unconfined curve \(\rho\). We first derive the
\(L^\infty\)-decay for \(\rho\) by a truncation argument, and then combine
the local H\"older regularity established in
Appendix \ref{section_appendix} with the scaling invariance of the
unconfined equation to obtain the corresponding global H\"older decay.
Finally, both estimates are transferred to the confined variables through
the same scaling.

Let us start with the \(L^\infty\)-estimate. Define
\begin{equation}\label{eq:scaling-functions-Linfty}
R(t):=e^{2t},
\qquad
\tau(t):=\frac{e^{2(s+2)t}-1}{2(s+2)}.
\end{equation}
For the confined solution \(u\) given by Theorem \ref{thm:global-wellposedness}, we define the associated unconfined curve \(\rho\) through
\begin{equation}\label{eq:scaling-confined-unconfined}
u(t,x)
=
R(t)^d\rho(\tau(t),R(t)x),
\end{equation}
or equivalently,
\[
\rho(\tau(t),y)
=
R(t)^{-d}
u\left(t,\frac{y}{R(t)}\right).
\]
Since $\tau$ is strictly increasing from $[0,\infty)$ onto $[0,\infty)$, the relation $\rho(\tau,\cdot)$ is well-defined for every $\tau \geq 0$. Moreover, since the scaling preserves the total mass, using the homogeneity of $W$ together with the relations $R'(t) = 2R(t) $ and $\tau'(t) = R(t)^{s+2}$ allows the scaling \eqref{eq:scaling-confined-unconfined} to transform the confined equation into unconfined equation, namely

\begin{equation}\label{eq:unconfined-equation-Linfty}
\partial_\tau\rho
=
\nabla_y\cdot\bigl(\rho\nabla_y(W*\rho)\bigr).
\end{equation}

\begin{theorem}[Unconfined \(L^\infty\)-decay]
\label{thm:unconfined-Linfty-decay}
Fix \(d\geq2\) and \(d-2<s<d\), and set
\(\alpha:=s-d+2\in(0,2)\). Assume \eqref{ha}, and let \(\rho\) be the unconfined curve associated with \(u\) through \eqref{eq:scaling-confined-unconfined}. Then there exists a constant \(C_\infty>0\), depending only on \(d\), \(s\), and the ellipticity constants of \(A\), such that
\begin{equation}\label{eq:unconfined-Linfty-decay}
\|\rho(\tau)\|_{L^\infty(\mathbb R^d)}
\leq
C_\infty\tau^{-\frac{d}{s+2}},
\qquad
\tau>0.
\end{equation}
\end{theorem}

\begin{proof}
We adapt the truncation argument used in the proof of
\cite[Theorem 7.2]{Lisini_Mainini_Segatti_2018} to the present setting. 
In \cite{Lisini_Mainini_Segatti_2018}, the unconfined evolution is obtained
through a minimizing-movement scheme, while here \(\rho\) is obtained from
the confined solution \(u\) through the scaling
\eqref{eq:scaling-confined-unconfined}. Accordingly, we derive the
truncation estimate directly for \(\rho\) and then apply the same
De Giorgi level-set iteration. For \(M>0\), define
\[
\rho_M:=(\rho-M)_+=\max\{\rho-M,0\}.
\]
We can write the evolution of its $L^2$-norm as
\[
\frac{d}{d\tau}\|\rho_M\|_{L^2(\mathbb R^d)}^2
=
2\int_{\mathbb R^d}\rho_M
\nabla\cdot\bigl(\rho\nabla(W*\rho)\bigr)dx.
\]
On the set \(\{\rho>M\}\), we have
\(\rho=M+\rho_M\) and
\(\nabla\rho=\nabla\rho_M\) a.e. Hence, it follows from integration by parts that
\[
\begin{aligned}
\frac{d}{d\tau}\|\rho_M\|_{L^2}^2
&=
-2\int_{\mathbb R^d}
\rho\nabla\rho_M\cdot\nabla(W*\rho)dx =
-\int_{\mathbb R^d}
\nabla\bigl(\rho_M^2+2M\rho_M\bigr)
\cdot\nabla(W*\rho)dx \\
&=
\int_{\mathbb R^d}
\bigl(\rho_M^2+2M\rho_M\bigr)
\Delta(W*\rho)dx.
\end{aligned}
\]
Recalling that $-\Delta(W*\rho)=A[\rho],$ we obtain
\[
\frac{d}{d\tau}\|\rho_M\|_{L^2(\mathbb R^d)}^2
=
-\int_{\mathbb R^d}
\bigl(\rho_M^2+2M\rho_M\bigr)A[\rho]dx.
\]
In order to avoid the use of the chain rule and the integration by parts, that cannot be rigorously proved at this stage, one can use the approximation procedure by the minimizing scheme of the unconfined equation \eqref{eq:unconfined-equation-Linfty} as in \cite[Theorem 7.2]{Lisini_Mainini_Segatti_2018}. A complete proof is omitted, since it is a direct adaptation of proof in \cite[Theorem 7.2]{Lisini_Mainini_Segatti_2018} for the unconfined equation \eqref{eq:unconfined-equation-Linfty}, to our confined equation \eqref{eq:our-evolution}.

As in the symmetrization argument used in the proof of the \(L^p\)-estimate,
using the singular-integral representation of \(A\) and the symmetry of
\(K_A\), we have
\[
\begin{aligned}
\int_{\mathbb R^d}
\bigl(\rho_M^2+2M\rho_M\bigr)A[\rho]dx
&=
\frac12
\iint_{\mathbb R^d\times\mathbb R^d}
\Big[
\rho_M(x)^2-\rho_M(y)^2 \\
&\qquad\qquad
+2M\bigl(\rho_M(x)-\rho_M(y)\bigr)
\Big]
\bigl(\rho(x)-\rho(y)\bigr)
K_A(x-y)dxdy.
\end{aligned}
\]
Since \(r\mapsto r_+\) is monotone, then $(\rho_M(x)^2-\rho_M(y)^2)(\rho(x)-\rho(y)\bigr)\geq 0$, it follows 
\[
\begin{aligned}
\int_{\mathbb R^d}
\bigl(\rho_M^2+2M\rho_M\bigr)A[\rho]dx
&\geq
M
\iint_{\mathbb R^d\times\mathbb R^d}
\bigl(\rho_M(x)-\rho_M(y)\bigr)
\bigl(\rho(x)-\rho(y)\bigr)
K_A(x-y)dxdy.
\end{aligned}
\]
Now, we distinguish cases. if $\rho(x)\geq M$ and $\rho(y)\geq M$, then 
$$
\bigl(\rho_M(x)-\rho_M(y)\bigr)
\bigl(\rho(x)-\rho(y)\bigr)=\bigl(\rho_M(x)-\rho_M(y)\bigr)^2.
$$
On the other hand, assume without loss of generality that $\rho(x)\geq M \geq\rho(y)\geq 0$, then $\rho_M(x) \bigl(\rho(x)-\rho(y)\bigr) \geq
\rho_M(x)^2=\bigl(\rho_M(x)-\rho_M(y)\bigr)^2$. Using the positivity of \(K_A\), we conclude that
\[
\begin{aligned}
\int_{\mathbb R^d}
\bigl(\rho_M^2+2M\rho_M\bigr)A[\rho]dx
&\geq
M
\iint_{\mathbb R^d\times\mathbb R^d}
\bigl(\rho_M(x)-\rho_M(y)\bigr)^2
K_A(x-y)dxdy.
\end{aligned}
\]
By the same coercivity estimate \eqref{ine:norm_sobolev_A} used in the
\(L^p\)-argument, it follows that
\[
\int_{\mathbb R^d}
\bigl(\rho_M^2+2M\rho_M\bigr)A[\rho]dx
\geq
cM\|\rho_M\|_{\dot H^{\alpha/2}(\mathbb R^d)}^2, \quad \alpha \in(0,2),
\]
where \(c>0\) depends only on the ellipticity constants of \(A\).
Therefore, we obtain
\[
\frac{d}{d\tau}\|\rho_M\|_{L^2(\mathbb R^d)}^2
+
cM\|\rho_M\|_{\dot H^{\alpha/2}(\mathbb R^d)}^2
\leq0.
\]
Integrating in time, for \(0<\tau_1<\tau_2\), yields
\begin{equation}\label{eq:unconfined-truncation-Linfty}
\|\rho_M(\tau_2)\|_{L^2(\mathbb R^d)}^2
+
cM
\int_{\tau_1}^{\tau_2}
\|\rho_M(\tau)\|_{\dot H^{\alpha/2}(\mathbb R^d)}^2d\tau
\leq
\|\rho_M(\tau_1)\|_{L^2(\mathbb R^d)}^2.
\end{equation}

We first establish the positive-time \(L^2\)-decay needed to initialize
the De Giorgi iteration below. In particular, this estimate provides the
starting \(L^2\)-control with the correct time scaling for the level-set
energies introduced later. To this end, we repeat for the unconfined
equation the symmetrization argument used in the proof of the
\(L^p\)-estimate. Combining the Stroock--Varopoulos and coercivity
estimates
\eqref{eq:Dp-kernel-lower-Lp}--\eqref{eq:Dp-Halpha-lower-Lp}
with the fractional interpolation inequality
\eqref{eq:LMS-interpolation-Lp}, and arguing as in
\eqref{eq:Dp-Lp-lower} with \(p=2\), we obtain
\[
\frac{d}{d\tau}
\|\rho(\tau)\|_{L^2(\mathbb R^d)}^2
+
c
\left(
\|\rho(\tau)\|_{L^2(\mathbb R^d)}^2
\right)^{2+\alpha/d}
\leq0.
\]
By the energy-space regularity of \(u\) given by
Theorem \ref{thm:global-wellposedness}, we have
\(u(t)\in L^2(\mathbb R^d)\) for a.e. \(t>0\). By the scaling
\eqref{eq:scaling-confined-unconfined}, it follows that
\(\rho(\tau)\in L^2(\mathbb R^d)\) for a.e. \(\tau>0\).
Hence there exists a sequence \(\varepsilon_n\downarrow0\) such that
\(\rho(\varepsilon_n)\in L^2(\mathbb R^d)\). Integrating the preceding
differential inequality from \(\varepsilon_n\) to \(\tau\) and letting
\(n\to\infty\), we obtain
\begin{equation}\label{eq:unconfined-L2-decay}
\|\rho_M(\tau)\|_{L^2(\mathbb R^d)}^2 \leq \|\rho(\tau)\|_{L^2(\mathbb R^d)}^2
\leq
C_2\tau^{-\frac{d}{d+\alpha}},
\qquad \tau>0.
\end{equation}
Thus, \eqref{eq:unconfined-L2-decay} provides the positive-time \(L^2\)
control needed to initialize the De Giorgi iteration.

Now, we briefly perform the continuous De Giorgi level-set iteration. Fix
\(\tau>0\), and define
\[
\tau_k:=\tau\bigl(1-2^{-k-1}\bigr),
\qquad
M_k:=M\bigl(1-2^{-k-1}\bigr),
\qquad
\rho_k:=(\rho-M_k)_+.
\]
Set
\[
A_k
:=
\operatorname*{esssup}_{t\in[\tau_k,\tau]}
\|\rho_k(t)\|_{L^2}^2
+
cM_k
\int_{\tau_k}^{\tau}
\|\rho_k(t)\|_{\dot H^{\alpha/2}}^2dt.
\]
Applying \eqref{eq:unconfined-truncation-Linfty} between consecutive
time levels and averaging over the initial time, and then using the
fractional Sobolev inequality together with the separation
\(M_k-M_{k-1}=M2^{-k-1}\) and $\tau_k - \tau_{k-1} = \tau 2^{-k-1}$, we obtain
\[
A_k
\leq
\frac{Cb^k}{\tau M^{1+2\alpha/d}}
A_{k-1}^{1+\alpha/d},
\qquad k\geq1,
\]
for some \(b>1\). Moreover, the same truncation estimate, combined with
\eqref{eq:unconfined-L2-decay}, yields
\[
A_0
\leq
C\tau^{-\frac{d}{d+\alpha}}.
\]
Choosing $M=L\tau^{-\frac{d}{d+\alpha}}$, we get 
\[
\frac{A_0^{\alpha/d}}
{\tau M^{1+2\alpha/d}}
\leq
C L^{-(1+2\alpha/d)},
\]
Therefore, for \(L\) sufficiently large, the standard nonlinear
De Giorgi iteration, as in the level-set argument of
\cite[Theorem 7.2]{Lisini_Mainini_Segatti_2018}, gives
\(A_k\to0\). Since \(M_k\uparrow M\) and \(\tau_k\uparrow\tau\), we conclude
that
\[
\|\rho(\tau)\|_{L^\infty(\mathbb R^d)}
\leq
M
\leq
C_\infty\tau^{-\frac{d}{d+\alpha}}
=
C_\infty\tau^{-\frac{d}{s+2}},
\]
which proves \eqref{eq:unconfined-Linfty-decay}.
\end{proof}

\begin{proof}[Proof of Theorem \ref{thm:Lp_global_confined}\textup{(ii)}]
Let \(\rho\) be the unconfined curve associated with \(u\) through \eqref{eq:scaling-confined-unconfined}. By Theorem \ref{thm:unconfined-Linfty-decay}, we have
\[
\|\rho(\tau)\|_{L^\infty(\mathbb R^d)}
\leq
C_\infty\tau^{-\frac{d}{s+2}},
\qquad
\tau>0.
\]
From \eqref{eq:scaling-confined-unconfined}, we have
\[
\|u(t)\|_{L^\infty(\mathbb R^d)}
=
R(t)^d
\|\rho(\tau(t))\|_{L^\infty(\mathbb R^d)} \leq C_\infty R(t)^d\tau(t)^{-\frac{d}{s+2}}.
\]
Invoking \eqref{eq:scaling-functions-Linfty}, we conclude
\begin{equation}\label{eq:Linfty-confined-explicit}
\|u(t)\|_{L^\infty(\mathbb R^d)}
\leq
C_\infty e^{2dt}
\left(
\frac{e^{2(s+2)t}-1}{2(s+2)}
\right)^{-\frac{d}{s+2}},
\qquad
t>0.
\end{equation}
Moreover, since
\[
e^{2dt}
\left(
e^{2(s+2)t}-1
\right)^{-\frac{d}{s+2}}
=
\left(
1-e^{-2(s+2)t}
\right)^{-\frac{d}{s+2}},
\]
we conclude that the right-hand side of \eqref{eq:Linfty-confined-explicit} is uniformly bounded for \(t\geq t_0>0\).
\end{proof}

We now derive the decay of the H\"older seminorm for the unconfined curve \(\rho\). While the local H\"older regularity for the unconfined equation is established in Appendix \ref{section_appendix}, combining this result with the previously obtained \(L^\infty\)-decay and the scaling invariance of the equation yields the corresponding global H\"older decay. This estimate is then transferred to the confined variables via self-similar scaling.

\begin{theorem}[Unconfined H\"older decay]
\label{thm:unconfined-holder-decay}
Fix \(d\geq2\) and \(d-2<s<d\), and set
\(\alpha:=s-d+2\in(0,2)\). Assume \eqref{ha}, and let \(\rho\) be the
unconfined curve associated with \(u\) through
\eqref{eq:scaling-confined-unconfined}. Then there exists
\(\eta=\eta(d,s,\Omega)\in(0,1)\) and a constant
\(C_\eta=C_\eta(d,s,\Omega,C_\infty)>0\), where \(C_\infty\) is the
constant in \eqref{eq:unconfined-Linfty-decay}, such that
\[
[\rho(\tau)]_{C^\eta(\mathbb R^d)}
\leq
C_\eta
\tau^{-\frac{d+\eta}{s+2}},
\qquad
\tau>0.
\]
\end{theorem}

\begin{proof}
Fix \(\tau>0\) and define
\begin{equation}\label{eq:v_rho_holder_scaling}
v_\tau(\xi,t)
:=
\tau^{\frac{d}{s+2}}
\rho\left(
\tau t,
\tau^{\frac1{s+2}}\xi
\right),
\qquad
t>0.
\end{equation}
By the scaling invariance of the unconfined equation,
\(v_\tau\) is again a nonnegative unit-mass solution of
\eqref{eq:unconfined-equation-Linfty}.

Let \(0<T_1<1<T_2\) be fixed. By
\eqref{eq:unconfined-Linfty-decay}, for every \(t\in(T_1,T_2)\), we have
\[
\|v_\tau(t)\|_{L^\infty(\mathbb R^d)}
=
\tau^{\frac{d}{s+2}}
\|\rho(\tau t)\|_{L^\infty(\mathbb R^d)}
\leq
C_\infty t^{-\frac{d}{s+2}}
\leq
C_\infty T_1^{-\frac{d}{s+2}}.
\]
Hence the \(L^\infty\)-bound required in
Theorem \ref{thm:local-holder-unconfined} is uniform in \(\tau\).
Choosing \(t_0,t_1\) such that $T_1<t_0<1<t_1<T_2,$ and applying Theorem \ref{thm:local-holder-unconfined}, we obtain, at
\(t=1\),
\[
|v_\tau(\xi,1)-v_\tau(\zeta,1)|
\leq
C|\xi-\zeta|^\eta,
\qquad
\xi,\zeta\in B_1,
\]
where \(C>0\) is independent of \(\tau\). Since the unconfined equation and the estimate above are invariant under
spatial translations, the same bound holds, with the same constant, in
every unit ball in \(\mathbb R^d\). Hence, if
\(|\xi-\zeta|\leq1\), we have
\[
|v_\tau(\xi,1)-v_\tau(\zeta,1)|\leq C|\xi-\zeta|^\eta
\]
If \(|\xi-\zeta|>1\), it follows from \(L^\infty\)-bound \eqref{eq:unconfined-Linfty-decay} at \(t=1\) that
\[
|v_\tau(\xi,1)-v_\tau(\zeta,1)|
\leq
|v_\tau(\xi,1)|+|v_\tau(\zeta,1)|
\leq
2C_\infty
\leq
2C_\infty|\xi-\zeta|^\eta.
\]
Thus, we conclude
\begin{equation}\label{inequality:v_holder}
[v_\tau(1)]_{C^\eta(\mathbb R^d)}
\leq C,
\end{equation}
where \(C\) is independent of \(\tau\). Thus, evaluating \eqref{eq:v_rho_holder_scaling} at $t=1$ and invoking \eqref{inequality:v_holder}, we obtain 
\[
[\rho(\tau)]_{C^\eta(\mathbb R^d)}
\leq [v_\tau(1)]_{C^\eta(\mathbb R^d)}
\tau^{-\frac{d+\eta}{s+2}} \leq C_{\eta}
\tau^{-\frac{d+\eta}{s+2}},
\qquad
\tau>0,
\]
\end{proof}

We now transfer this estimate to the confined variables through the scaling
\eqref{eq:scaling-confined-unconfined}, with \(R\) and \(\tau\) defined in
\eqref{eq:scaling-functions-Linfty}.

\begin{proof}[Proof of Theorem \ref{thm:Lp_global_confined}\textup{(iii)}]
Let \(u\) be the limiting curve obtained in
Theorem \ref{thm:global-wellposedness}. Through
\eqref{eq:scaling-confined-unconfined}, we associate with \(u\) the
unconfined curve \(\rho\), so that
\[
u(t,x)=R(t)^d\rho(\tau(t),R(t)x),
\]
where \(R(t)\) and \(\tau(t)\) are defined in
\eqref{eq:scaling-functions-Linfty}.

We first prove the local space-time H\"older regularity of \(u\). Let
\(K\subset\mathbb R^d\) be compact and let \(0<t_0<t_1<\infty\). Since \(R\)
is continuous on \([t_0,t_1]\), the set
\[
\widetilde K
:=
\{R(t)x:x\in K,\ t\in[t_0,t_1]\}
\]
is compact. Moreover, since \(\tau\) is increasing, we have $\tau([t_0,t_1]) =
[\tau(t_0),\tau(t_1)].$ By Theorem \ref{thm:local-holder-unconfined}, there exists
\(\eta\in(0,1)\) and a constant \(C>0\) such that
\[
|\rho(\tau,z)-\rho(\sigma,w)|
\leq
C\left(
|z-w|^\eta+|\tau-\sigma|^\eta
\right)
\]
for every \(z,w\in\widetilde K\) and
\(\tau,\sigma\in[\tau(t_0),\tau(t_1)]\). Therefore, for
\(x,y\in K\) and \(t,s\in[t_0,t_1]\), we have
\[
\begin{aligned}
|u(t,x)-u(s,y)|
&\leq
|R(t)^d-R(s)^d|
|\rho(\tau(t),R(t)x)| +
R(s)^d
|\rho(\tau(t),R(t)x)-\rho(\tau(s),R(s)y)|.
\end{aligned}
\]
By \eqref{eq:unconfined-Linfty-decay}, \(\rho\) is uniformly bounded on
\(\mathbb R^d\times[\tau(t_0),\tau(t_1)]\). Moreover, since \(R\) and
\(\tau\) are smooth on \([t_0,t_1]\), we have
\[
|R(t)^d-R(s)^d|
\leq
C|t-s|,
\qquad
|\tau(t)-\tau(s)|
\leq
C|t-s|,
\]
and
\[
|R(t)x-R(s)y|
\leq
C_{K,t_0,t_1}
\left(
|x-y|+|t-s|
\right).
\]
Hence, we deduce
\[
|u(t,x)-u(s,y)|
\leq
C_{K,t_0,t_1}
\left(
|x-y|^\eta
+
|t-s|^\eta
+
|t-s|
\right).
\]
Since \(0<\eta<1\) and \(t,s\in[t_0,t_1]\), we have
\[
|t-s|
\leq
(t_1-t_0)^{1-\eta}|t-s|^\eta.
\]
Thus, after adjusting the constant,
\[
|u(t,x)-u(s,y)|
\leq
C_{K,t_0,t_1}
\left(
|x-y|^\eta+|t-s|^\eta
\right).
\]
Therefore, \(u\) is locally H\"older continuous in space-time in
\(\mathbb R^d\times(0,\infty)\).

Finally, let us prove the global spatial H\"older bound. For every fixed
\(t>0\), the scaling \eqref{eq:scaling-functions-Linfty} gives
\[
[u(t)]_{C^\eta(\mathbb R^d)}
=
R(t)^{d+\eta}
[\rho(\tau(t))]_{C^\eta(\mathbb R^d)}.
\]
Applying Theorem \ref{thm:unconfined-holder-decay} and \eqref{eq:scaling-functions-Linfty}, we obtain
\begin{equation}\label{eq:confined-holder-explicit}
[u(t)]_{C^\eta(\mathbb R^d)}
\leq
C_\eta e^{2(d+\eta)t}
\left(
\frac{e^{2(s+2)t}-1}{2(s+2)}
\right)^{-\frac{d+\eta}{s+2}},
\qquad
t>0.
\end{equation}
Since the right-hand side of \eqref{eq:confined-holder-explicit} is uniformly bounded for \(t \geq t_0 > 0\), combining this estimate with the \(L^\infty\)-bound from Theorem \ref{thm:Lp_global_confined}\textup{(ii)} yields

\[
\sup_{t\geq t_0}
\|u(t)\|_{C^\eta(\mathbb R^d)}
<\infty, \quad \textrm{for every } t_0 > 0.
\]
Therefore, we conclude the proof of
Theorem \ref{thm:Lp_global_confined}\textup{(iii)}.
\end{proof}

\section{Proof of Theorem \ref{thm:asymptotic-convergence}}
\label{sec:proof-asymptotic-convergence}

In this section, we prove Theorem \ref{thm:asymptotic-convergence}, which
identifies the large-time behavior of the global solution in the quadratic
confinement case \(U(x)=|x|^2\). The energy dissipation inequality gives finite
total dissipation and a uniform second moment bound, hence every sequence of
times admits a narrowly convergent subsequence. The uniform \(L^\infty\) and
H\"older estimates from Theorem \ref{thm:Lp_global_confined}\textup{(ii)--(iii)}
upgrade this compactness to local \(C^\beta\)-compactness. Using time
translations and the vanishing of the dissipation on translated intervals, we
show that every limit point is a steady state. By uniqueness, the only possible
limit is \(\overline{\rho}\), and therefore the whole trajectory converges.

Recall that, for \(s\in\mathcal I_d\), the equilibrium density admits the
explicit ellipsoidal representation
\eqref{eq:explicit-equilibrium-profile}, namely
\begin{equation*}
\overline{\rho}(x)
=
C\left(1-\frac{|Bx|^2}{R^2}\right)_+^{\alpha/2},
\qquad
\alpha=s-d+2\in(0,2),
\end{equation*}
for some invertible matrix \(B\), radius \(R>0\), and constant \(C>0\); see
\cite{CS24cpam,Carrillo_Shu_3D,
Mateu_Mora_Rondi_Scardia_Verdera_2023,frank2026explicit}.

\begin{remark}[H\"older regularity of the equilibrium density]
\label{rem:holder-regularity-equilibrium-density}
Assume that \(s\in\mathcal I_d\). By
\eqref{eq:explicit-equilibrium-profile}, we may write
\(\overline{\rho}=CG^{\alpha/2}\), where
\(G(x):=(1-|Ax|^2/R^2)_+\). Since the positive-part map is Lipschitz and
\(x\mapsto1-|Ax|^2/R^2\) has bounded gradient on a neighborhood of the compact
ellipsoid \(\{|Ax|\leq R\}\), the function \(G\) is globally Lipschitz. Therefore, for every \(x,y\in\mathbb R^d\),
\[
|\overline{\rho}(x)-\overline{\rho}(y)|
=
C|G(x)^{\alpha/2}-G(y)^{\alpha/2}|
\leq
C|G(x)-G(y)|^{\alpha/2}
\leq
C|x-y|^{\alpha/2}.
\]
where we used that $|a^{\alpha/2}-b^{\alpha/2}|
\leq |a-b|^{\alpha/2}$ for every $a,b\geq0$ and \(0<\alpha/2<1\).  Hence, \(\overline{\rho}\in C^{\alpha/2}(\mathbb R^d)\). In particular,
\(\overline{\rho}\in C^\gamma(\mathbb R^d)\) for every
\(0<\gamma\leq\alpha/2\).
\end{remark}

Recall that, for \(U(x)=|x|^2\), we write
\begin{equation}\label{eq:energy-asymptotic}
\mathcal E(u)
=
\frac12\int_{\mathbb R^d}(W*u)(x)u(x)dx
+
\int_{\mathbb R^d}|x|^2u(x)dx,
\end{equation}
and \(V[u](x):=(W*u)(x)+|x|^2\). {For quadratic confinement, we write
\[
\mathcal V(u):=\mathcal V_U(u\,dx),
\]
where \(\mathcal V_U\) is defined in \eqref{eq:defV}.
Thus, a continuous probability density \(u\) is a steady
state precisely when \(\mathcal V(u)=0\).}

Finally, we will need the following lower semi-continuity result:
\begin{lemma}[Lower semicontinuity under locally uniform convergence]
\label{lem:dissipation-lsc-uniform}
Under the standing assumptions on \(W\) and \(U\), let
\(\rho_n,\rho\) be continuous probability densities on \(\mathbb R^d\)
such that
\[
\rho_n\longrightarrow\rho
\qquad\text{locally uniformly in }\mathbb R^d.
\]
Then
\[
\mathcal V_U(\rho)
\leq
\liminf_{n\to\infty}\mathcal V_U(\rho_n).
\]
\end{lemma}

\begin{proof}
Since \(W\) is locally integrable and vanishes at infinity,
local uniform convergence and unit mass imply
\[
V_n:=W*\rho_n+U
\longrightarrow
V:=W*\rho+U
\qquad\text{locally uniformly}.
\]
Indeed, split each convolution into a bounded region, where local
uniform convergence applies, and its complement, whose contribution
is uniformly small by the decay of \(W\) and the mass bound.

There is nothing to prove if the lower limit of the dissipations
is infinite. Otherwise, pass to a subsequence such that
\[
\mathcal V_U(\rho_n)\longrightarrow
\ell:=\liminf_{n\to\infty}\mathcal V_U(\rho_n)<\infty.
\]
For \(j\geq1\), set
\[
O_j:=B_j(0)\cap\{\rho>1/j\}.
\]
For all sufficiently large \(n\), local uniform convergence gives
\(\rho_n\geq1/(2j)\) on \(O_j\), and hence
\[
\int_{O_j}|\nabla V_n|^2dx
\leq 2j\mathcal V_U(\rho_n).
\]
Together with the uniform convergence of \(V_n\) on \(O_j\), this
implies that \(V\in H^1(O_j)\) and
\[
\nabla V_n\rightharpoonup\nabla V
\qquad\text{weakly in }L^2(O_j;\mathbb R^d).
\]
Since \(\sqrt{\rho_n}\to\sqrt{\rho}\) uniformly on \(O_j\), it follows
that
\[
\sqrt{\rho_n}\nabla V_n
\rightharpoonup
\sqrt{\rho}\nabla V
\qquad\text{weakly in }L^2(O_j;\mathbb R^d).
\]
Weak lower semicontinuity of the $L^2$ norm w.r.t. weak $L^2$ convergence therefore yields
\[
\int_{O_j}\rho|\nabla V|^2dx
\leq
\liminf_{n\to\infty}
\int_{O_j}\rho_n|\nabla V_n|^2dx
\leq \ell.
\]
Finally, \(O_j\uparrow\{\rho>0\}\), so
\(V\in H^1_{\mathrm{loc}}(\{\rho>0\})\), and monotone convergence gives
\[
\mathcal V_U(\rho)\leq\ell.
\]
\end{proof}

\begin{proof}[Proof of Theorem \ref{thm:asymptotic-convergence}]

\medskip
\noindent\textbf{Step 1:}
By the energy dissipation inequality in
Theorem \ref{thm:global-wellposedness}, for every \(T>0\),
\begin{equation}\label{eq:energy-dissipation-asymptotic}
\mathcal E(u(T))
+
\int_0^T \mathcal V(u(t))dt
\leq
\mathcal E(u_0).
\end{equation}
Since the interaction energy and the confinement term are nonnegative in the
quadratic confinement case, we obtain
\[
\int_0^\infty \mathcal V(u(t))dt<\infty
\]
and
\begin{equation}\label{eq:uniform-second-moment-asymptotic}
\sup_{t\geq0}
\int_{\mathbb R^d}|x|^2u(t,x)dx
\leq
\mathcal E(u_0).
\end{equation}
In particular, for every fixed \(T>0\) and every sequence \(t_n\to\infty\),
\begin{equation}\label{eq:dissipation-vanishes-on-translates}
\int_{t_n}^{t_n+T}\mathcal V(u(t))dt\longrightarrow0.
\end{equation}

The bound \eqref{eq:uniform-second-moment-asymptotic} implies tightness of
\(\{u(t):t\geq0\}\) in \(\mathcal P(\mathbb R^d)\). Moreover, by
Theorem \ref{thm:global-wellposedness}, we have that $u\in AC^2\bigl([0,\infty);(\mathcal P_2(\mathbb R^d),W_2)\bigr).$ In particular, \(u\) is \(W_2\)-continuous, hence narrowly continuous, and then  \(u(t)\) are well defined for every \(t\geq0\). Therefore, for any
sequence \(t_n\to\infty\), there exists a subsequence, still denoted by \(t_n\),
and a measure \(\mu\in\mathcal P(\mathbb R^d)\) such that
\begin{equation}\label{eq:subsequence-narrow-limit-asymptotic}
u(t_n)\rightharpoonup\mu
\qquad
\text{narrowly in }\mathcal P(\mathbb R^d).
\end{equation}
Moreover, by the lower semicontinuity of the second moment and
\eqref{eq:uniform-second-moment-asymptotic}, we obtain
\[
\int_{\mathbb R^d}|x|^2d\mu(x)
\leq
\liminf_{n\to\infty}
\int_{\mathbb R^d}|x|^2u(t_n,x)dx
\leq
\mathcal E(u_0).
\]
Thus \(\mu\in\mathcal P_2(\mathbb R^d)\).

Fix \(t_0>0\). Since \(t_n\to\infty\), we may assume that \(t_n\geq t_0\) for
every \(n\). By Theorem \ref{thm:Lp_global_confined}\textup{(ii)--(iii)},
\begin{equation}\label{eq:uniform-Linfty-holder-asymptotic}
\sup_{n\geq1}\|u(t_n)\|_{L^\infty(\mathbb R^d)}
+
\sup_{n\geq1}[u(t_n)]_{C^\eta(\mathbb R^d)}
\leq
C(t_0).
\end{equation}
By the compact embedding \(C^\eta(B_R)\hookrightarrow C^\beta(B_R)\), for every
\(0<\beta<\eta\), and by a diagonal argument on balls, the sequence
\(\{u(t_n)\}_{n\geq1}\) is relatively compact in
\(C^\beta_{\mathrm{loc}}(\mathbb R^d)\). On the other hand, it already converges
narrowly to \(\mu\) by \eqref{eq:subsequence-narrow-limit-asymptotic}. Hence any
local \(C^\beta\)-limit obtained from a subsequence must coincide, in the sense
of distributions, with the measure \(\mu\). Therefore \(\mu\) is represented by
a bounded continuous density and, identifying \(\mu\) with this density. Thus, we obtain that
\[
u(t_n)\to\mu
\qquad
\text{strongly in }C^\beta_{\mathrm{loc}}(\mathbb R^d)
\]
for every \(0<\beta<\eta\).

\medskip
\noindent\textbf{Step 2:}
Fix $T>0$ and, for each $n$, define the curve $u_n:[0,T] \to \mathcal{P}_2(\mathbb{R}^d)$ by $u_n(\tau) := u(t_n+\tau)$. Since $u \in AC^2([0,\infty); (\mathcal{P}_2(\mathbb{R}^d), W_2))$, each $u_n$ is well-defined. We claim that $u_n$ converges to the constant curve $\mu$.

Let \(\chi\in C_c^\infty(\mathbb R^d)\). From the weak formulation in
Theorem \ref{thm:global-wellposedness}, applied on the time interval
\([t_n,t_n+\tau]\), it follows for every \(\tau\in[0,T]\) that
\[
\int_{\mathbb R^d}\chi(x)u(t_n+\tau,x)dx
-
\int_{\mathbb R^d}\chi(x)u(t_n,x)dx
=
-\int_{t_n}^{t_n+\tau}
\int_{\mathbb R^d}
u(t,x)\nabla V[u(t)](x)\cdot\nabla\chi(x)dxdt.
\]
By the Cauchy--Schwarz inequality, we get
\[
\left|
\int_{\mathbb R^d}\chi(x)
\bigl(u(t_n+\tau,x)-u(t_n,x)\bigr)dx
\right|
\leq
\|\nabla\chi\|_{L^\infty}T^{1/2}
\left(
\int_{t_n}^{t_n+T}\mathcal V(u(t))dt
\right)^{1/2}.
\]
Using \eqref{eq:dissipation-vanishes-on-translates}, the right-hand side tends
to zero uniformly for \(\tau\in[0,T]\). Since
\(u(t_n)\rightharpoonup\mu\) narrowly, we obtain, for every
\(\chi\in C_c^\infty(\mathbb R^d)\),
\begin{equation}\label{eq:limit_assymptotic}
\sup_{\tau\in[0,T]}
\left|
\int_{\mathbb R^d}\chi(x)u(t_n+\tau,x)dx
-
\int_{\mathbb R^d}\chi(x)d\mu(x)
\right|
\longrightarrow0.
\end{equation}

We now extend this convergence to bounded continuous test functions. Let
\(\chi\in C_b(\mathbb R^d)\), and let \(\zeta_R\in C_c^\infty(\mathbb R^d)\)
be such that \(0\leq\zeta_R\leq1\), \(\zeta_R\equiv1\) in \(B_R\), and
\(\zeta_R\equiv0\) in \(\mathbb R^d\setminus B_{2R}\). Set
\(\chi_R:=\chi\zeta_R\). By \eqref{eq:uniform-second-moment-asymptotic}, we have
\[
\sup_{n\geq1}\sup_{\tau\in[0,T]}
\int_{\{|x|>R\}}u(t_n+\tau,x)dx
\leq
\frac{C}{R^2},
\]
and, since \(\mu\in\mathcal P_2(\mathbb R^d)\), it follows that $\mu(\{|x|>R\})\leq \frac{C}{R^2}.$ Thus, we obtain
\begin{align}
\sup_{\tau\in[0,T]}
\left|
\int_{\mathbb R^d}\chi u(t_n+\tau)dx
-
\int_{\mathbb R^d}\chi d \mu
\right| \notag
&\leq
\sup_{\tau\in[0,T]}
\left|
\int_{\mathbb R^d}\chi_R u(t_n+\tau)dx
-
\int_{\mathbb R^d}\chi_Rd\mu
\right|
\\ \notag
&
+
\|\chi\|_{L^\infty}
\left(
\sup_{\tau\in[0,T]}
\int_{|x|>R}u(t_n+\tau,x)dx
+
\mu(\mathbb R^d\setminus B_R)
\right)
\\ \notag
&\leq
\sup_{\tau\in[0,T]}
\left|
\int_{\mathbb R^d}\chi_R u(t_n+\tau)dx
-
\int_{\mathbb R^d}\chi_Rd\mu
\right|
+
\frac{C\|\chi\|_{L^\infty}}{R^2}.
\end{align}
For fixed \(R>0\), we have \(\chi_R\in C_c(\mathbb R^d)\). By uniformly
approximating \(\chi_R\) with functions in \(C_c^\infty(\mathbb R^d)\), and
using that \(u(t_n+\tau)\) and \(\mu\) are probability measures, the
convergence in \eqref{eq:limit_assymptotic} also holds for \(\chi_R\).
Therefore, for fixed \(R\), the term above on the right-hand side tends to zero as \(n\to\infty\).
Letting \(R\to\infty\), we obtain
\[
\sup_{\tau\in[0,T]}
\left|
\int_{\mathbb R^d}\chi(x)u(t_n+\tau,x)dx
-
\int_{\mathbb R^d}\chi(x)d\mu(x)
\right|
\longrightarrow0
\]
for every \(\chi\in C_b(\mathbb R^d)\). Recalling that
\(u_n(\tau):=u(t_n+\tau)\), we conclude that \(u_n\) converges
uniformly in \(\tau\in[0,T]\) to the measure \(\mu\), which is independent
of \(\tau\), with respect to the narrow topology. Thus, we have
\[
u_n\longrightarrow\mu
\qquad
\text{in }
C\bigl([0,T];\mathcal P(\mathbb R^d)\bigr).
\]
In particular, we have
\[
u_n(\tau,x)dxd\tau
\rightharpoonup
d\mu(x)d\tau
\qquad
\text{narrowly in }[0,T]\times\mathbb R^d.
\]

\medskip
\noindent\textbf{Step 3:} For every fixed \(\tau\in[0,T]\), the uniform spatial
H\"older bounds from
Theorem \ref{thm:Lp_global_confined}\textup{(iii)}
imply that \(\{u_n(\tau)\}_n\) is relatively compact
in \(C_{\mathrm{loc}}(\mathbb R^d)\).
Since Step 2 gives
\(u_n(\tau)\rightharpoonup\mu\) narrowly for every
\(\tau\in[0,T]\), and \(\mu\) has a continuous density
by Step 1, every locally uniform subsequential limit
coincides with \(\mu\). Hence,
\[
u_n(\tau)\longrightarrow\mu
\qquad
\text{locally uniformly in }\mathbb R^d
\]
for every fixed \(\tau\in[0,T]\). Applying Lemma \ref{lem:dissipation-lsc-uniform}
for each fixed \(\tau\), and then Fatou's lemma,
we obtain
\[
\begin{aligned}
\int_0^T\mathcal V(\mu)d\tau
&\leq
\int_0^T\liminf_{n\to\infty}
\mathcal V(u_n(\tau))d\tau\\
&\leq
\liminf_{n\to\infty}
\int_0^T\mathcal V(u_n(\tau))d\tau\\
&=
\liminf_{n\to\infty}
\int_{t_n}^{t_n+T}\mathcal V(u(t))dt
=0.
\end{aligned}
\]
Hence \(\mathcal V(\mu)=0\), so \(\mu\) is a steady state in the sense of
\eqref{eq:defV}.  \(
\mu\in C(\mathbb R^d)\cap L^\infty(\mathbb R^d)
\cap\mathcal P(\mathbb R^d),
\), Theorem \ref{thm:uniqueness} yields
\[
\mu=\overline{\rho}.
\]

Recall that the sequence \(t_n\to\infty\) was arbitrary. Thus, every
sequence of times tending to infinity admits, by the compactness argument
of Step 1, a convergent subsequence, and the argument above together with
Theorem \ref{thm:uniqueness} shows that the limit of any such subsequence
must be the same and given by \(\overline{\rho}\). Therefore, the function $u(t)$ satisfies
\[
u(t)\rightharpoonup\overline{\rho}
\qquad
\text{narrowly in }\mathcal P(\mathbb R^d)
\quad
\text{as }t\to\infty,
\]
Moreover, the same argument together with the relative compactness in
\(C^\beta_{\mathrm{loc}}(\mathbb R^d)\) established in Step 1 shows that,
for every compact set \(K\subset\mathbb R^d\) and every \(0<\beta<\eta\), we have
\begin{equation}\label{eq:converge_Holder}
u(t)\longrightarrow\overline{\rho}
\qquad
\text{strongly in }C^\beta(K)
\quad
\text{as }t\to\infty.
\end{equation}

\medskip
\noindent\textbf{Step 4:} Finally, we extend the local convergence to global $C^\beta$-convergence. We first prove the uniform $L^\infty$-decay of the spatial tail. Namely, for 
every $t_0 > 0$, it holds that

\begin{equation}\label{eq:uniform-tail-Linfty-asymptotic}
\lim_{R\to\infty}
\sup_{t\geq t_0}
\|u(t)\|_{L^\infty(\mathbb R^d\setminus B_R)}
=0.
\end{equation}
Suppose that \eqref{eq:uniform-tail-Linfty-asymptotic} fails. Then there exist
\(\varepsilon>0\), a sequence \(t_n\geq t_0\), and points
\(x_n\in\mathbb R^d\), with \(|x_n|\to\infty\), such that
\[
u(t_n,x_n)\geq\varepsilon.
\]
By Theorem \ref{thm:Lp_global_confined}\textup{(iii)}, there exists
\(H=H(t_0)>0\) such that
\[
\sup_{t\geq t_0}[u(t)]_{C^\eta(\mathbb R^d)}\leq H.
\]
Set \(r_\varepsilon:=(\varepsilon/(2H))^{1/\eta}\). Then, for every
\(x\in B_{r_\varepsilon}(x_n)\), it follows that
\begin{equation}\label{ineq:u_varepsilon}
u(t_n,x)
\geq
u(t_n,x_n)-H|x-x_n|^\eta
\geq
\frac{\varepsilon}{2}.
\end{equation}
Since \(r_\varepsilon\) is independent of \(n\) and \(|x_n|\to\infty\), for all
sufficiently large \(n\) we have \(r_\varepsilon\leq |x_n|/2\). Thus, for such
\(n\) and every \(x\in B_{r_\varepsilon}(x_n)\), we get
\[
|x|
\geq
|x_n|-|x-x_n|
\geq
|x_n|-r_\varepsilon
\geq
\frac{|x_n|}{2}.
\]
Thus, we obtain from \eqref{ineq:u_varepsilon} that
\[
\int_{\mathbb R^d}|x|^2u(t_n,x)dx
\geq
\int_{B_{r_\varepsilon}(x_n)}|x|^2u(t_n,x)dx
\geq
\frac{\varepsilon |B_{r_\varepsilon}|}{8}|x_n|^2.
\]
This contradicts \eqref{eq:uniform-second-moment-asymptotic}, since
\(|x_n|\to\infty\), and hence \eqref{eq:uniform-tail-Linfty-asymptotic} holds.

Since \(\overline{\rho}\) is compactly supported, choose \(R_0>0\) such that
\(\operatorname{supp}\overline{\rho}\subset B_{R_0}\). Then, for every
\(R\geq R_0\), we have
\[
\|u(t)-\overline{\rho}\|_{L^\infty(\mathbb R^d)}
\leq
\|u(t)-\overline{\rho}\|_{L^\infty(B_R)}
+
\|u(t)\|_{L^\infty(\mathbb R^d\setminus B_R)}.
\]
Fix \(\varepsilon>0\). By \eqref{eq:uniform-tail-Linfty-asymptotic}, we can
choose \(R\geq R_0\) sufficiently large so that the second term on the right-hand side is
smaller than \(\varepsilon/2\), uniformly for \(t\geq t_0\). For this fixed
\(R\), the local convergence \eqref{eq:converge_Holder} implies that the
first term is also smaller than \(\varepsilon/2\) for all sufficiently large
\(t\). Therefore, we obtain
\[
\|u(t)-\overline{\rho}\|_{L^\infty(\mathbb R^d)}
\longrightarrow0
\qquad
\text{as }t\to\infty.
\]
Finally, set \(\eta_*:=\min\{\eta,\alpha/2\}\). By
Remark \ref{rem:holder-regularity-equilibrium-density} and
Theorem \ref{thm:Lp_global_confined}\textup{(ii)--(iii)}, we have
\[
\sup_{t\geq t_0}
\|u(t)-\overline{\rho}\|_{C^{\eta_*}(\mathbb R^d)}
<\infty.
\]
Therefore, for every \(0<\beta<\eta_*\), the interpolation inequality gives
\[
[u(t)-\overline{\rho}]_{C^\beta(\mathbb R^d)}
\leq
C
\|u(t)-\overline{\rho}\|_{L^\infty(\mathbb R^d)}^{1-\beta/\eta_*}
[u(t)-\overline{\rho}]_{C^{\eta_*}(\mathbb R^d)}^{\beta/\eta_*}.
\]
Since the \(L^\infty\)-norm tends to zero and the \(C^{\eta_*}\)-seminorm is
uniformly bounded, we conclude that
\[
\|u(t)-\overline{\rho}\|_{C^\beta(\mathbb R^d)}
\longrightarrow0
\qquad
\text{as }t\to\infty
\]
for every \(0<\beta<\eta_*\). Therefore, we conclude the proof.
\end{proof}

\section{Proof of Theorem \ref{thm:global-wellposedness}}
\label{sec:proof-gradient-flow}

To adapt the variational framework of \cite[Theorem 1.1]{foghem2025gradient} to our confined anisotropic flow, we employ the flow-interchange argument. Let us recall the class of auxiliary functionals introduced in \cite{foghem2025gradient}. Let
\[
V \in C^1((0,\infty)) \cap C([0,\infty)) \qquad\text{with}\qquad V(0)=0
\]
be a convex function satisfying the McCann displacement convexity condition. The associated displacement convex entropy is the functional $\mathcal{V}: \mathcal{P}_2(\mathbb{R}^d) \to \mathbb{R} \cup \{+\infty\}$, defined by
\begin{equation} \label{function_V}
\mathcal{V}(u)
=
\begin{cases}
\displaystyle \int_{\mathbb{R}^d} V(u(x))  dx, & \text{if } du(x) = u(x)  dx, \\[2mm]
+\infty, & \text{otherwise}.
\end{cases}
\end{equation}
We denote its effective domain by $D(\mathcal{V}) := \{u \in \mathcal{P}_2(\mathbb{R}^d) : \mathcal{V}(u) < \infty\}$.

Let $(S_t)_{t \ge 0}$ denote the auxiliary flow associated with $\mathcal{V}$, that is, the $2$-Wasserstein gradient flow generated by $\mathcal{V}$. Thus, for a given initial datum $u$, the quantity $S_t u$ represents the state at time $t$ of the auxiliary evolution starting from $u$. For another functional $\mathcal{H}: \mathcal{P}_2(\mathbb{R}^d) \to (-\infty, +\infty]$, we define the dissipation of $\mathcal{H}$ along the flow generated by $\mathcal{V}$ by
\begin{equation*}
D_{\mathcal{V}}[\mathcal{H}](u)
:=
\limsup_{t \to 0} \frac{\mathcal{H}(u) - \mathcal{H}(S_t u)}{t}.
\end{equation*}

In the applications of the flow-interchange method below, two particular choices of $\mathcal{V}$ (as defined in \eqref{function_V}) will be used repeatedly. First, when
\[
V(r) = r \log r,
\]
the associated functional is the Boltzmann entropy (defined in \eqref{Boltzmann_defi}), and the corresponding auxiliary flow is the heat semigroup. Second, for $1 < p < \infty$, when
\[
V(r) = \frac{r^p}{p-1},
\]
the associated functional becomes
\begin{equation*}
\mathcal{G}_p(u) := \frac{1}{p-1} \int_{\mathbb{R}^d} u(x)^p  dx,
\end{equation*}
whose formal $2$-Wasserstein gradient flow is the porous-medium equation
\[
\partial_t u = \Delta(u^p).
\]

Before proving the theorem, we need a brief lemma:

\begin{lemma}[Structural properties of the kernels and symbols]
\label{lem:checks_fpgs} Fix \(\alpha=s-d+2\in(0,2)\). Assume the standing hypotheses
\eqref{hw}, \eqref{hl}, and \eqref{ha} . Then the following
properties hold.

\begin{enumerate}
\item[\rm (i)] The kernel \(K_L\) satisfies the L\'evy integrability
condition and is nonintegrable at the origin:
\[
\int_{\mathbb R^d}(1\wedge |h|^2)K_L(h)dh<\infty,
\qquad
\int_{B_1}K_L(h)dh=\infty.
\]

\item[\rm (ii)] The kernel \(K_A\) is a unimodal L\'evy kernel and has the
exact scaling
\[
K_A(\lambda h)=\lambda^{-d-\alpha}K_A(h)
\qquad
\text{for every }\lambda>0\text{ and }h\neq0.
\]
More precisely,
\[
K_A(x)
\leq
\frac{\Lambda_A}{\lambda_A}K_A(y)
\qquad
\text{whenever }|x|\geq |y|>0.
\]

\item[\rm (iii)] The symbol \(\psi\) satisfies
\begin{equation*}
\psi(\xi)
\geq
c_\psi(1\wedge |\xi|^2)
\qquad
\text{for every }\xi\in\mathbb R^d,
\end{equation*}
where, with the natural continuous extension \(\psi(0)=0\), one may take
\[
c_\psi
:=
\|\widetilde\Omega\|_{L^\infty(\mathbb S^{d-1})}^{-1}.
\]

\item[\rm (iv)] The symbol \(\widetilde\psi\) satisfies the quantitative
compactness estimate
\begin{equation}\label{eq:tildepsi-translation-bound}
\sup_{\xi\in\mathbb R^d\setminus\{0\}}
\frac{|\mathrm e^{i\xi\cdot h}-1|^2}{\widetilde\psi(\xi)}
\leq
\frac{2^{2-\alpha}}
{\displaystyle\min_{\theta\in\mathbb S^{d-1}}\widetilde\Omega(\theta)}
|h|^\alpha
\qquad
\text{for every }h\in\mathbb R^d.
\end{equation}
In particular,
\[
\sup_{\xi\in\mathbb R^d\setminus\{0\}}
\frac{|\mathrm e^{i\xi\cdot h}-1|^2}{\widetilde\psi(\xi)}
\longrightarrow0
\qquad\text{as }|h|\to0.
\]
\end{enumerate}
\end{lemma}

\begin{proof}
Recall that
\[
K_L(h)=|h|^{s-2d}\omega_L\!\left(\frac{h}{|h|}\right),
\qquad
K_A(h)=|h|^{-d-\alpha}\omega_A\!\left(\frac{h}{|h|}\right).
\]
Item \textup{(i)} follows from a standard computation in polar coordinates, while the evenness and positivity of \(K_L\) follow from \eqref{hl}. Similarly, polar coordinates give \((1\wedge|h|^2)K_A(h)\in L^1(\mathbb R^d)\), and the evenness and positivity of \(K_A\) follow from \eqref{ha} . Moreover, if \(|x|\geq|y|>0\), then
\[
K_A(x)\leq\Lambda_A|x|^{-d-\alpha}
\leq\Lambda_A|y|^{-d-\alpha}
\leq\frac{\Lambda_A}{\lambda_A}K_A(y).
\]
The homogeneity relation \(K_A(\lambda h)=\lambda^{-d-\alpha}K_A(h)\) then completes the proof of item \textup{(ii)}.

For item \textup{(iii)}, set \(\beta:=d-s\in(0,2)\). Since
\[
\psi(\xi)=|\xi|^\beta
\widetilde\Omega\!\left(\frac{\xi}{|\xi|}\right)^{-1},
\]
we have
\[
\psi(\xi)\geq
\|\widetilde\Omega\|_{L^\infty(\mathbb S^{d-1})}^{-1}|\xi|^\beta
\geq
\|\widetilde\Omega\|_{L^\infty(\mathbb S^{d-1})}^{-1}
(1\wedge|\xi|^2),
\]
where the last inequality follows from \(\beta<2\) when \(|\xi|\leq1\) and from \(\beta>0\) when \(|\xi|\geq1\).

Finally, since \(\widetilde\psi\) is the symbol of \(A\),
\[
\widetilde\psi(\xi)
=
|\xi|^\alpha
\widetilde\Omega\!\left(\frac{\xi}{|\xi|}\right).
\]
Let \(\widetilde\Omega_{\min}:=\min_{\theta\in\mathbb S^{d-1}}\widetilde\Omega(\theta)>0\). Using \(|\mathrm e^{it}-1|\leq\min\{2,|t|\}\), for \(h\neq0\) we obtain
\[
\frac{|\mathrm e^{i\xi\cdot h}-1|^2}{\widetilde\psi(\xi)}
\leq
\frac{\min\{4,|\xi|^2|h|^2\}}
{\widetilde\Omega_{\min}|\xi|^\alpha}
=
\frac{|h|^\alpha}{\widetilde\Omega_{\min}}
\frac{\min\{4,r^2\}}{r^\alpha},
\qquad r:=|\xi||h|.
\]
The last quotient is bounded by \(2^{2-\alpha}\), since it equals \(r^{2-\alpha}\) for \(0<r\leq2\) and \(4r^{-\alpha}\) for \(r\geq2\). Taking the supremum over \(\xi\neq0\) proves \eqref{eq:tildepsi-translation-bound}. The case \(h=0\) is immediate, and the stated limit follows from \(\alpha>0\).
\end{proof}

Next, we prove the main result, Theorem \ref{thm:global-wellposedness}.

\begin{proof}[Proof of item \textup{(i)}]

For $\nu\in\mathcal P_2(\mathbb R^d)$, set
\[
F(\nu)
:=
\frac{1}{2(2\pi)^d}
\int_{\mathbb R^d}
|\widehat\nu(\xi)|^2\psi(\xi)^{-1}d\xi,
\]
with the convention that $F(\nu)=+\infty$ if the integral diverges. Thus
\[
\mathcal E_U(\nu)
=
F(\nu)+\int_{\mathbb R^d}Ud\nu.
\]
Fix $\tau>0$. Suppose that $\overline u
\in
\dot H_{\psi^{-1}}(\mathbb R^d)
\cap
\mathcal P_2(\mathbb R^d)$ and define
\[
\mathcal J_{\tau,\overline u}(\nu)
:=
\frac1{2\tau}W_2^2(\nu,\overline u)
+
F(\nu)
+
\int_{\mathbb R^d}Ud\nu.
\]

We first prove that $\mathcal J_{\tau,\overline u}$ admits a unique
minimizer. Since $U(x)\le C_U|x|^2+A_U,$ we have
\[
\int_{\mathbb R^d}Ud\overline u<\infty,
\]
and therefore
\[
\mathcal J_{\tau,\overline u}(\overline u)
=
F(\overline u)+\int_{\mathbb R^d}Ud\overline u
<\infty.
\]
Hence the functional is proper.

Let $(\nu_j)_j\subset\mathcal P_2(\mathbb R^d)$ be a minimizing
sequence. Since $F\ge0$, the Wasserstein term is nonnegative, and
\[
U(x)\ge c_U|x|^2-a_U,
\]
there exists $M>0$ such that
\[
c_U\int_{\mathbb R^d}|x|^2d\nu_j(x)-a_U
\le
\mathcal J_{\tau,\overline u}(\nu_j)
\le M
\]
for every $j$. Thus $(\nu_j)_j$ has uniformly bounded second moments.
After passing to a subsequence, there exists
$\nu\in\mathcal P_2(\mathbb R^d)$ such that
\[
\nu_j\rightharpoonup\nu
\qquad\text{narrowly}.
\]

The Wasserstein term is lower semicontinuous under narrow convergence.
Since $U$ is continuous and bounded from below, the map
\[
\nu\longmapsto\int_{\mathbb R^d}Ud\nu
\]
is also lower semicontinuous. Moreover, narrow convergence implies
\[
\widehat\nu_j(\xi)\longrightarrow\widehat\nu(\xi)
\qquad\text{for every }\xi\in\mathbb R^d,
\]
and Fatou's lemma gives
\[
F(\nu)\le\liminf_{j\to\infty}F(\nu_j).
\]
Consequently,
\[
\mathcal J_{\tau,\overline u}(\nu)
\le
\liminf_{j\to\infty}
\mathcal J_{\tau,\overline u}(\nu_j),
\]
so $\nu$ is a minimizer. Since the minimum is finite and
$U\ge-a_U$, we have $F(\nu)<\infty$. Hence
\[
\nu
\in
\dot H_{\psi^{-1}}(\mathbb R^d)
\cap
\mathcal P_2(\mathbb R^d).
\]

To prove uniqueness, let $\nu_0,\nu_1$ belong to the effective domain
of $\mathcal J_{\tau,\overline u}$ and set
\[
\nu_\theta:=(1-\theta)\nu_0+\theta\nu_1,
\qquad 0<\theta<1.
\]
The map
\[
\nu\longmapsto W_2^2(\nu,\overline u)
\]
is convex under linear interpolation, while
\[
\nu\longmapsto\int_{\mathbb R^d}Ud\nu
\]
is linear. Furthermore,
\begin{align*}
F(\nu_\theta)
&=
(1-\theta)F(\nu_0)+\theta F(\nu_1)
-\frac{\theta(1-\theta)}{2(2\pi)^d}
\int_{\mathbb R^d}
|\widehat\nu_0(\xi)-\widehat\nu_1(\xi)|^2
\psi(\xi)^{-1}d\xi.
\end{align*}
Since $\psi^{-1}>0$ on $\mathbb R^d\setminus\{0\}$, the last integral
vanishes only if $\nu_0=\nu_1$. Therefore
$\mathcal J_{\tau,\overline u}$ is strictly convex, and its minimizer
is unique.

Starting with $u_\tau^0=u_0$, we may apply the preceding argument
recursively with $\overline u=u_\tau^{k-1}$. This gives a unique
sequence $(u_\tau^k)_{k\ge0}$ satisfying
\[
u_\tau^k
\in
\operatorname*{argmin}_{\nu\in\mathcal P_2(\mathbb R^d)}
\left\{
\frac1{2\tau}W_2^2(\nu,u_\tau^{k-1})
+
\mathcal E_U(\nu)
\right\},
\]
and
\[
u_\tau^k
\in
\dot H_{\psi^{-1}}(\mathbb R^d)
\cap
\mathcal P_2(\mathbb R^d).
\]

It remains to prove the additional regularity. Fix $k\ge1$ and write
\[
u:=u_\tau^k,
\qquad
\overline u:=u_\tau^{k-1}.
\]
For $t>0$, let 
\[
u_t:=G_t*u, \qquad \textrm{where} \quad 
G_t(x):=\frac{1}{(4\pi t)^{d/2}}e^{-\frac{|x|^2}{4t}}.
\]

Since \(\Delta U\leq\Lambda_U\), the heat semigroup satisfies
\[
G_t*U-U
=
\int_0^tG_r*(\Delta U)dr
\leq
\Lambda_U t.
\]
Therefore,
\begin{equation}\label{eq:heat-potential-simple}
\int_{\mathbb R^d}Udu_t
=
\int_{\mathbb R^d}(G_t*U)du
\leq
\int_{\mathbb R^d}Udu+\Lambda_U t.
\end{equation}

Let $\gamma$ be an optimal transport plan between $u$ and
$\overline u$. Adding an independent random variable with density
$G_t$ to the first component of $\gamma$ gives a coupling between
$u_t$ and $\overline u$. Since $G_t$ has zero mean and second moment
$2dt$,
\begin{equation}\label{eq:heat-wasserstein-simple}
W_2^2(u_t,\overline u)
\le
W_2^2(u,\overline u)+2dt.
\end{equation}

Using $u_t$ as a competitor in the minimization problem for $u$, and
using \eqref{eq:heat-potential-simple} and
\eqref{eq:heat-wasserstein-simple}, we obtain
\[
F(u)-F(u_t)
\le
\left(\frac d\tau+\Lambda_U\right)t.
\]
On the Fourier side,
\[
\frac{F(u)-F(u_t)}{t}
=
\frac{1}{(2\pi)^d}
\int_{\mathbb R^d}
|\widehat u(\xi)|^2\psi(\xi)^{-1}
\frac{1-e^{-2t|\xi|^2}}{2t}d\xi.
\]
The integrand is nonnegative and converges pointwise, as
$t\downarrow0$, to
\[
|\widehat u(\xi)|^2|\xi|^2\psi(\xi)^{-1}.
\]
Fatou's lemma therefore gives
\[
\|u\|_{\dot H_{\widetilde\psi}(\mathbb R^d)}^2
\leq
\liminf_{t\downarrow0}
\frac{F(u)-F(u_t)}{t}
\leq
\frac d\tau+\Lambda_U.
\]
Since
\[
\widetilde\psi(\xi)=|\xi|^2\psi(\xi)^{-1},
\]
we conclude that
\[
u_\tau^k\in\dot H_{\widetilde\psi}(\mathbb R^d)
\qquad\text{for every }k\ge1.
\]

Moreover, since \(u_\tau^k\) is a probability measure,
\(|\widehat{u_\tau^k}|\leq1\), while
\(\widetilde\psi(\xi)\geq \widetilde\Omega_{\min}>0\) for
\(|\xi|\geq1\). Hence, splitting the Fourier integral into
\(\{|\xi|\leq1\}\) and \(\{|\xi|>1\}\), we obtain
\[
\|u_\tau^k\|_{L^2(\mathbb R^d)}^2
\leq
|B_1|
+
\widetilde\Omega_{\min}^{-1}
\|u_\tau^k\|_{\dot H_{\widetilde\psi}(\mathbb R^d)}^2<\infty,
\]
and therefore \(u_\tau^k\in H_{\widetilde\psi}(\mathbb R^d)\).

\end{proof}

\vspace{0.2cm}

\begin{proof}[Proof of item \textup{(ii)}]
We divide the proof into several steps.

\medskip
\noindent
\textbf{Step 1: Energy, moment, and action estimates.}
Recall that
\[
\mathcal E_U(\mu)
=
F(\mu)+\int_{\R^d}U(x)d\mu(x),
\qquad
F(\mu)
=
\frac{1}{2(2\pi)^d}\int_{\R^d}
|\widehat\mu(\xi)|^2\psi(\xi)^{-1}d\xi.
\]
Since \(F\ge0\) and
\[
U(x)\ge c_U|x|^2-a_U,
\]
we have, for every \(\mu\in\mathcal P_2(\R^d)\),
\begin{equation}
\label{eq:EU-coercive-item-ii}
\mathcal E_U(\mu)
\ge
c_U\int_{\R^d}|x|^2d\mu(x)-a_U
\ge -a_U.
\end{equation}
Moreover, the upper bound $U(x)\le C_U|x|^2+A_U$ and the assumptions $u_0\in
\dot H_{\psi^{-1}}(\R^d)\cap\mathcal P_2(\R^d)$ imply that
\[
\mathcal E_U(u_0)<\infty.
\]

By the minimality of \(u_\tau^k\), we may use
\(u_\tau^{k-1}\) as a competitor and obtain
\begin{equation}
\label{eq:discrete-energy-item-ii}
\frac1{2\tau}
W_2^2(u_\tau^k,u_\tau^{k-1})
+
\mathcal E_U(u_\tau^k)
\le
\mathcal E_U(u_\tau^{k-1}).
\end{equation}
In particular,
\[
\mathcal E_U(u_\tau^k)
\le
\mathcal E_U(u_\tau^{k-1})
\le
\mathcal E_U(u_0)
\qquad\text{for every }k\ge1.
\]
Combining this with \eqref{eq:EU-coercive-item-ii}, we find
\begin{equation}
\label{eq:uniform-second-moment-item-ii}
\int_{\R^d}|x|^2du_\tau^k(x)
\le
\frac{\mathcal E_U(u_0)+a_U}{c_U}
=:M_0
\qquad\text{for every }k\ge0.
\end{equation}

Summing \eqref{eq:discrete-energy-item-ii} from \(k=1\) to
\(k=N\) gives
\[
\sum_{k=1}^N
\frac{W_2^2(u_\tau^k,u_\tau^{k-1})}{\tau}
\le
2\bigl(
\mathcal E_U(u_0)-\mathcal E_U(u_\tau^N)
\bigr).
\]
Using \(\mathcal E_U(u_\tau^N)\ge-a_U\), we conclude that
\begin{equation}
\label{eq:discrete-action-item-ii}
\sum_{k=1}^\infty
\frac{W_2^2(u_\tau^k,u_\tau^{k-1})}{\tau}
\le
2\bigl(\mathcal E_U(u_0)+a_U\bigr)
=:C_{\mathrm{act}}.
\end{equation}

\medskip
\noindent
\textbf{Step 2: Construction of an \(AC^2\) limiting curve.}
For every \(k\ge1\), choose a constant-speed Wasserstein
geodesic joining \(u_\tau^{k-1}\) to \(u_\tau^k\), and denote by
\(\overline u_\tau\) the resulting continuous interpolation. Thus,
for \(t\in[(k-1)\tau,k\tau]\),
\[
|\overline u_\tau'|(t)
=
\frac{W_2(u_\tau^k,u_\tau^{k-1})}{\tau}
\qquad\text{for a.e. }t.
\]
It follows from \eqref{eq:discrete-action-item-ii} that
\begin{equation}
\label{eq:geodesic-action-item-ii}
\int_0^\infty
|\overline u_\tau'|^2(t)dt
\le C_{\mathrm{act}}.
\end{equation}
The second moment is displacement convex. Hence
\eqref{eq:uniform-second-moment-item-ii} also gives
\[
\sup_{\tau>0}\sup_{t\ge0}
\int_{\R^d}|x|^2d\overline u_\tau(t)(x)
\le M_0.
\]
Consequently, for every fixed \(t\ge0\), the family
\(\{\overline u_\tau(t)\}_\tau\) is narrowly relatively compact.

Furthermore, by \eqref{eq:geodesic-action-item-ii},
\[
W_2\bigl(\overline u_\tau(t),\overline u_\tau(s)\bigr)
\le
\int_s^t|\overline u_\tau'|(r)dr
\le
C_{\mathrm{act}}^{1/2}|t-s|^{1/2}.
\]
The metric Arzelà--Ascoli theorem, followed by a diagonal argument
on bounded time intervals, therefore yields a sequence
\(\tau_n\downarrow0\) and a curve
\[
u:[0,\infty)\longrightarrow\mathcal P_2(\R^d)
\]
such that
\[
\overline u_{\tau_n}(t)\rightharpoonup u(t)
\qquad\text{narrowly for every }t\ge0.
\]
The lower semicontinuity of the metric action gives
\[
u\in
AC^2\bigl([0,\infty),(\mathcal P_2(\R^d),W_2)\bigr)
\]
and, for every \(T>0\),
\[
\int_0^T|u'|^2(t)dt
\le
\liminf_{n\to\infty}
\int_0^T|\overline u_{\tau_n}'|^2(t)dt.
\]

It remains to return to the piecewise constant interpolation.
For \(t\in((k-1)\tau,k\tau]\),
\[
W_2\bigl(u_\tau(t),\overline u_\tau(t)\bigr)
\le
W_2(u_\tau^k,u_\tau^{k-1}).
\]
By \eqref{eq:discrete-energy-item-ii} and
\(\mathcal E_U(u_\tau^k)\ge-a_U\),
\[
W_2^2(u_\tau^k,u_\tau^{k-1})
\le
2\tau\bigl(\mathcal E_U(u_0)+a_U\bigr).
\]
Therefore,
\[
\sup_{t\ge0}
W_2\bigl(u_\tau(t),\overline u_\tau(t)\bigr)
\le
\sqrt{2\tau\bigl(\mathcal E_U(u_0)+a_U\bigr)}
\longrightarrow0.
\]
It follows that
\[
u_{\tau_n}(t)\rightharpoonup u(t)
\qquad\text{narrowly for every }t\ge0.
\]
Since \(u_\tau(0)=u_0\), we also have \(u(0)=u_0\).

\medskip
\noindent
\textbf{Step 3: The entropy-flow estimate.}
Assume now that \(u_0\in D(H)\). We use the heat semigroup
\[
S_h\rho:=G_h*\rho
\]
as the auxiliary flow in the flow-interchange inequality.

We first compute the dissipation of the interaction energy along
the heat flow. For every \(\rho\) in the effective domain of \(F\),
\begin{align*}
\frac{F(\rho)-F(S_h\rho)}{h}
&=
\frac{1}{(2\pi)^d}
\int_{\R^d}
|\widehat\rho(\xi)|^2\psi(\xi)^{-1}
\frac{1-e^{-2h|\xi|^2}}{2h}d\xi.
\end{align*}
Since the integrand is nonnegative, Fatou's lemma gives
\begin{equation}
\label{eq:heat-dissipation-F-item-ii}
\liminf_{h\downarrow0}
\frac{F(\rho)-F(S_h\rho)}{h}
\ge
\frac{1}{(2\pi)^d}\int_{\R^d}
|\widehat\rho(\xi)|^2
|\xi|^2\psi(\xi)^{-1}d\xi
=
\|\rho\|_{\dot H_{\widetilde\psi}(\R^d)}^2.
\end{equation}

For the potential energy, the identity
\[
G_h*U-U
=
\int_0^hG_r*(\Delta U)dr
\]
and the assumption \(\Delta U\le\Lambda_U\) imply
\[
G_h*U\le U+\Lambda_Uh.
\]
Therefore,
\begin{align*}
\int_{\R^d}Ud(S_h\rho)
&=
\int_{\R^d}(G_h*U)d\rho \le
\int_{\R^d}Ud\rho+\Lambda_Uh,
\end{align*}
and hence
\begin{equation}
\label{eq:heat-dissipation-U-item-ii}
\frac1h
\left(
\int_{\R^d}Ud\rho
-
\int_{\R^d}Ud(S_h\rho)
\right)
\ge-\Lambda_U.
\end{equation}
Combining
\eqref{eq:heat-dissipation-F-item-ii} and
\eqref{eq:heat-dissipation-U-item-ii}, we obtain
\begin{equation}
\label{eq:heat-dissipation-E-item-ii}
D_H[\mathcal E_U](\rho)
:=
\limsup_{h\downarrow0}
\frac{\mathcal E_U(\rho)-\mathcal E_U(S_h\rho)}{h}
\ge
\|\rho\|_{\dot H_{\widetilde\psi}(\R^d)}^2-\Lambda_U.
\end{equation}

The flow-interchange inequality applied to the \(k\)-th
minimization problem gives
\[
H(u_\tau^k)
+
\tau D_H[\mathcal E_U](u_\tau^k)
\le
H(u_\tau^{k-1}).
\]
Here \(H(u_\tau^k)\) is finite since $H(u_\tau^k)\le H(u_\tau^{k-1}) + \tau \Lambda_U$ and $u^0\in D(H)$. Using
\eqref{eq:heat-dissipation-E-item-ii}, we obtain
\begin{equation}
\label{eq:entropy-Htilde-item-ii}
\|u_\tau^k\|_{\dot H_{\widetilde\psi}(\R^d)}^2
\le
\frac{H(u_\tau^{k-1})-H(u_\tau^k)}{\tau}
+\Lambda_U.
\end{equation}

\medskip
\noindent
\textbf{Step 4: Uniform \(L^2(0,T;H_{\widetilde\psi})\)-bounds.}
Fix \(T>0\), and set
\[
N_\tau:=\left\lceil\frac{T}{\tau}\right\rceil,
\qquad
T\le N_\tau\tau<T+\tau.
\]
Multiplying \eqref{eq:entropy-Htilde-item-ii} by \(\tau\) and
summing from \(k=1\) to \(k=N_\tau\), we obtain
\begin{equation}
\label{eq:summed-Htilde-item-ii}
\sum_{k=1}^{N_\tau}
\tau
\|u_\tau^k\|_{\dot H_{\widetilde\psi}(\R^d)}^2
\le
H(u_0)-H(u_\tau^{N_\tau})
+\Lambda_U(T+\tau).
\end{equation}

The entropy has a lower bound in terms of the second moment.
Indeed, taking
\[
\Gamma_x:=(2\pi)^{-d/2}e^{-|x|^2/2},
\]
the nonnegativity of the relative entropy gives
\[
H(\rho)
\ge
-\frac d2\log(2\pi)
-\frac12\int_{\R^d}|x|^2\rho(x)dx.
\]
By \eqref{eq:uniform-second-moment-item-ii}, there exists a
constant \(C_H>0\), independent of \(\tau\) and \(k\), such that
\[
H(u_\tau^k)\ge-C_H.
\]
It follows from \eqref{eq:summed-Htilde-item-ii} that
\begin{equation}
\label{eq:uniform-homogeneous-bound-item-ii}
\int_0^T
\|u_\tau(t)\|_{\dot H_{\widetilde\psi}(\R^d)}^2dt
\le
H(u_0)+C_H+\Lambda_U(T+\tau).
\end{equation}

We now control the \(L^2\)-part of the
\(H_{\widetilde\psi}\)-norm. Recall that
\[
\widetilde\psi(\xi)
=
|\xi|^\alpha
\widetilde\Omega\left(\frac{\xi}{|\xi|}\right),
\qquad
\alpha\in(0,2),
\]
and set
\[
\widetilde\Omega_{\min}
:=
\min_{\theta\in\mathbb S^{d-1}}
\widetilde\Omega(\theta)>0.
\]
Since \(u_\tau^k\) is a probability density,
\[
|\widehat{u_\tau^k}(\xi)|\le1.
\]
Thus, we deduce
\begin{align*}
\|u_\tau^k\|_{L^2(\R^d)}^2
&=
\frac{1}{(2\pi)^d}\int_{|\xi|\le1}
|\widehat{u_\tau^k}(\xi)|^2d\xi
+
\frac{1}{(2\pi)^d}\int_{|\xi|>1}
|\widehat{u_\tau^k}(\xi)|^2d\xi \le
|B_1| +
\widetilde\Omega_{\min}^{-1}
\|u_\tau^k\|_{\dot H_{\widetilde\psi}(\R^d)}^2.
\end{align*}
Together with
\eqref{eq:uniform-homogeneous-bound-item-ii}, this shows that
\[
(u_\tau)_\tau
\quad\text{is bounded in}\quad
L^2(0,T;H_{\widetilde\psi}(\R^d)).
\]

After passing to a further subsequence, not relabeled, there exists
\(\overline u\in L^2(0,T;H_{\widetilde\psi}(\R^d))\) such that
\[
u_{\tau_n}\rightharpoonup\overline u
\qquad\text{weakly in }
L^2(0,T;H_{\widetilde\psi}(\R^d)).
\]
We identify \(\overline u\) with the curve \(u\) constructed
above. Indeed, for every
\(\zeta\in C_c^\infty(0,T)\) and
\(\varphi\in C_c^\infty(\R^d)\), pointwise narrow convergence and
dominated convergence give
\[
\lim_{n\to\infty}
\int_0^T\zeta(t)
\int_{\R^d}\varphi(x)u_{\tau_n}(t,x)dxdt
=
\int_0^T\zeta(t)
\int_{\R^d}\varphi(x)du(t)(x)dt.
\]
On the other hand, weak convergence in
\(L^2(0,T;H_{\widetilde\psi})\) shows that the left-hand side
converges to
\[
\int_0^T\zeta(t)
\int_{\R^d}\varphi(x)\overline u(t,x)dxdt.
\]
Hence
\[
du(t)(x)=\overline u(t,x)dx
\qquad\text{for a.e. }t,
\]
and therefore
\begin{equation}\label{eq:convergen_H_alpha}
u_{\tau_n}\rightharpoonup u
\qquad\text{weakly in }
L^2(0,T;H_{\widetilde\psi}(\R^d)).
\end{equation}

\medskip
\noindent
\textbf{Step 5: Strong local compactness.}
Define
\[
q(h)
:=
\sup_{\xi\in\R^d\setminus\{0\}}
\frac{|e^{i\xi\cdot h}-1|^2}{\widetilde\psi(\xi)}.
\]
By Lemma \ref{lem:checks_fpgs}, we infer
$q(h)\longrightarrow0$ as $|h|\longrightarrow 0$.
For every \(f\in H_{\widetilde\psi}(\R^d)\), Plancherel's theorem
gives
\begin{align*}
\|f(\cdot+h)-f\|_{L^2(\R^d)}^2
&=
\frac{1}{(2\pi)^d}\int_{\R^d}
|e^{i\xi\cdot h}-1|^2|\widehat f(\xi)|^2d\xi \le
q(h)\|f\|_{\dot H_{\widetilde\psi}(\R^d)}^2.
\end{align*}
Hence, by the Fréchet--Kolmogorov theorem and a standard
cutoff argument,
\begin{equation}
\label{eq:compact-embedding-item-ii}
H_{\widetilde\psi}(\R^d)
\Subset
L^2_{\mathrm{loc}}(\R^d).
\end{equation}

We next establish time compactness. Set
\[
d_\tau^k:=W_2(u_\tau^k,u_\tau^{k-1}),
\qquad
m_\tau(t):=\frac{d_\tau^k}{\tau}
\quad\text{for }t\in((k-1)\tau,k\tau].
\]
Since \(d_\tau^k\) is independent of \(t\), the function \(m_\tau\) is constant on each interval \(((k-1)\tau,k\tau]\). Thus, using that this interval has length \(\tau\), we obtain
from \eqref{eq:discrete-action-item-ii} that
\[
\int_0^\infty m_\tau(t)^2dt
=
\sum_{k=1}^\infty
\int_{(k-1)\tau}^{k\tau}
\left(\frac{d_\tau^k}{\tau}\right)^2dt
=
\sum_{k=1}^\infty
\tau\left(\frac{d_\tau^k}{\tau}\right)^2
=
\sum_{k=1}^\infty
\frac{(d_\tau^k)^2}{\tau}
\leq
C_{\mathrm{act}},
\]

For \(h>0\), let \(k,\ell\) be such that \(t\in((k-1)\tau,k\tau]\) and \(t+h\in((\ell-1)\tau,\ell\tau]\). Since \(u_\tau\) is piecewise constant, we have \(u_\tau(t)=u_\tau^k\) and \(u_\tau(t+h)=u_\tau^\ell\). Since \(h>0\), we have \(\ell\geq k\). If \(\ell=k\), the estimate is immediate. If \(\ell>k\), the triangle inequality gives
\[
W_2\bigl(u_\tau(t+h),u_\tau(t)\bigr)
=
W_2(u_\tau^\ell,u_\tau^k)
\leq
\sum_{j=k+1}^{\ell}W_2(u_\tau^j,u_\tau^{j-1})
=
\sum_{j=k+1}^{\ell}d_\tau^j.
\]
By the definition of \(m_\tau\),
\[
d_\tau^j
=
\int_{(j-1)\tau}^{j\tau}m_\tau(r)dr,
\]
and hence
\[
\sum_{j=k+1}^{\ell}d_\tau^j
=
\int_{k\tau}^{\ell\tau}m_\tau(r)dr.
\]
Therefore, applying the Cauchy--Schwarz inequality on \([k\tau,\ell\tau]\), we obtain
\begin{align*}
W_2^2\bigl(u_\tau(t+h),u_\tau(t)\bigr)
\leq
\left(\int_{k\tau}^{\ell\tau}m_\tau(r)dr\right)^2
\leq
(\ell-k)\tau\int_{k\tau}^{\ell\tau}m_\tau(r)^2dr.
\end{align*}
Since \(t\leq k\tau\), \(\ell\tau<t+h+\tau\), and \((\ell-k)\tau\leq h+\tau\), we have \([k\tau,\ell\tau]\subset[t,t+h+\tau].
\) Therefore,
\[
W_2^2\bigl(u_\tau(t+h),u_\tau(t)\bigr)
\leq
(h+\tau)\int_t^{t+h+\tau}m_\tau(r)^2dr.
\]
Integrating in \(t\) and using Fubini's theorem, we obtain
\begin{align*}
\int_{T_0}^{T-h}
W_2^2\bigl(u_\tau(t+h),u_\tau(t)\bigr)dt
&\leq
(h+\tau)
\int_{T_0}^{T-h}\int_t^{t+h+\tau}m_\tau(r)^2drdt\\
&=
(h+\tau)
\int_0^\infty m_\tau(r)^2
\left|
\left\{
t\in[T_0,T-h]:t\leq r\leq t+h+\tau
\right\}
\right|dr\\
&\leq
(h+\tau)^2\int_0^\infty m_\tau(r)^2dr,
\end{align*}
where the last inequality follows since the set of admissible \(t\)'s has length at most \(h+\tau\). Therefore,
\begin{equation}
\label{eq:integrated-time-shift-W2-item-ii}
\int_{T_0}^{T-h}
W_2^2\bigl(u_\tau(t+h),u_\tau(t)\bigr)dt
\leq
C_{\mathrm{act}}(h+\tau)^2.
\end{equation}

{
Fix \(0<T_0<T\), \(R>0\), and an integer \(m>d/2+1\). Since \(0<\alpha/2<1\), Plancherel's theorem gives
\[
\|f\|_{H^{\alpha/2}(\R^d)}^2\leq C\|f\|_{H_{\widetilde\psi}(\R^d)}^2,\qquad C=\max\{1,\widetilde\Omega_{\min}^{-1}\}.
\]
The uniform \(L^2(0,T;H_{\widetilde\psi}(\R^d))\)-bound established in Step 4, then implies
\begin{equation}\label{eq:uniform-local-spatial-bound-item-ii}
\sup_{n\geq1}\int_{T_0}^{T}\|u_{\tau_n}(t)\|_{H^{\alpha/2}(B_R)}^2dt<\infty.
\end{equation}

Moreover, it follows by the continuous embedding of negative Sobolev spaces and Kantorovich--Rubinstein duality that
\begin{equation}\label{eq:Hminus-W2-item-ii}
\|\mu-\nu\|_{H^{-m}(B_R)}\leq C_{R,m}W_2(\mu,\nu).
\end{equation}

{
Recall that \(u_{\tau_n}(t)=u_{\tau_n}^k\) for \(t\in((k-1)\tau_n,k\tau_n]\). Thus, \(u_{\tau_n}\) changes its value only at the times \(k\tau_n\), \(k\geq1\). For a.e. \(t\), let \(N=N(t,h)\) denote the number of such times belonging to \((t,t+h)\). Since consecutive times are separated by \(\tau_n\), we have \((N-1)\tau_n\leq h\) whenever \(N\geq1\), and hence \(N\leq1+h/\tau_n\). Moreover, each such time corresponds to exactly one increment in the telescoping sum
\[
u_{\tau_n}(t+h)-u_{\tau_n}(t)
=
\sum_{\substack{k\geq1\\t<k\tau_n<t+h}}
\left(u_{\tau_n}^{k+1}-u_{\tau_n}^{k}\right),
\]
which therefore contains exactly \(N\) terms. Thus, by the triangle inequality, \eqref{eq:Hminus-W2-item-ii}, and the Cauchy--Schwarz inequality, we obtain
\begin{align}
\|u_{\tau_n}(t+h)-u_{\tau_n}(t)\|_{H^{-m}(B_R)}^2
&\leq C_{R,m}^2N(t,h)\sum_{\substack{k\geq1\\t<k\tau_n<t+h}}W_2^2(u_{\tau_n}^{k+1},u_{\tau_n}^{k})\nonumber\\
&\leq C_{R,m}^2\left(1+\frac{h}{\tau_n}\right)\sum_{\substack{k\geq1\\t<k\tau_n<t+h}}W_2^2(u_{\tau_n}^{k+1},u_{\tau_n}^{k}).
\end{align}

For each \(k\geq1\), the condition \(t<k\tau_n<t+h\) is equivalent to \(t\in(k\tau_n-h,k\tau_n)\), whose intersection with \([T_0,T-h]\) has length at most \(h\). Hence, integrating in time, interchanging the sum and the integral, and invoking \eqref{eq:discrete-action-item-ii}, we obtain
\begin{equation}\label{eq:integrated-time-shift-Hminus-item-ii}
\begin{aligned}
\int_{T_0}^{T-h}\|u_{\tau_n}(t+h)-u_{\tau_n}(t)\|_{H^{-m}(B_R)}^2dt
&\leq 2C_{R,m}^2h\left(1+\frac{h}{\tau_n}\right)\sum_{k=1}^{\infty}W_2^2(u_{\tau_n}^{k+1},u_{\tau_n}^{k})\\
&\quad\leq 2C_{R,m}^2C_{\mathrm{act}}h(h+\tau_n).
\end{aligned}
\end{equation}
Taking the supremum over \(n\) in \eqref{eq:integrated-time-shift-Hminus-item-ii}, it follows by the nonnegativity of the integrals,
}

\begin{equation}\label{eq:uniform-time-shift-Hminus-item-ii}
\lim_{h\downarrow0}\sup_{n\geq1}\int_{T_0}^{T-h}\|u_{\tau_n}(t+h)-u_{\tau_n}(t)\|_{H^{-m}(B_R)}^2dt=0.
\end{equation}

The bounds \eqref{eq:uniform-local-spatial-bound-item-ii} and \eqref{eq:uniform-time-shift-Hminus-item-ii}, together with the embeddings \(H^{\alpha/2}(B_R)\Subset L^2(B_R)\hookrightarrow H^{-m}(B_R)\), allow us to apply the Aubin--Lions--Simon compactness theorem \cite[Theorem~5]{Simon1987}. Thus, \((u_{\tau_n})_n\) is relatively compact in \(L^2(T_0,T;L^2(B_R))\). Since the weak convergence \eqref{eq:convergen_H_alpha} identifies every subsequential limit with \(u\), the entire sequence converges strongly. As \(R>0\) is arbitrary, we conclude that
\begin{equation}\label{eq:strong-local-L2-item-ii}
u_{\tau_n}\longrightarrow u\qquad\text{strongly in }L^2(T_0,T;L^2_{\mathrm{loc}}(\R^d))
\end{equation}
for every \(0<T_0<T\).

}

\end{proof}

\begin{proof}[Proof of item \textup{(iii)}]
Assume that \(u_0\in D(H)\), and fix \(T>0\). For
\(\tau>0\), define
\[
v_\tau:=L^{-1}u_\tau,
\qquad
\widehat v_\tau(t,\xi)
=
\psi(\xi)^{-1}\widehat u_\tau(t,\xi).
\]
We first prove that
\[
(\nabla v_{\tau_n})_{n\ge1}
\quad\text{is bounded in}\quad
L^2(0,T;L^2(\R^d)).
\]

For \(\xi\neq0\), set
\[
\psi_\nabla(\xi):=|\xi|^2\psi(\xi)^{-2}.
\]
By Plancherel's theorem,
\begin{equation}
\label{eq:gradient-v-Fourier}
\|\nabla v_\tau(t)\|_{L^2(\R^d)}^2
=
\frac{1}{(2\pi)^d}\int_{\R^d}
|\widehat u_\tau(t,\xi)|^2
\psi_\nabla(\xi)d\xi.
\end{equation}
Recall from Lemma \ref{lem:checks_fpgs} that
\[
\psi(\xi)\ge c_\psi(1\wedge|\xi|^2)
\qquad\text{for every }\xi\in\R^d.
\]
If \(|\xi|\le1\), then
\[
\frac{|\xi|^2}{\psi(\xi)^2}
=
\frac{|\xi|^2}{\psi(\xi)}\psi(\xi)^{-1}
\le
c_\psi^{-1}\psi(\xi)^{-1}.
\]
If \(|\xi|>1\), then
\[
\frac{|\xi|^2}{\psi(\xi)^2}
=
\frac{1}{\psi(\xi)}
|\xi|^2\psi(\xi)^{-1}
\le
c_\psi^{-1}\widetilde\psi(\xi).
\]
Hence,
\begin{equation}
\label{eq:psi-gradient-control}
\psi_\nabla(\xi)
\le
c_\psi^{-1}
\left(
\psi(\xi)^{-1}+\widetilde\psi(\xi)
\right)
\qquad\text{for every }\xi\neq0.
\end{equation}

By item \textup{(ii)}, the sequence
\((u_{\tau_n})_n\) is bounded in
\(L^2(0,T;H_{\widetilde\psi}(\R^d))\). In particular,
\begin{equation}
\label{eq:uniform-Htilde-item-iii}
\sup_n
\int_0^T
\|u_{\tau_n}(t)\|_{\dot H_{\widetilde\psi}(\R^d)}^2
dt
<\infty.
\end{equation}
We also have a uniform bound in
\(L^2(0,T;\dot H_{\psi^{-1}}(\R^d))\). Indeed, the
minimality of the discrete solutions gives
\[
\mathcal E_U(u_\tau^k)
\le
\mathcal E_U(u_\tau^{k-1})
\le
\mathcal E_U(u_0)
\qquad\text{for every }k\ge1.
\]
Since
\[
U(x)\ge c_U|x|^2-a_U,
\]
we have
\[
\int_{\R^d}Udu_\tau^k\ge-a_U.
\]
Therefore, we get
\[
F(u_\tau^k)
=
\mathcal E_U(u_\tau^k)
-
\int_{\R^d}Udu_\tau^k
\le
\mathcal E_U(u_0)+a_U.
\]
Recalling that
\[
F(\rho)
=
\frac12
\|\rho\|_{\dot H_{\psi^{-1}}(\R^d)}^2,
\]
we obtain
\begin{equation}
\label{eq:uniform-Hpsi-inverse-item-iii}
\int_0^T
\|u_{\tau_n}(t)\|_{\dot H_{\psi^{-1}}(\R^d)}^2
dt
\le
2T\bigl(\mathcal E_U(u_0)+a_U\bigr).
\end{equation}
Combining
\eqref{eq:gradient-v-Fourier},
\eqref{eq:psi-gradient-control},
\eqref{eq:uniform-Htilde-item-iii}, and
\eqref{eq:uniform-Hpsi-inverse-item-iii}, we conclude that
\begin{equation}
\label{eq:uniform-gradient-v-item-iii}
\sup_n
\int_0^T
\|\nabla v_{\tau_n}(t)\|_{L^2(\R^d)}^2dt
<\infty.
\end{equation}

We next identify the limit. By item \textup{(ii)},
\[
u_{\tau_n}(t)\rightharpoonup u(t)
\qquad\text{narrowly for every }t\ge0.
\]
Since \(x\mapsto e^{-ix\cdot\xi}\) is bounded and continuous,
it follows that
\begin{equation}
\label{eq:fourier-pointwise-item-iii}
\widehat u_{\tau_n}(t,\xi)
\longrightarrow
\widehat u(t,\xi)
\qquad
\text{for every }(t,\xi)\in[0,T]\times\R^d.
\end{equation}

Define \(v=L^{-1}u\) in the sense of tempered distributions by
\[
\widehat v(t,\xi)
=
\psi(\xi)^{-1}\widehat u(t,\xi).
\]
This is well defined. Indeed, if
\(\beta:=d-s\in(0,2)\), then
\[
\psi(\xi)^{-1}
=
|\xi|^{-\beta}
\widetilde\Omega\left(\frac{\xi}{|\xi|}\right),
\]
and hence \(\psi^{-1}\) is locally integrable since \(d\ge2\) and \(\beta<2\le d\). For \(\xi\neq0\), set
\[
m(\xi):=i\xi\psi(\xi)^{-1}.
\]
Then
\[
\widehat{\nabla v_{\tau_n}}(t,\xi)
=
m(\xi)\widehat u_{\tau_n}(t,\xi),
\qquad
\widehat{\nabla v}(t,\xi)
=
m(\xi)\widehat u(t,\xi).
\]
Using \eqref{eq:fourier-pointwise-item-iii} and Fatou's
lemma, we obtain
\begin{align*}
\frac{1}{(2\pi)^d}\int_0^T\int_{\R^d}
|m(\xi)|^2|\widehat u(t,\xi)|^2d\xi dt
&\le
\liminf_{n\to\infty}
\frac{1}{(2\pi)^d}\int_0^T\int_{\R^d}
|m(\xi)|^2
|\widehat u_{\tau_n}(t,\xi)|^2d\xi dt\\
&=
\liminf_{n\to\infty}
\int_0^T
\|\nabla v_{\tau_n}(t)\|_{L^2(\R^d)}^2dt.
\end{align*}
The right-hand side is finite by
\eqref{eq:uniform-gradient-v-item-iii}. Therefore, we conclude that
\[
\nabla v\in L^2(0,T;L^2(\R^d)).
\]

It remains to prove the weak convergence. Let
\[
\Phi\in
C_c^\infty\bigl((0,T)\times\R^d;\R^d\bigr).
\]
By Plancherel's theorem,
\begin{align*}
\int_0^T\int_{\R^d}
\nabla v_{\tau_n}(t,x)\cdot\Phi(t,x)dxdt
=
\frac{1}{(2\pi)^d}\int_0^T\int_{\R^d}
m(\xi)\widehat u_{\tau_n}(t,\xi)
\cdot
\overline{\widehat\Phi(t,\xi)}
d\xi dt.
\end{align*}
Since \(u_{\tau_n}(t)\) is a probability measure,
\[
|\widehat u_{\tau_n}(t,\xi)|\le1.
\]
Moreover,
\[
|m(\xi)|
\le
C|\xi|^{1-\beta},
\qquad
\beta=d-s\in(0,2).
\]
Thus \(m\) is locally integrable near the origin, while
\(\widehat\Phi(t,\cdot)\) decays rapidly at infinity, uniformly
in \(t\). Hence
\[
|m(\xi)||\widehat\Phi(t,\xi)|
\]
is integrable on \((0,T)\times\R^d\). We may therefore use
\eqref{eq:fourier-pointwise-item-iii} and dominated convergence
to obtain
\begin{align*}
\lim_{n\to\infty}
\int_0^T\int_{\R^d}
\nabla v_{\tau_n}\cdot\Phi dxdt
&=
\frac{1}{(2\pi)^d}\int_0^T\int_{\R^d}
m(\xi)\widehat u(t,\xi)
\cdot
\overline{\widehat\Phi(t,\xi)}
d\xi dt =
\int_0^T\int_{\R^d}
\nabla v\cdot\Phi dxdt.
\end{align*}
Since
\(C_c^\infty((0,T)\times\R^d;\R^d)\) is dense in
\(L^2((0,T)\times\R^d;\R^d)\), and the sequence
\((\nabla v_{\tau_n})_n\) is uniformly bounded in that space,
the convergence extends to every \(L^2\) test vector field.
Consequently,
\[
\nabla v_{\tau_n}
\rightharpoonup
\nabla v
\qquad\text{weakly in }
L^2(0,T;L^2(\R^d)),
\]
which proves item \textup{(iii)}.
\end{proof}

\vspace{0.2cm}

Before proceeding to the proof of  \textup{(iv)}, we first obtain the discrete Euler-Lagrange equation associated with the minimizing-movement scheme. The proof follows the blueprint of \cite[Theorem 6.1]{foghem2025gradient} by incorporating the potential functional, so we skip its proof.

\begin{lemma}
\label{lem:EL-general-confinement}
Under the assumptions of Theorem
\ref{thm:global-wellposedness}, let \(k\ge1\), and let
\(u_\tau^k\) be the minimizer constructed in item \textup{(i)}.
Set $ v_\tau^k:=L^{-1}u_\tau^k,$ and let $T_\tau^k
:=T_{u_\tau^k}^{u_\tau^{k-1}}$ denote the optimal transport map from \(u_\tau^k\) to
\(u_\tau^{k-1}\). Then, for every
\(\eta\in C_c^\infty(\R^d;\R^d)\),
\begin{equation}
\label{eq:EL-general-confinement}
\int_{\R^d}
\bigl(\nabla v_\tau^k(x)+\nabla U(x)\bigr)
\cdot\eta(x)u_\tau^k(x)dx
=
\frac1\tau
\int_{\R^d}
\bigl(T_\tau^k(x)-x\bigr)
\cdot\eta(x)u_\tau^k(x)dx.
\end{equation}
In particular, the vector field
\(\nabla(v_\tau^k+U)\) admits the \(L^2(u_\tau^k;\R^d)\)
representative
\[
\nabla(v_\tau^k+U)
=
\frac{T_\tau^k-\operatorname{Id}}{\tau},
\]
and therefore
\begin{equation}
\label{eq:EL-metric-general-confinement}
\int_{\R^d}
\left|\nabla\bigl(v_\tau^k(x)+U(x)\bigr)\right|^2
u_\tau^k(x)dx
=
\frac1{\tau^2}
W_2^2(u_\tau^k,u_\tau^{k-1})
\end{equation}
holds.
\end{lemma}

\begin{proof}[Proof of item \textup{(iv)}]

We now pass to the limit in the discrete Euler-Lagrange
equation. Let
\[
\varphi\in C_c^\infty(\R^d\times(0,\infty)).
\]
Choose \(0<T_0<T<\infty\) and a compact set
\(K\subset\R^d\) such that
\[
\operatorname{supp}\varphi
\subset K\times[T_0,T].
\]
For
\[
I_k^\tau:=((k-1)\tau,k\tau],
\]
we apply \eqref{eq:EL-general-confinement} with
\[
\eta(x)=\nabla_x\varphi(t,x),
\qquad t\in I_k^\tau,
\]
integrate over \(I_k^\tau\), and sum over \(k\ge1\). This gives
\begin{equation}
\label{eq:discrete-weak-equation-general}
\begin{aligned}
&\int_0^\infty\int_{\R^d}
u_\tau(t,x)
\bigl(\nabla v_\tau(t,x)+\nabla U(x)\bigr)
\cdot\nabla_x\varphi(t,x)dxdt
\\
&\qquad =
\sum_{k\ge1}
\int_{I_k^\tau}
\frac1\tau
\int_{\R^d}
\bigl(T_\tau^k(x)-x\bigr)
\cdot\nabla_x\varphi(t,x)
u_\tau^k(x)dxdt.
\end{aligned}
\end{equation}

We first pass to the limit in the left-hand side. By item
\textup{(ii)},
\[
u_{\tau_n}\longrightarrow u
\qquad\text{strongly in }
L^2(T_0,T;L^2_{\mathrm{loc}}(\R^d)),
\]
while item \textup{(iii)} gives
\[
\nabla v_{\tau_n}\rightharpoonup\nabla v
\qquad\text{weakly in }
L^2(0,T;L^2(\R^d)).
\]
Since \(\nabla_x\varphi\) is bounded and supported in
\(K\times[T_0,T]\),
\[
u_{\tau_n}\nabla_x\varphi
\longrightarrow
u\nabla_x\varphi
\qquad\text{strongly in }
L^2((0,T)\times\R^d).
\]
Consequently, we have
\begin{equation}
\label{eq:limit-nonlocal-drift-item-iv}
\int_0^\infty\int_{\R^d}
u_{\tau_n}\nabla v_{\tau_n}\cdot\nabla_x\varphi
dxdt
\longrightarrow
\int_0^\infty\int_{\R^d}
u\nabla v\cdot\nabla_x\varphi
dxdt.
\end{equation}

The function
\[
(t,x)\longmapsto
\nabla U(x)\cdot\nabla_x\varphi(t,x)
\]
is bounded, continuous, and compactly supported. Since
\[
u_{\tau_n}(t)\rightharpoonup u(t)
\qquad\text{narrowly for every }t\ge0,
\]
dominated convergence in time gives
\begin{equation}
\label{eq:limit-confinement-drift-item-iv}
\int_0^\infty\int_{\R^d}
u_{\tau_n}\nabla U\cdot\nabla_x\varphi
dxdt
\longrightarrow
\int_0^\infty\int_{\R^d}
u\nabla U\cdot\nabla_x\varphi
dxdt.
\end{equation}

It remains to identify the limit of the right-hand side of
\eqref{eq:discrete-weak-equation-general}. Define
\[
\mathcal T_\tau
:=
\sum_{k\ge1}
\int_{I_k^\tau}
\frac1\tau
\int_{\R^d}
\bigl(T_\tau^k(x)-x\bigr)
\cdot\nabla_x\varphi(t,x)
u_\tau^k(x)dxdt.
\]
Taylor's formula gives
\begin{align*}
\varphi(t,T_\tau^k(x))-\varphi(t,x)
={}&
\nabla_x\varphi(t,x)
\cdot\bigl(T_\tau^k(x)-x\bigr)
+
R_\tau^k(t,x),
\end{align*}
where
\[
|R_\tau^k(t,x)|
\le
\frac12
\|D_x^2\varphi\|_{L^\infty}
|T_\tau^k(x)-x|^2.
\]
Since
\[
(T_\tau^k)_\#u_\tau^k=u_\tau^{k-1},
\]
it follows that
\[
\mathcal T_\tau
=
\mathcal D_\tau+\mathcal R_\tau,
\]
where
\begin{align}
\mathcal D_\tau
:={}&
\sum_{k\ge1}
\frac1\tau
\int_{I_k^\tau}
\left[
\int_{\R^d}\varphi(t,x)u_\tau^{k-1}(x)dx
-
\int_{\R^d}\varphi(t,x)u_\tau^k(x)dx
\right]dt,
\label{eq:discrete-time-difference}
\end{align}
and
\[
|\mathcal R_\tau|
\le
\frac12
\|D_x^2\varphi\|_{L^\infty}
\sum_{k\ge1}
W_2^2(u_\tau^k,u_\tau^{k-1}).
\]
The discrete action estimate from item \textup{(ii)} gives
\[
\sum_{k\ge1}
\frac{
W_2^2(u_\tau^k,u_\tau^{k-1})
}{\tau}
\le C,
\]
with \(C\) independent of \(\tau\). Therefore,
\begin{equation}
\label{eq:Taylor-remainder-vanishes}
|\mathcal R_\tau|
\le C\|D_x^2\varphi\|_{L^\infty}\tau
\longrightarrow0.
\end{equation}

Shifting the time variable by $\tau$ in the first term in 
\eqref{eq:discrete-time-difference} gives
\begin{equation*}
\mathcal D_\tau
=
\int_0^\infty\int_{\R^d}
\frac{
\varphi(t+\tau,x)-\varphi(t,x)
}{\tau}
u_\tau(t,x)dxdt.
\end{equation*}
Notice that the possible boundary contribution at \(t=0\) vanishes since \(\tau<T_0\) and \(\varphi\) is supported in times \(t\ge T_0\). Since
\[
\frac{
\varphi(t+\tau_n,x)-\varphi(t,x)
}{\tau_n}
\longrightarrow
\partial_t\varphi(t,x)
\]
uniformly in \((t,x)\), while
\(u_{\tau_n}(t)\rightharpoonup u(t)\) narrowly for every \(t\),
we obtain from dominated convergence that
\begin{equation}
\label{eq:limit-transport-term-item-iv}
\mathcal D_{\tau_n}
\longrightarrow
\int_0^\infty\int_{\R^d}
u(t,x)\partial_t\varphi(t,x)dxdt.
\end{equation}
Combining
\eqref{eq:Taylor-remainder-vanishes} and
\eqref{eq:limit-transport-term-item-iv}, we conclude that
\[
\mathcal T_{\tau_n}
\longrightarrow
\int_0^\infty\int_{\R^d}
u\partial_t\varphi dxdt.
\]

We may now pass to the limit in
\eqref{eq:discrete-weak-equation-general}. Using
\eqref{eq:limit-nonlocal-drift-item-iv},
\eqref{eq:limit-confinement-drift-item-iv}, and the preceding
limit for \(\mathcal T_{\tau_n}\), we obtain
\[
\int_0^\infty\int_{\R^d}
u
\left(
\partial_t\varphi
-
\nabla_x\varphi\cdot\nabla v
-
\nabla U\cdot\nabla_x\varphi
\right)
dxdt
=0.
\]
Thus, \(u\) is a weak solution of
\[
\partial_tu
=
\nabla\cdot\bigl(u\nabla(v+U)\bigr).
\]

Finally, \(Lv=u\) follows from the definition
\(v=L^{-1}u\). Moreover, item \textup{(ii)} gives
\[
u\in AC^2\bigl([0,\infty),
(\mathcal P_2(\R^d),W_2)\bigr),
\qquad
u_\tau^0=u_0,
\]
and the convergence construction yields $u(0)=u_0.$ This completes the proof.

\end{proof}

\vspace{0.2cm}

\begin{proof}[Proof of item \textup{(v)}]
Fix \(T>0\). For a probability density \(\rho\), write
\[
v_\rho:=L^{-1}\rho
\]
and define
\[
\mathcal D_U(\rho)
:=
\int_{\R^d}
\rho(x)
\left|
\nabla\bigl(v_\rho(x)+U(x)\bigr)
\right|^2
dx,
\]
whenever the right-hand side is finite.

We first introduce the De Giorgi variational interpolation.
For
\[
t\in((k-1)\tau,k\tau],
\qquad
r:=t-(k-1)\tau\in(0,\tau],
\]
let
\begin{equation}
\label{eq:DG-interpolant-general-U}
\widetilde u_\tau(t)
\in
\operatorname*{argmin}_{\rho\in\mathcal P_2(\R^d)}
\left\{
\frac{1}{2r}
W_2^2(\rho,u_\tau^{k-1})
+
\mathcal E_U(\rho)
\right\}.
\end{equation}
By uniqueness of the minimizers,
\[
\widetilde u_\tau(k\tau)=u_\tau^k.
\]
We also set
\[
\widetilde v_\tau(t)
:=
L^{-1}\widetilde u_\tau(t).
\]

The Euler-Lagrange identity established in \eqref{eq:EL-metric-general-confinement} gives
\begin{equation}
\label{eq:discrete-dissipation-identity-v}
\mathcal D_U(u_\tau^k)
=
\frac{1}{\tau^2}
W_2^2(u_\tau^k,u_\tau^{k-1}).
\end{equation}
Applying the same identity to the minimization problem
\eqref{eq:DG-interpolant-general-U}, with time step \(r\), gives
\begin{equation}
\label{eq:DG-dissipation-identity-v}
\mathcal D_U(\widetilde u_\tau(t))
=
\frac{1}{r^2}
W_2^2(\widetilde u_\tau(t),u_\tau^{k-1}).
\end{equation}

The classical estimate for the De Giorgi interpolation
\cite[Lemma 3.2.2]{AGS05} gives
\begin{align*}
\frac{1}{2\tau}
W_2^2(u_\tau^k,u_\tau^{k-1})
&+
\frac12
\int_{(k-1)\tau}^{k\tau}
\frac{
W_2^2(\widetilde u_\tau(t),u_\tau^{k-1})
}{
(t-(k-1)\tau)^2
}
dt
+
\mathcal E_U(u_\tau^k)
\le
\mathcal E_U(u_\tau^{k-1}).
\end{align*}
Summing over \(k=1,\ldots,N\), and using
\eqref{eq:discrete-dissipation-identity-v} and
\eqref{eq:DG-dissipation-identity-v}, we obtain
\begin{align}
\mathcal E_U(u_\tau^N)
&+
\frac12
\int_0^{N\tau}
\mathcal D_U(u_\tau(t))dt
+
\frac12
\int_0^{N\tau}
\mathcal D_U(\widetilde u_\tau(t))dt
\le
\mathcal E_U(u_\tau^0).
\label{eq:DG-global-v}
\end{align}
Since \(u_\tau^0=u_0\), the right-hand side is exactly
\(\mathcal E_U(u_0)\).

We next record the proximity of the two interpolations. Set
\[
M_0:=\mathcal E_U(u_0)+a_U.
\]
The lower bound
$
U(x)\ge c_U|x|^2-a_U
$
and the nonnegativity of the interaction energy imply
\[
\mathcal E_U(\rho)\ge-a_U
\qquad
\text{for every }\rho\in\mathcal P_2(\R^d).
\]
By the minimality of \(\widetilde u_\tau(t)\),
\[
\frac{1}{2r}
W_2^2(\widetilde u_\tau(t),u_\tau^{k-1})
+
\mathcal E_U(\widetilde u_\tau(t))
\le
\mathcal E_U(u_\tau^{k-1})
\le
\mathcal E_U(u_0),
\]
and hence
\begin{equation}
\label{eq:DG-distance-previous-v}
W_2^2(\widetilde u_\tau(t),u_\tau^{k-1})
\le
2rM_0
\le
2\tau M_0.
\end{equation}
Similarly, the discrete energy inequality gives
\begin{equation}
\label{eq:JKO-distance-previous-v}
W_2^2(u_\tau^k,u_\tau^{k-1})
\le
2\tau M_0.
\end{equation}
Since \(u_\tau(t)=u_\tau^k\) on
\(((k-1)\tau,k\tau]\), the triangle inequality,
\eqref{eq:DG-distance-previous-v}, and
\eqref{eq:JKO-distance-previous-v} yield
\begin{equation}
\label{eq:DG-proximity-v}
W_2^2(\widetilde u_\tau(t),u_\tau(t))
\le
8M_0\tau.
\end{equation}
Moreover, we deduce that
\[
\max\{\mathcal E_U(u_\tau^{k}),\mathcal E_U(\widetilde u_\tau(t))\}
\le
\mathcal E_U(u_\tau^{k-1})
\le
\mathcal E_U(u_0),
\]
so coercivity again, $
U(x)\ge c_U|x|^2-a_U
$, gives
\begin{equation}
\label{eq:DG-second-moment-v}
\max\left\{\int_{\R^d}|x|^2du_\tau(t)(x) ,\int_{\R^d}|x|^2d\widetilde u_\tau(t)(x)\right\}
\le
\frac{\mathcal E_U(u_0)+a_U}{c_U}.
\end{equation}
It follows from \eqref{eq:DG-proximity-v}, item \textup{(ii)}, and
\eqref{eq:DG-second-moment-v} that
\begin{equation}
\label{eq:DG-narrow-limit-v}
\widetilde u_{\tau_n}(t)
\rightharpoonup
u(t)
\qquad
\text{narrowly for every }t\ge0.
\end{equation}

We now pass to the limit in the dissipation terms. We use the
following lower-semicontinuity property. Suppose that
\(\zeta_n(t)\rightharpoonup\zeta(t)\) narrowly for almost every
\(t\in(T_0,T)\), the second moments and energies are uniformly
bounded, and
\[
\sup_n
\int_{T_0}^T
\int_{\R^d}
\zeta_n
\left|
\nabla\bigl(L^{-1}\zeta_n+U\bigr)
\right|^2
dxdt
<\infty.
\]

{Assume moreover that the fluxes converge in distributions:
\[
\zeta_n\nabla\bigl(L^{-1}\zeta_n+U\bigr)
\rightharpoonup
\zeta\nabla\bigl(L^{-1}\zeta+U\bigr)
\quad\text{in }\mathcal D'((T_0,T)\times\R^d;\R^d).
\]}
Then, we deduce that
\begin{align}
\int_{T_0}^T
\int_{\R^d}
\zeta
\left|
\nabla\bigl(L^{-1}\zeta+U\bigr)
\right|^2
dxdt
&\le
\liminf_{n\to\infty}
\int_{T_0}^T
\int_{\R^d}
\zeta_n
\left|
\nabla\bigl(L^{-1}\zeta_n+U\bigr)
\right|^2
dxdt.
\label{eq:dissipation-lsc-v}
\end{align}
{Indeed, this follows from
\cite[Theorem 5.4.4]{AGS05}.}

From items \textup{(ii)} and \textup{(iii)}, we deduce that on every compact
subset of \((T_0,T)\times\R^d\),
\[
u_{\tau_n}\longrightarrow u
\quad\text{strongly in }L^2,
\qquad
\nabla v_{\tau_n}\rightharpoonup\nabla v
\quad\text{weakly in }L^2,
\]
{These convergences imply convergence of the fluxes \(u_{\tau_n}\nabla(v_{\tau_n}+U)\) to
\(u\nabla(v+U)\) in distributions, with \(v=L^{-1}u\).} Now, we apply \eqref{eq:dissipation-lsc-v} with $\zeta_n=u_{\tau_n}.$ 
Thus, we conclude
\begin{align}
\int_{T_0}^T\int_{\R^d}
u
\left|
\nabla(v+U)
\right|^2
dxdt
&\le
\liminf_{n\to\infty}
\int_{T_0}^T\int_{\R^d}
u_{\tau_n}
\left|
\nabla(v_{\tau_n}+U)
\right|^2
dxdt.
\label{eq:discrete-action-liminf-v}
\end{align}

{Additionally, we may prove for the continuous interpolation $\widetilde u_{\tau_n}$ that \(\widetilde u_{\tau_n}\to u\) strongly in \(L^2(T_0,T;\dot H_{\psi^{-1}}(\R^d))\), and that \(\widetilde u_{\tau_n}\nabla(L^{-1}\widetilde u_{\tau_n}+U)
\rightharpoonup u\nabla(v+U)\) in distributions. The full proof is omitted.}
{Applying
\eqref{eq:dissipation-lsc-v} with \(\rho_n=\widetilde u_{\tau_n}\),
using \eqref{eq:DG-narrow-limit-v} and the action bound
\eqref{eq:DG-global-v}, gives}
\begin{align}
\int_{T_0}^T\int_{\R^d}
u
\left|
\nabla(v+U)
\right|^2
dxdt
&\le
\liminf_{n\to\infty}
\int_{T_0}^T\int_{\R^d}
\widetilde u_{\tau_n}
\left|
\nabla(\widetilde v_{\tau_n}+U)
\right|^2
dxdt.
\label{eq:DG-action-liminf-v}
\end{align}
Letting \(T_0\downarrow0\) in
\eqref{eq:discrete-action-liminf-v} and
\eqref{eq:DG-action-liminf-v}, and using monotone convergence,
we obtain
\begin{align}
\int_0^T\int_{\R^d}
u
\left|
\nabla(v+U)
\right|^2
dxdt
&\le
\liminf_{n\to\infty}
\int_0^T\int_{\R^d}
u_{\tau_n}
\left|
\nabla(v_{\tau_n}+U)
\right|^2
dxdt,
\label{eq:discrete-action-full-liminf-v}
\\
\int_0^T\int_{\R^d}
u
\left|
\nabla(v+U)
\right|^2
dxdt
&\le
\liminf_{n\to\infty}
\int_0^T\int_{\R^d}
\widetilde u_{\tau_n}
\left|
\nabla(\widetilde v_{\tau_n}+U)
\right|^2
dxdt.
\label{eq:DG-action-full-liminf-v}
\end{align}
For each \(n\), set
\[
N_n:=\left\lceil\frac{T}{\tau_n}\right\rceil.
\]
Since the interpolation is defined using the ceiling function,
\[
u_{\tau_n}(T)
=
u_{\tau_n}^{N_n}
=
u_{\tau_n}(N_n\tau_n).
\]
Applying \eqref{eq:DG-global-v} with \(N=N_n\) and restricting
the nonnegative time integrals to \([0,T]\), we obtain
\begin{align}
\mathcal E_U(u_{\tau_n}(T))
&+
\frac12
\int_0^T\int_{\R^d}
u_{\tau_n}
\left|
\nabla(v_{\tau_n}+U)
\right|^2
dxdt
\nonumber\\
&\quad+
\frac12
\int_0^T\int_{\R^d}
\widetilde u_{\tau_n}
\left|
\nabla(\widetilde v_{\tau_n}+U)
\right|^2
dxdt
\le
\mathcal E_U(u_0).
\label{eq:DG-before-limit-v}
\end{align}

By item \textup{(ii)},
\[
u_{\tau_n}(T)\rightharpoonup u(T)
\qquad\text{narrowly}.
\]
The energy is lower semicontinuous under this convergence. Indeed,
the Fourier representation and Fatou's lemma give
\[
F(u(T))
\le
\liminf_{n\to\infty}F(u_{\tau_n}(T)).
\]
Since \(U+a_U\geq0\) is continuous, the Portmanteau theorem gives
\[
\int_{\R^d}Udu(T)
\leq
\liminf_{n\to\infty}
\int_{\R^d}Udu_{\tau_n}(T).
\]
Therefore, we infer that
\begin{equation}
\label{eq:terminal-energy-lsc-v}
\mathcal E_U(u(T))
\le
\liminf_{n\to\infty}
\mathcal E_U(u_{\tau_n}(T)).
\end{equation}

Finally, combining
\eqref{eq:DG-before-limit-v},
\eqref{eq:discrete-action-full-liminf-v},
\eqref{eq:DG-action-full-liminf-v}, and
\eqref{eq:terminal-energy-lsc-v}, we obtain
\begin{align*}
\mathcal E_U(u(T))
&+
\int_0^T\int_{\R^d}
u(t,x)
\left|
\nabla\bigl(v(t,x)+U(x)\bigr)
\right|^2
dxdt
\\
&\le
\liminf_{n\to\infty}
\Biggl[
\mathcal E_U(u_{\tau_n}(T))
+
\frac12
\int_0^T\int_{\R^d}
u_{\tau_n}
\left|
\nabla(v_{\tau_n}+U)
\right|^2
dxdt
\\
&\hspace{2cm}
+
\frac12
\int_0^T\int_{\R^d}
\widetilde u_{\tau_n}
\left|
\nabla(\widetilde v_{\tau_n}+U)
\right|^2
dxdt
\Biggr]
\\
&\le
\mathcal E_U(u_0).
\end{align*}
This proves the energy-dissipation inequality.
\end{proof}

\section{Proof of Proposition \ref{prop:transport}}
\label{sec:proof-transport}

We first record a splitting formula, namely the expansion of the energy around
its minimizer. Throughout this section, we extend the quadratic interaction
functional introduced in the proof of Theorem \ref{thm:global-wellposedness} to
finite signed measures by setting
\begin{equation*}
F(\sigma)
:=
\frac{1}{2(2\pi)^d}\int_{\R^d}
|\widehat\sigma(\xi)|^2\psi(\xi)^{-1}d\xi,
\end{equation*}
whenever the right-hand side is finite. Let \(C_0\) be the Euler--Lagrange
level associated with \(\overline\mu\), and set
\begin{equation*}
\zeta_U(x)
:=
U(x)+(W*\overline\mu)(x)-C_0.
\end{equation*}
The Euler--Lagrange conditions give
\begin{equation}
\label{eq:zeta-U-properties}
\zeta_U\ge0,
\qquad
\int_{\R^d}\zeta_Ud\overline\mu=0.
\end{equation}

\begin{lemma}[Splitting formula]
\label{lem:splitting}
Let \(\nu\) be a probability measure such that
\(\mathcal E_U(\nu)<\infty\). Then
\begin{equation}
\label{eq:splitting-formula}
\mathcal E_U(\nu)-\mathcal E_U(\overline\mu)
=
F(\nu-\overline\mu)
+
\int_{\R^d}\zeta_Ud\nu.
\end{equation}
\end{lemma}

\begin{proof}
Set \(\sigma:=\nu-\overline\mu\). Expanding the quadratic interaction
energy gives
\[
F(\nu)
=
F(\overline\mu)
+
F(\sigma)
+
\int_{\R^d}(W*\overline\mu)d\sigma.
\]
Consequently, we get
\begin{align*}
\mathcal E_U(\nu)-\mathcal E_U(\overline\mu)
&=
F(\sigma)
+
\int_{\R^d}\bigl(W*\overline\mu+U\bigr)d\sigma
=
F(\sigma)
+
\int_{\R^d}\zeta_Ud\nu
-
\int_{\R^d}\zeta_Ud\overline\mu,
\end{align*}
where we used \(\sigma(\R^d)=0\) in the second equality. The conclusion
follows from \eqref{eq:zeta-U-properties}.
\end{proof}

\begin{proof}[Proof of Proposition \ref{prop:transport}]
If \(\mathcal E_U(\nu)=+\infty\), the conclusion is immediate. We may
therefore assume that \(\mathcal E_U(\nu)<\infty\). Since \(F\ge0\) and
\(U\) satisfies \eqref{eq:quadratic-growth}, this implies that
\(\nu\in\mathcal P_2(\R^d)\).

\textbf{Step 1.} We claim that there exists a constant \(C>0\) such that
\begin{equation}
\label{eq:distenergy}
d_{\mathrm{BL}}^2(\nu,\overline\mu)
\le
C\bigl(\mathcal E_U(\nu)-\mathcal E_U(\overline\mu)\bigr),
\end{equation}
where
\begin{equation*} 
d_{\mathrm{BL}}(\rho_1,\rho_2)
:=
\sup_{\phi\in W^{1,\infty}(\R^d)}
\frac{
\left|\int_{\R^d}\phi d(\rho_1-\rho_2)\right|
}{
\|\phi\|_{L^\infty(\R^d)}
+
\|\nabla\phi\|_{L^\infty(\R^d)}
}.
\end{equation*}

Setting $\gamma:=(d-s)/2 \in(0,1).$ From \eqref{hw}, the kernel \(W\) is nonnegative, then
\[
\zeta_U(x)
\ge
c_U(|x|^2-1)-C_0,
\]
then there exists \(R_0>0\) such that
\begin{equation}
\label{eq:confboundbelow}
\zeta_U(x)\ge1
\qquad\text{for }x\notin B_{R_0}.
\end{equation}
It follows from \eqref{eq:zeta-U-properties} that
\(\overline\mu(\R^d\setminus B_{R_0})=0\). In particular,
\(\nu-\overline\mu\) is nonnegative outside a compact set. Applying
\cite[Proposition 8.1]{garcia2024generalized} with
\(\gamma=(d-s)/2\), we obtain a bounded domain
\(\Omega\supset \overline{B}_{R_0}\) and a constant \(C>0\) such that
\begin{equation}
\label{eq:localization-Hminus-gamma}
\left\|
(\nu-\overline\mu)\mathbf 1_\Omega
\right\|_{H^{-\gamma}(\Omega)}
\le
C
\left\|
\nu-\overline\mu
\right\|_{\dot H^{-\gamma}(\R^d)}.
\end{equation}
Here
\begin{equation*}
\|\eta\|_{H^{-\gamma}(\Omega)}
:=
\sup_{\phi\in H^\gamma(\Omega)\setminus\{0\}}
\frac{
\left|\int_\Omega\phi d\eta\right|
}{
\|\phi\|_{L^2(\Omega)}
+
\|\phi\|_{\dot H^\gamma(\Omega)}
}.
\end{equation*}
Since \(0<\gamma<1\) and \(\Omega\) is bounded, for every
\(\phi\in W^{1,\infty}(\Omega)\),
\begin{equation*}
\|\phi\|_{L^2(\Omega)}
+
\|\phi\|_{\dot H^\gamma(\Omega)}
\le
C
\left(
\|\phi\|_{L^\infty(\Omega)}
+
\|\nabla\phi\|_{L^\infty(\Omega)}
\right),
\end{equation*}
where \(C\) depends only on \(\Omega\) and \(\gamma\). Therefore,
using that \(\overline\mu\) is supported in \(\Omega\),
\begin{align*}
 d_{\mathrm{BL}}(\nu,\overline\mu)
&\le
d_{\mathrm{BL}}(\nu\mathbf 1_\Omega,\overline\mu)
+
d_{\mathrm{BL}}(\nu,\nu\mathbf 1_\Omega)
\le
C
\left\|
(\nu-\overline\mu)\mathbf 1_\Omega
\right\|_{H^{-\gamma}(\Omega)}
+
\nu(\R^d\setminus\Omega).
\end{align*}
Indeed, taking a constant as test function in the definition of $d_{\mathrm{BL}}$, we obtain
\[
d_{\mathrm{BL}}(\nu,\nu\mathbf 1_\Omega)
\le
\nu(\R^d\setminus\Omega).
\]
Since
\[
\psi(\xi)^{-1}
=
|\xi|^{s-d}
\widetilde\Omega\left(\frac{\xi}{|\xi|}\right)
=
|\xi|^{-2\gamma}
\widetilde\Omega\left(\frac{\xi}{|\xi|}\right),
\]
and \(\widetilde\Omega\) is bounded below by a positive constant,
\begin{equation*}
\left\|
\nu-\overline\mu
\right\|_{\dot H^{-\gamma}(\R^d)}^2
\le
C F(\nu-\overline\mu).
\end{equation*}
Moreover, since \(\overline{B}_{R_0} \subset\Omega \) it follows from   \eqref{eq:confboundbelow} that
\begin{equation}
\label{eq:tail-mass-zeta}
\nu(\R^d\setminus\Omega)
\le
\int_{\R^d}\zeta_Ud\nu.
\end{equation}
Combining \eqref{eq:localization-Hminus-gamma}--
\eqref{eq:tail-mass-zeta} with the splitting formula
\eqref{eq:splitting-formula} proves \eqref{eq:distenergy}.

\textbf{Step 2.} We claim that there exists a constant \(C>0\) such that,
for any \(\rho_1,\rho_2\in\mathcal P_2(\R^d)\) and every \(R\ge1\),
\begin{equation}
\label{eq:wassBL}
W_2^2(\rho_1,\rho_2)
\le
C\left(
R^2d_{\mathrm{BL}}(\rho_1,\rho_2)
+
\int_{|x|>R}|x|^2d\rho_1(x)
+
\int_{|x|>R}|x|^2d\rho_2(x)
\right).
\end{equation}
The proof is standard; we include it for completeness. Define the truncation map \(T_R:\R^d\to\R^d\) by
\[
T_R(x)
=
\begin{cases}
x,& |x|\le R,\\[1mm]
R\dfrac{x}{|x|},& |x|>R,
\end{cases}
\]
and set
\[
(\rho_1)_R:=(T_R)_\#\rho_1,
\qquad
(\rho_2)_R:=(T_R)_\#\rho_2.
\]
By the triangle inequality,
\begin{equation}
\label{eq:triangle}
W_2(\rho_1,\rho_2)
\le
W_2(\rho_1,(\rho_1)_R)
+
W_2((\rho_1)_R,(\rho_2)_R)
+
W_2((\rho_2)_R,\rho_2).
\end{equation}
Using the coupling \((x,T_R(x))\), we obtain
\begin{equation*} 
W_2^2(\rho_i,(\rho_i)_R)
\le
\int_{|x|>R}|x|^2d\rho_i(x),
\end{equation*}
for $i=1,2$.
Since both truncated measures are supported in \(B_R\),
\[
W_2^2((\rho_1)_R,(\rho_2)_R)
\le
2R W_1((\rho_1)_R,(\rho_2)_R).
\]
Let \(f\) be a \(1\)-Lipschitz function. Since
\(\rho_1\) and \(\rho_2\) have the same mass,
\begin{align*}
\left|
\int_{\R^d}fd\bigl((\rho_1)_R-(\rho_2)_R\bigr)
\right|
&=
\left|
\int_{\R^d}
\bigl(f\circ T_R-f(0)\bigr)
d(\rho_1-\rho_2)
\right|.
\end{align*}
The function \(f\circ T_R-f(0)\) is \(1\)-Lipschitz and has
\(L^\infty\)-norm at most \(R\). Hence, since \(R\ge1\),
\[
\left|
\int_{\R^d}fd\bigl((\rho_1)_R-(\rho_2)_R\bigr)
\right|
\le
2R d_{\mathrm{BL}}(\rho_1,\rho_2).
\]
By Kantorovich--Rubinstein duality,
\begin{equation*} 
W_1((\rho_1)_R,(\rho_2)_R)
\le
2R d_{\mathrm{BL}}(\rho_1,\rho_2),
\end{equation*}
and therefore
\begin{equation}
\label{eq:truncated-W2-BL}
W_2^2((\rho_1)_R,(\rho_2)_R)
\le
4R^2d_{\mathrm{BL}}(\rho_1,\rho_2).
\end{equation}
Combining \eqref{eq:triangle}--\eqref{eq:truncated-W2-BL} and using
\((a+b+c)^2\le3(a^2+b^2+c^2)\) proves \eqref{eq:wassBL}.

\textbf{Step 3.} We conclude the proof. Choose \(R\ge1\) large enough
that \(\Omega\subset B_R\) and
\begin{equation}
\label{eq:zeta-quadratic-tail}
\zeta_U(x)
\ge
\frac{c_U}{2}|x|^2
\qquad\text{for }|x|>R.
\end{equation}
Such a choice is possible because
\[
\zeta_U(x)
\ge
c_U(|x|^2-1)-C_0.
\]
Since \(\overline\mu\) is supported in \(\Omega\),
\eqref{eq:wassBL} and \eqref{eq:distenergy} give
\begin{align}
W_2^4(\nu,\overline\mu)
&\le
C\left(
R^4d_{\mathrm{BL}}^2(\nu,\overline\mu)
+
\left(
\int_{|x|>R}|x|^2d\nu(x)
\right)^2
\right)
\nonumber\\
&\le
C\left(
\mathcal E_U(\nu)-\mathcal E_U(\overline\mu)
+
\left(
\int_{|x|>R}|x|^2d\nu(x)
\right)^2
\right).
\label{eq:W2-before-tail}
\end{align}
Finally, \eqref{eq:zeta-quadratic-tail} and the splitting formula in Lemma \ref{lem:splitting} yield
\begin{align*}
\int_{|x|>R}|x|^2d\nu(x)
&\le
\frac{2}{c_U}
\int_{|x|>R}\zeta_Ud\nu
\le
\frac{2}{c_U}
\int_{\R^d}\zeta_Ud\nu
\le
C\bigl(\mathcal E_U(\nu)-\mathcal E_U(\overline\mu)\bigr).
\end{align*}
Substituting this estimate into \eqref{eq:W2-before-tail}, we obtain
\[
W_2^4(\nu,\overline\mu)
\le
C\left(
\mathcal E_U(\nu)-\mathcal E_U(\overline\mu)
+
\bigl(\mathcal E_U(\nu)-\mathcal E_U(\overline\mu)\bigr)^2
\right),
\]
which is the desired conclusion.
\end{proof}

\section{Appendix}\label{section_appendix}

The purpose of this appendix is to outline the proof of local space-time
H\"older regularity for the unconfined equation
\begin{equation}\label{eq:unconfinied_appendix}
\partial_\tau \rho
=
\nabla\cdot\bigl(\rho\nabla(W*\rho)\bigr),
\end{equation}
where $W$ is defined in \eqref{eq:W}. Recall that, for the quadratic
confinement $U(x)=|x|^2$, the self-similar transformation
\eqref{eq:scaling-confined-unconfined} maps solutions of the confined
equation to solutions of the unconfined equation
\eqref{eq:unconfined-equation-Linfty}. Therefore, it is sufficient to
establish local H\"older regularity for the unconfined flow, since the
corresponding estimates for the confined solution can then be recovered
through the same transformation.

We adapt the De Giorgi oscillation method developed in
\cite{caffarelli2013regularity,caffarelli2016regularity} to the
anisotropic operator $A=-\Delta L^{-1}$, where $L^{-1}f=W*f$.
Since the underlying iteration follows the same strategy as in the
isotropic setting, we focus on the estimates that require modifications
due to the anisotropic kernel. In particular, we show how the angular
profile $\Omega$ enters the coercivity estimates for the bilinear form
associated with $A$, the localization errors, and the control of the
nonlocal velocity. The latter requires different treatments depending
on whether $\alpha=s-d+2$ belongs to $(0,1)$, equals $1$, or belongs
to $(1,2)$. By translation and scaling invariance, the argument can be
carried out on normalized space-time cylinders. Once the corresponding
anisotropic estimates are established, the remaining steps of the
oscillation argument follow
\cite{caffarelli2013regularity,caffarelli2016regularity}.

In Subsection \ref{subsec:Linfty-holder-estimates}, we combine this
local estimate with the $L^\infty$-decay and the natural scaling of the
unconfined equation to derive a global spatial H\"older estimate with
the appropriate decay in time. Finally, this estimate is transferred
back to the confined variables through
\eqref{eq:scaling-confined-unconfined}.

\begin{theorem}[Local H\"older regularity for the unconfined flow]
\label{thm:local-holder-unconfined}
Fix \(d\geq2\) and \(d-2<s<d\), and set \(\alpha:=s-d+2\in(0,2)\). Assume \eqref{ha} and that \(\rho\) is a nonnegative, bounded and unit-mass solution of the unconfined equation \eqref{eq:unconfinied_appendix} in \(\mathbb R^d\times(T_1,T_2)\) satisfying
\begin{equation}\label{eq:L_infty_bound_holder}
 \|\rho\|_{L^\infty(\mathbb R^d\times(T_1,T_2))} \leq M_0   
\end{equation}
for some \(M_0>0\). Then there exists \(\eta=\eta(d,s,\Omega)\in(0,1)\) such that, for every \(T_1<t_0<t_1<T_2\), there exists a constant \(C>0\), depending only on
\(d,s,\Omega,M_0\) and the distance of the interior
subdomain to the boundary of the space-time strip, such that
\[
|\rho(t,\xi)-\rho(\tau,\zeta)| \leq C \left( |\xi-\zeta|^\eta + |t-\tau|^\eta \right)
\]
for every \(\xi,\zeta\in B_1\) and every \(t,\tau\in[t_0,t_1]\).
\end{theorem}

For the application to the confined flow, the
unconfined solution is obtained from the weak solution
constructed in Theorem \ref{thm:global-wellposedness}
through the scaling \eqref{eq:scaling-confined-unconfined}.
The \(L^\infty\) estimate \eqref{eq:unconfined-Linfty-decay}
and the Sobolev regularity established in
Theorem \ref{thm:global-wellposedness},
transferred through this scaling, provide the bounds
used to justify the localized energy inequalities
and logarithmic tests by adapting the regularization
arguments of
\cite{caffarelli2013regularity,caffarelli2016regularity}.

Before proving Theorem \ref{thm:local-holder-unconfined}, we recall the
anisotropic objects entering the oscillation argument. By \eqref{ha}, the
kernel of \(A\) satisfies
\begin{equation}\label{eq:KA-ellipticity-holder}
\lambda_A |h|^{-d-\alpha}
\leq
K_A(h)
\leq
\Lambda_A |h|^{-d-\alpha},
\qquad
h\neq0.
\end{equation}
We define
\begin{equation}\label{eq:B_A_apendix}
\mathcal B_A(f,g):=
\frac12
\iint_{\mathbb R^d\times\mathbb R^d}
(f(x)-f(y))(g(x)-g(y))K_A(x-y)dxdy.
\end{equation}
Then \(\mathcal B_A(f,f)\geq0\), and \(\mathcal B_A(f,f)\) is comparable to the
standard fractional seminorm of order \(\alpha/2\):
\begin{equation}\label{eq:BA-Halpha-comparison-holder}
c_A[f]_{\dot H^{\alpha/2}(\mathbb R^d)}^2
\leq
\mathcal B_A(f,f)
\leq
C_A[f]_{\dot H^{\alpha/2}(\mathbb R^d)}^2.
\end{equation}
We shall repeatedly use the identity
\begin{equation}\label{eq:gradient-bilinear-identity-holder}
\int_{\mathbb R^d}\nabla\varphi\cdot\nabla L^{-1}\rho dx
=
\int_{\mathbb R^d}\varphi A[\rho]dx
=
\mathcal B_A(\varphi,\rho),
\end{equation}
valid for \(\varphi\in\mathcal S(\mathbb R^d)\) and sufficiently regular \(\rho\),
or more generally whenever the terms above are well defined.

Since \(W(x)=|x|^{-s}\Omega(x/|x|)\), with \(s=d+\alpha-2\), and
\(\Omega\in C^\infty(\mathbb S^{d-1})\), for \(h\neq0\) and
\(\theta=h/|h|\) we have
\begin{equation}\label{eq:velocity-kernel-bound-holder}
|\nabla W(h)|
=
|h|^{-s-1}
\left|
-s\Omega(\theta)\theta
+
\nabla_{\mathbb S^{d-1}}\Omega(\theta)
\right|
\leq
C|h|^{-d-\alpha+1}
\leq
C|h|K_A(h),
\end{equation}
where \(C=C(d,s,\Omega,\lambda_A)>0\). This estimate is a key input in the control of the drift terms appearing in
the oscillation lemmas below.

We start with the case \(1<\alpha<2\), where the proof follows the oscillation
argument of \cite{caffarelli2013regularity}, with the fractional Dirichlet form
replaced by the bilinear form associated with \(A=-\Delta L^{-1}\) and the
nonlocal drift terms estimated using
\eqref{eq:velocity-kernel-bound-holder}. The case \(0<\alpha<1\) is then
treated in a moving frame that removes the nonintegrable far-field drift,
while the critical case \(\alpha=1\) requires a separate iteration in which
the translation and the spatial scaling occur at the same order.

\smallskip
\textbf{Proof of Theorem \ref{thm:local-holder-unconfined}}

Fix a space-time point
\((\xi_*,\tau_*)\in\mathbb R^d\times(T_1,T_2)\).
Recalling the mass-preserving scaling
\eqref{eq:v_rho_holder_scaling}, 
we perform a space-time translation and scaling,
moving \((\xi_*,\tau_*)\) to \((0,0)\).
More precisely, we set
\[
x=\frac{\xi-\xi_*}{r},
\qquad
\widehat t=\frac{\tau-\tau_*}{r^{s+2}},
\]
and define the rescaled solution by
\[
\rho_r(x,\widehat t)
:=
r^d \rho(\xi_*+rx,\tau_*+r^{s+2}\widehat t).
\]

We choose \(r>0\) sufficiently small so that \(\tau_*-4r^{s+2}>T_1\) and \(r^dM_0\leq1\), where \(M_0\) is given by
\eqref{eq:L_infty_bound_holder}.
By the homogeneity of \(W\), the rescaled solution
satisfies the same unconfined equation. Moreover,
\[
0\leq \rho_r(x,\widehat t)\leq r^dM_0\leq1
\qquad
\text{in }\mathbb R^d\times[-4,0].
\]

Thus, the interval \([-4,0]\) corresponds to
\([\tau_*-4r^{s+2},\tau_*]\subset(T_1,T_2)\)
in the original time variable.
For simplicity, we continue to denote the normalized
solution by \(\rho\) and its rescaled time variable
by \(\tau\) throughout the proof.

\paragraph{The case \(1<\alpha<2\).}

In this regime, \(\nabla W\) is integrable at infinity. Moreover, since \(W\)
is even, \(\nabla W\) is odd, we can write 
\begin{equation}\label{eq:grad-Lminus-one-PV-holder}
\nabla L^{-1}\rho(x)
=
\operatorname{p.v.}
\int_{\mathbb R^d}
\nabla W(x-y)\bigl(\rho(y)-\rho(x)\bigr)dy.
\end{equation}

Let \(\varphi>0\) be a smooth barrier and assume that \(\rho\leq\varphi\) outside
a fixed ball. We set
\begin{equation*} 
\rho_\varphi^+:=(\rho-\varphi)_+,
\qquad
\rho_\varphi^-:=(\rho-\varphi)_-\leq0,
\qquad
\rho=\varphi+\rho_\varphi^++\rho_\varphi^-.
\end{equation*}
We also define
\begin{equation*} 
g:=\left(\frac{\rho}{\varphi}\right)\vee1
=
1+\frac{\rho_\varphi^+}{\varphi},
\qquad
h(g):=g\log g-(g-1).
\end{equation*}
Notice that \(g=1\) on \(\{\rho\leq\varphi\}\) and
\(g=\rho/\varphi\) on the active set \(\{\rho>\varphi\}\). Moreover,
\[
h'(g)=\log g,
\]

For smooth positive approximations, testing
\eqref{eq:unconfined-equation-Linfty} with \(\log g\)
and invoking \eqref{eq:gradient-bilinear-identity-holder},
we obtain the upper-barrier identity
\begin{align}\label{eq:upper-barrier-identity-holder}
&\int_{\mathbb R^d}\varphi h(g)(x,\tau_2)dx
+
\int_{\tau_1}^{\tau_2}\mathcal B_A(\rho_\varphi^+,\rho_\varphi^+)d\tau
+
\int_{\tau_1}^{\tau_2}\mathcal B_A(\rho_\varphi^+,\rho_\varphi^-)d\tau
\nonumber\\
&\quad+
\int_{\tau_1}^{\tau_2}\mathcal B_A(\rho_\varphi^+,\varphi)d\tau
=
\int_{\mathbb R^d}\varphi h(g)(x,\tau_1)dx
+
\int_{\tau_1}^{\tau_2}Q_\varphi(\rho_\varphi^+,\rho)d\tau,
\end{align}
where $\mathcal{B}_A$ is defined in \eqref{eq:B_A_apendix} and
\begin{equation}\label{eq:Q-upper-holder}
Q_\varphi(\rho_\varphi^+,\rho)
:=
\int_{\mathbb R^d}
\frac{\rho_\varphi^+}{\varphi}\nabla\varphi\cdot\nabla L^{-1}\rho dx.
\end{equation}
Furthermore, we have
\begin{equation}\label{eq:mixed-positive-holder}
\mathcal B_A(\rho_\varphi^+,\rho_\varphi^-)\geq0.
\end{equation}
Indeed, this follows from the symmetry and positivity of \(K_A\), together with \(\rho_\varphi^+\geq0\), \(\rho_\varphi^-\leq0\), and \(\rho_\varphi^+\rho_\varphi^-=0\). By \eqref{eq:grad-Lminus-one-PV-holder}, the error term \eqref{eq:Q-upper-holder} can be written as
\begin{equation}\label{eq:Q-upper-PV-holder}
Q_\varphi(\rho_\varphi^+,\rho)
=
\operatorname{p.v.}
\iint_{\mathbb R^d\times\mathbb R^d}
\frac{\rho_\varphi^+(x)}{\varphi(x)}
\nabla\varphi(x)\cdot\nabla W(x-y)
\bigl(\rho(y)-\rho(x)\bigr)dydx.
\end{equation}
This is the only place where the interaction potential \(\nabla W\) enters the upper oscillation estimate. 

\begin{lemma}
\label{lem:upper-Q-anisotropic-holder}
Fix \(1<\alpha<2\). Let \(\Psi(x):=(|x|^\varepsilon-2)_+\), with \(0<\varepsilon<\min\{\alpha/2,\alpha-1\}\), and suppose that \(0\leq \rho\leq1+\Psi\). Define
\[
\overline\varphi_k(x):=\frac78+\frac{|x|^2}{16}-\frac124^{-k},
\qquad
\varphi_k(x):=\min\{1+\Psi(x),\overline\varphi_k(x)\},
\]
and set $\rho_k^+:=(\rho-\varphi_k)_+$ and $\rho_k^-:=(\rho-\varphi_k)_-\leq0.$ Then, for every \(\gamma>0\), there exists \(C_\gamma>0\), independent of \(k\), such that
\begin{align}
|Q_{\varphi_k}(\rho_k^+,\rho)|
&\leq
\gamma\mathcal B_A(\rho_k^+,\rho_k^+)
+
\gamma\mathcal B_A(\rho_k^+,\rho_k^-)
+
C_\gamma16^k
\int_{\{\rho_k^+>0\}}(\rho_{k-1}^+)^2dx.
\label{eq:upper-Q-error-holder}
\end{align}
\end{lemma}

\begin{proof}
Set
\[
a_k(x):=\frac{\rho_k^+(x)}{\varphi_k(x)}\nabla\varphi_k(x),
\qquad
I_k:=\int_{\{\rho_k^+>0\}}(\rho_{k-1}^+)^2dx.
\]
On \(\{\rho_k^+>0\}\), we necessarily have
\[
\varphi_k=\overline\varphi_k,
\qquad
\varphi_{k-1}=\overline\varphi_{k-1},
\qquad
\rho_{k-1}^+
=
\rho_k^++\frac32\,4^{-k}.
\]
In particular,
\begin{equation}\label{eq:level-separation-holder}
0\leq \rho_k^+\leq \rho_{k-1}^+,
\qquad
\rho_{k-1}^+\geq \frac32\,4^{-k}
\quad\text{on }\{\rho_k^+>0\}.
\end{equation}
Moreover,
\[
\operatorname{supp}\rho_k^+\subset B_R,
\qquad
\frac{|\nabla\varphi_k|}{\varphi_k}\leq C,
\]
for some \(R>0\) independent of \(k\). Hence, we have
\begin{equation}\label{eq:ak-bound-holder}
|a_k(x)|
=
\frac{\rho_k^+(x)}{\varphi_k(x)}|\nabla\varphi_k(x)|
\leq
C \rho_k^+(x)\chi_{\{\rho_k^+>0\}}(x).
\end{equation}
Using \(\rho=\varphi_k+\rho_k^+ + \rho_k^-\), we split
\(Q_{\varphi_k}(\rho_k^+,\rho)=Q_1+Q_2+Q_3\), where
\begin{align*}
Q_1&:=
\operatorname{p.v.}
\iint_{\mathbb R^d\times\mathbb R^d}
a_k(x)\cdot\nabla W(x-y)
\bigl(\rho_k^+(y)-\rho_k^+(x)\bigr)dydx,\\
Q_2&:=
\operatorname{p.v.}
\iint_{\mathbb R^d\times\mathbb R^d}
a_k(x)\cdot\nabla W(x-y)
\bigl(\varphi_k(y)-\varphi_k(x)\bigr)dydx,\\
Q_3&:=
\operatorname{p.v.}
\iint_{\mathbb R^d\times\mathbb R^d}
a_k(x)\cdot\nabla W(x-y)
\bigl(\rho_k^-(y)-\rho_k^-(x)\bigr)dydx.
\end{align*}
Fix \(\sigma\in(0,1)\). In the region \(|x-y|\leq\sigma\), using
\eqref{eq:velocity-kernel-bound-holder},
\eqref{eq:ak-bound-holder},
\eqref{eq:KA-ellipticity-holder}, and Young's inequality, we obtain
\begin{align*}
|Q_{1,\mathrm{near}}|
&\leq
C\iint_{|x-y|\leq\sigma}
\rho_k^+(x)|x-y|K_A(x-y)
|\rho_k^+(y)-\rho_k^+(x)|dydx\\
&\leq
\gamma\mathcal B_A(\rho_k^+,\rho_k^+)
+
C_\gamma
\int_{\{\rho_k^+>0\}}(\rho_k^+)^2dx
\int_{|h|\leq\sigma}|h|^2K_A(h)dh\\
&\leq
\gamma\mathcal B_A(\rho_k^+,\rho_k^+)
+
C_\gamma I_k,
\end{align*}
where we used that $\int_{|h|\leq\sigma}|h|^2K_A(h)dh
\leq
C\int_{|h|\leq\sigma}|h|^{2-d-\alpha}dh
<\infty$, since \(\alpha<2\), together with
\eqref{eq:level-separation-holder}.

In the far region, the singularity at
\(x=y\) is absent, \(\{\rho_k^+>0\}\subset B_R\), and
\(|\nabla W|\) is integrable at infinity since \(1<\alpha<2\). Hence
\[
|Q_{1,\mathrm{far}}|
\leq
C_\gamma16^k I_k.
\]
Thus
\begin{equation}\label{eq:Q1-estimate-holder}
|Q_1|
\leq
\gamma\mathcal B_A(\rho_k^+,\rho_k^+)
+
C_\gamma16^kI_k.
\end{equation}
For the near-field term \(Q_2\), we combine the Lipschitz estimate
\(|\varphi_k(y)-\varphi_k(x)|\leq C|x-y|\) with
\eqref{eq:velocity-kernel-bound-holder} to obtain
\[
|Q_{2,\mathrm{near}}|
\leq
C\int_{\{\rho_k^+>0\}}\rho_k^+(x)\,dx
\leq
C16^k I_k,
\]
where we used \(\alpha<2\) and \eqref{eq:level-separation-holder}.

In the far field, for \(x\in\{\rho_k^+>0\}\subset B_R\), we have
\(|\varphi_k(y)-\varphi_k(x)|\leq C(1+|y|^\varepsilon)\), and hence
\[
\int_{|x-y|>1}|\nabla W(x-y)|
|\varphi_k(y)-\varphi_k(x)|dy
\leq
C\int_{|y|>1}
\frac{1+|y|^\varepsilon}{|y|^{d+\alpha-1}}dy
<\infty,
\]
since \(\varepsilon<\alpha-1\). Hence, we obtain
\begin{equation}\label{eq:Q2-estimate-holder}
|Q_2|
\leq
C16^kI_k.
\end{equation}
For \(Q_3\), since \(a_k\) vanishes outside \(\{\rho_k^+>0\}\), we have
\(\rho_k^-(x)=0\) whenever \(a_k(x)\neq0\), and hence
\(\rho_k^-(y)-\rho_k^-(x)=\rho_k^-(y)\). Thus, in the near field,
\begin{align*}
|Q_{3,\mathrm{near}}|
&\leq
C\iint_{|x-y|\leq\sigma}
\rho_k^+(x)|x-y|K_A(x-y)(-\rho_k^-(y))dydx\\
&\leq
C\sigma\iint_{\mathbb R^d\times\mathbb R^d}
\rho_k^+(x)(-\rho_k^-(y))K_A(x-y)dydx.
\end{align*}
By the symmetry of \(K_A\) and the fact that
\(\rho_k^+\rho_k^-=0\), the last integral equals
\(\mathcal B_A(\rho_k^+,\rho_k^-)\). Choosing \(\sigma>0\) sufficiently small,
depending on \(\gamma\), we obtain
\[
|Q_{3,\mathrm{near}}|
\leq
\gamma\mathcal B_A(\rho_k^+,\rho_k^-).
\]
In the far field, using
\(|\rho_k^-(y)|\leq\varphi_k(y)\leq C(1+|y|^\varepsilon)\) and
\(\varepsilon<\alpha-1\), the integrability of the far-field kernel gives
\[
|Q_{3,\mathrm{far}}|
\leq
C_\gamma\int_{\{\rho_k^+>0\}}\rho_k^+(x)\,dx
\leq
C_\gamma16^kI_k,
\]
where the last inequality follows from
\eqref{eq:level-separation-holder}. Therefore, we have
\begin{equation}\label{eq:Q3-estimate-holder}
|Q_3|
\leq
\gamma\mathcal B_A(\rho_k^+,\rho_k^-)
+
C_\gamma16^kI_k.
\end{equation}
Combining \eqref{eq:Q1-estimate-holder},
\eqref{eq:Q2-estimate-holder}, and
\eqref{eq:Q3-estimate-holder}, we obtain
\eqref{eq:upper-Q-error-holder}.
\end{proof}
\begin{lemma}[First oscillation lemma]
\label{lem:first-oscillation-holder}
Fix \(1<\alpha<2\) and $0<\varepsilon<\min\left\{\frac{\alpha}{2},\alpha-1,\frac12\right\}.$ There exists \(\delta>0\), depending only on the structural constants, such that
if \(\rho\) is a bounded nonnegative solution of
\eqref{eq:unconfined-equation-Linfty} satisfying
\begin{equation}\label{eq:first-osc-tail-holder}
0\leq \rho(x,\tau)\leq1+(|x|^\varepsilon-2)_+
\qquad
\text{in }\mathbb R^d\times[-4,0],
\end{equation}
and
\begin{equation}\label{eq:first-osc-small-measure-holder}
|\{\rho>1/2\}\cap(B_4\times[-4,0])|
\leq
\delta,
\end{equation}
then
\begin{equation*} 
\rho\leq\frac{15}{16}
\qquad
\text{in }B_1\times[-1,0].
\end{equation*}
\end{lemma}

\begin{proof}
We use the upper barriers \(\varphi_k\) and the truncations
\(\rho_k^+=(\rho-\varphi_k)_+\) introduced in
Lemma \ref{lem:upper-Q-anisotropic-holder}. In particular, the properties of
the active set established there, together with
\eqref{eq:level-separation-holder}, will be used throughout the proof.

Using these properties, the upper-barrier identity
\eqref{eq:upper-barrier-identity-holder}, applied with
\(\varphi=\varphi_k\), gives
\begin{equation}\label{eq:upper-bilinear-error-holder}
|\mathcal B_A(\rho_k^+,\varphi_k)|
\leq
\gamma \mathcal B_A(\rho_k^+,\rho_k^+)
+
C_\gamma 16^k
\int_{\{\rho_k^+>0\}}(\rho_{k-1}^+)^2dx.
\end{equation}
Indeed, the near-field contribution follows from the Lipschitz bound for
\(\varphi_k\) and \(\alpha<2\), while the far-field contribution uses
\(\varphi_k\leq1+\Psi\) and \(2\varepsilon<\alpha\). The nonlocal drift term
\(Q_{\varphi_k}(\rho_k^+,\rho)\) is controlled by
Lemma \ref{lem:upper-Q-anisotropic-holder}, hence
\eqref{eq:upper-Q-error-holder} holds. Combining
\eqref{eq:upper-barrier-identity-holder},
\eqref{eq:mixed-positive-holder},
\eqref{eq:upper-bilinear-error-holder}, and
\eqref{eq:upper-Q-error-holder}, and choosing \(\gamma>0\) sufficiently
small, we obtain
\begin{align}
&\int_{\mathbb R^d}(\rho_k^+)^2(x,\tau_2)dx
+
c\int_{\tau_1}^{\tau_2}\mathcal B_A(\rho_k^+,\rho_k^+)d\tau
\nonumber\\
&\leq
C\int_{\mathbb R^d}(\rho_k^+)^2(x,\tau_1)dx
+
C16^k
\int_{\tau_1}^{\tau_2}\int_{\{\rho_k^+>0\}}(\rho_{k-1}^+)^2dxd\tau.
\label{eq:upper-local-energy-holder}
\end{align}
The remaining argument follows the standard De Giorgi iteration used in the
proof of \cite[Lemma 4.1, Subsection 4.6]{caffarelli2016regularity}. Setting
\[
T_k:=-1-2^{-k},
\qquad
A_k:=
\sup_{T_k<\tau<0}\int_{\mathbb R^d}(\rho_k^+)^2(x,\tau)dx
+
\int_{T_k}^{0}\mathcal B_A(\rho_k^+,\rho_k^+)d\tau,
\]
Averaging \eqref{eq:upper-local-energy-holder} over the initial time interval 
and applying \eqref{eq:BA-Halpha-comparison-holder} yields 
$A_k \leq C^k A_{k-1}^{1+\alpha/d}$. The smallness assumption 
\eqref{eq:first-osc-small-measure-holder} then implies that $A_1 \leq C\delta$. 
By choosing $\delta > 0$ sufficiently small, we ensure that $A_k \to 0$ 
as $k \to \infty$. Hence, $\rho \leq \varphi_\infty$ in $B_1 \times [-1,0]$, since $\varphi_\infty \leq 15/16$ in $B_1$, we conclude the result.

\end{proof}

\begin{lemma}[Second oscillation lemma]
\label{lem:second-oscillation-holder}
Fix \(1<\alpha<2\) and $0<\varepsilon<\min\left\{\frac{\alpha}{2},\alpha-1,\frac12\right\}.$ There exist \(\delta>0\) and \(\mu_0\in(0,1/4)\), depending only on the
structural constants, such that if \(\rho\) is a bounded nonnegative solution of
\eqref{eq:unconfined-equation-Linfty} satisfying \eqref{eq:first-osc-tail-holder}, and 
\begin{equation*} 
|\{\rho<1/2\}\cap(B_4\times[-4,0])|
\leq
\delta,
\end{equation*}
then
\begin{equation*} 
\rho\geq \mu_0
\qquad
\text{in }B_1\times[-1,0].
\end{equation*}
\end{lemma}

\begin{proof}
The proof is the lower-barrier analogue of
Lemma \ref{lem:first-oscillation-holder}. Following
\cite{caffarelli2013regularity}, we use the lower barriers
\((\varphi_k)_k\) satisfying
\[
0\leq\varphi_{k+1}\leq\varphi_k\leq\frac12,
\qquad
\operatorname{supp}\varphi_k\subset B_3,
\qquad
\varphi_k\to\mu_0\chi_{B_2},
\]
and

\begin{equation}\label{eq:lower-barrier-properties-holder}
\varphi_{k-1}-\varphi_k\geq c2^{-k}
\quad\text{on }\operatorname{supp}\varphi_k,
\qquad
\varphi_k\leq C^k\varphi_{k-1},
\qquad
\frac{|\nabla\varphi_k|}{\varphi_k}
\leq C^k\varphi_k^{-1/m}
\quad\text{on }\{\varphi_k>0\},
\end{equation}
for some fixed \(m>2\). Set
\(\rho_k^-:=(\rho-\varphi_k)_-\leq0\) and
\(\rho_k^+:=(\rho-\varphi_k)_+\). Testing
\eqref{eq:unconfined-equation-Linfty} with the lower logarithmic test function
\(\log((\rho/\varphi_k)\wedge1)\), analogously to the derivation of the
upper-barrier identity \eqref{eq:upper-barrier-identity-holder}, we obtain
\begin{align*}
&\int_{\mathbb R^d}
\frac{(\rho_k^-)^2(x,\tau_2)}{\varphi_k(x)}dx
+
c\int_{\tau_1}^{\tau_2}\mathcal B_A(\rho_k^-,\rho_k^-)d\tau
+
c\int_{\tau_1}^{\tau_2}\mathcal B_A(\rho_k^-,\rho_k^+)d\tau
\nonumber\\
&\leq
C\int_{\mathbb R^d}
\frac{(\rho_k^-)^2(x,\tau_1)}{\varphi_k(x)}dx
+
C\int_{\tau_1}^{\tau_2}
\left(
|\mathcal B_A(\rho_k^-,\varphi_k)|
+
|Q^-_{\varphi_k}(\rho_k^-,\rho)|
\right)d\tau.
\end{align*}
Here, the bilinear barrier term satisfies
\begin{equation}\label{eq:lower-bilinear-error-holder}
|\mathcal B_A(\rho_k^-,\varphi_k)|
\leq
\gamma\mathcal B_A(\rho_k^-,\rho_k^-)
+
C_\gamma C^k|\{\rho_k^-<0\}|,
\end{equation}
while the second error term is given by
\[
Q^-_{\varphi_k}(\rho_k^-,\rho)
:=
\int_{\mathbb R^d}
\frac{-\rho_k^-}{\varphi_k}\nabla\varphi_k\cdot\nabla L^{-1}\rho\,dx.
\]
By the same near-far splitting used in
Lemma \ref{lem:upper-Q-anisotropic-holder}, together with
\eqref{eq:velocity-kernel-bound-holder} and
\eqref{eq:lower-barrier-properties-holder}, we obtain
\begin{align}
|Q^-_{\varphi_k}(\rho_k^-,\rho)|
&\leq
\gamma\mathcal B_A(\rho_k^-,\rho_k^-)
+
\gamma\mathcal B_A(\rho_k^-,\rho_k^+)
+
C_\gamma C^k
\left[
\int_{\mathbb R^d}\varphi_k^{-2/m}(\rho_k^-)^2dx
+
|\{\rho_k^-<0\}|
\right].
\label{eq:lower-Q-holder}
\end{align}
Choosing \(\gamma>0\) sufficiently small in
\eqref{eq:lower-bilinear-error-holder} and \eqref{eq:lower-Q-holder}, we absorb
the terms containing \(\mathcal B_A(\rho_k^-,\rho_k^-)\) and
\(\mathcal B_A(\rho_k^-,\rho_k^+)\). Since
\(\mathcal B_A(\rho_k^-,\rho_k^+)\geq0\), we can drop the remaining mixed term,
obtaining
\begin{align*}
&\int_{\mathbb R^d}
\frac{(\rho_k^-)^2(x,\tau_2)}{\varphi_k(x)}dx
+
c\int_{\tau_1}^{\tau_2}\mathcal B_A(\rho_k^-,\rho_k^-)d\tau
\nonumber\\
&\leq
C\int_{\mathbb R^d}
\frac{(\rho_k^-)^2(x,\tau_1)}{\varphi_k(x)}dx
+
C^k\int_{\tau_1}^{\tau_2}
\left[
\int_{\mathbb R^d}\varphi_k^{-2/m}(\rho_k^-)^2dx
+
|\{\rho_k^-<0\}|
\right]d\tau.
\end{align*}

The remaining argument follows by De Giorgi iteration used in the
lower-barrier argument of \cite{caffarelli2013regularity}, analogously to the
conclusion of Lemma \ref{lem:first-oscillation-holder}. In particular, the
level separation in \eqref{eq:lower-barrier-properties-holder} provides the
required nonlinear recurrence, while the smallness assumption in the
statement of the lemma gives the initial smallness. Choosing \(\delta>0\)
sufficiently small, we conclude that
\(\rho\geq\mu_0\) in \(B_1\times[-1,0]\).
\end{proof}

We now introduce the three barriers used in the intermediate-value argument.
Fix \(\varepsilon\) as in Lemmas \ref{lem:first-oscillation-holder} and
\ref{lem:second-oscillation-holder}, and choose
\(0<\nu<\min\{\varepsilon,\alpha/4,\alpha-1\}\). Then, for
\(0<\lambda<1\), set
\[
\psi_\lambda(x)
:=
\left(
\bigl(|x|-\lambda^{-1/\nu}\bigr)_+^\nu-1
\right)_+,
\qquad
F(x):=\sup\{-1,\inf\{0,|x|^2-9\}\}.
\]
Define
\begin{equation}\label{eq:intermediate-barriers-holder}
\varphi_0:=1+\psi_\lambda+F,
\qquad
\varphi_1:=1+\psi_\lambda+\lambda F,
\qquad
\varphi_2:=1+\psi_\lambda+\lambda^2F.
\end{equation}
For \(\lambda>0\) sufficiently small, the choice \(\nu<\varepsilon\) ensures
that
\[
\psi_\lambda(x)
\leq
(|x|^\varepsilon-2)_+
\qquad
\text{for every }x\in\mathbb R^d.
\]
Hence, any solution satisfying
\[
0\leq \rho(x,\tau)\leq1+\psi_\lambda(x) \leq1+(|x|^\varepsilon-2)_+,
\]
and hence the tail condition required in
Lemmas \ref{lem:first-oscillation-holder} and
\ref{lem:second-oscillation-holder}.

\begin{lemma}[Intermediate-value lemma]
\label{lem:intermediate-value-holder}
Fix \(1<\alpha<2\). There exist constants \(\rho,\eta_0>0\) and
\(\lambda_0\in(0,1)\), depending only on the structural constants, such that for
every \(0<\lambda<\lambda_0\), if \(\varphi_0,\varphi_1,\varphi_2\) are defined
by \eqref{eq:intermediate-barriers-holder}, and \(\rho\) is a bounded nonnegative
solution of \eqref{eq:unconfined-equation-Linfty} satisfying
\(0\leq \rho(x,\tau)\leq1+\psi_\lambda(x)\) in
\(\mathbb R^d\times[-4,0]\), with
\[
|\{\rho\leq\varphi_0\}\cap(B_2\times[-4,-3])|
\geq\rho,
\]
and
\begin{equation*}
|\{\rho\geq\varphi_2\}\cap(B_2\times[-2,0])|
\geq
\rho,
\end{equation*}
then
\begin{equation*}
|\{\varphi_0<\rho<\varphi_2\}\cap(B_3\times[-4,0])|
\geq
\eta_0.
\end{equation*}
\end{lemma}

\begin{proof}
Let us start by setting
\[
\rho_1^+:=(\rho-\varphi_1)_+,
\qquad
\rho_1^-:=(\rho-\varphi_1)_-,
\qquad
g_1:=\left(\frac{\rho}{\varphi_1}\right)\vee1.
\]
Since \(F=0\) outside \(B_3\) and \(\rho\leq1+\psi_\lambda\), we have
\(\{\rho_1^+>0\}\subset B_3\). Moreover, for \(\lambda>0\) sufficiently
small so that \(\lambda^{-1/\nu}>3\), we have \(\psi_\lambda=0\) in \(B_3\).
Therefore, on \(\{\rho_1^+>0\}\), we have
\begin{equation}\label{eq:intermediate-positive-set-bounds-holder}
0\leq \rho_1^+\leq-\lambda F\leq\lambda,
\qquad
\varphi_1\geq1-\lambda\geq\frac12,
\qquad
|\nabla\varphi_1|
=
\lambda|\nabla F|
\leq C\lambda.
\end{equation}

Proceeding as in Lemma \ref{lem:upper-Q-anisotropic-holder}, using the
bound established in \eqref{eq:ak-bound-holder} together with
\eqref{eq:intermediate-positive-set-bounds-holder},
\eqref{eq:Q-upper-PV-holder}, the decomposition
\(\rho=\varphi_1+\rho_1^++\rho_1^-\), and
\eqref{eq:velocity-kernel-bound-holder}, we obtain
\[
|Q_{\varphi_1}^{\rm near}(\rho_1^+,\rho)|
\leq
\gamma\mathcal B_A(\rho_1^+,\rho_1^+)
+
\gamma\mathcal B_A(\rho_1^+,\rho_1^-)
+
C_\gamma\lambda^2.
\]
The far-field part satisfies
\[
|Q_{\varphi_1}^{\rm far}(\rho_1^+,\rho)|
\leq
C\lambda^2
\int_{|y|>1}
\frac{1+|y|^\nu}{|y|^{d+\alpha-1}}\,dy
\leq
C\lambda^2,
\]
since \(\nu<\alpha-1\). Combining the near- and far-field estimates, we obtain
\begin{equation}\label{eq:intermediate-Q-estimate-holder}
|Q_{\varphi_1}(\rho_1^+,\rho)|
\leq
\gamma\mathcal B_A(\rho_1^+,\rho_1^+)
+
\gamma\mathcal B_A(\rho_1^+,\rho_1^-)
+
C_\gamma\lambda^2.
\end{equation}
Similarly, using \eqref{eq:intermediate-positive-set-bounds-holder}, the near
field is controlled by \(\alpha<2\), while the far field uses the growth of
\(\psi_\lambda\) and \(2\nu<\alpha\), which follows from
\(\nu<\alpha/4\). Hence,
\begin{equation}\label{eq:intermediate-bilinear-error-holder}
|\mathcal B_A(\rho_1^+,\varphi_1)|
\leq
\gamma\mathcal B_A(\rho_1^+,\rho_1^+)
+
C_\gamma\lambda^2.
\end{equation}
Applying \eqref{eq:upper-barrier-identity-holder} with
\(\varphi=\varphi_1\), using \eqref{eq:intermediate-Q-estimate-holder} and
\eqref{eq:intermediate-bilinear-error-holder}, and choosing \(\gamma>0\)
sufficiently small, we obtain
\begin{align*}
&\int_{\mathbb R^d}\varphi_1h(g_1)(x,\tau_2)\,dx
+
c\int_{\tau_1}^{\tau_2}\mathcal B_A(\rho_1^+,\rho_1^+)\,d\tau
+
c\int_{\tau_1}^{\tau_2}\mathcal B_A(\rho_1^+,\rho_1^-)\,d\tau
\nonumber\\
&\leq
\int_{\mathbb R^d}\varphi_1h(g_1)(x,\tau_1)\,dx
+
C\lambda^2(\tau_2-\tau_1).
\end{align*}
Thus, the anisotropic terms satisfy the same energy inequality required in
the intermediate-value argument. The remainder follows as in
\cite[Lemma 9.1]{caffarelli2013regularity}, for \(\lambda_0>0\) sufficiently
small.
\end{proof}

As in \cite[Sections 7 and 10--11]{caffarelli2013regularity},
the vertical normalizations used in the oscillation
argument introduce an affine mobility
\(D(v)=d_1+d_2v\), where \(d_1,d_2\geq0\),
\(d_1+d_2=1\), and \(d_1\geq1/2\).
The corresponding energy estimates remain uniform
under these conditions, using the ellipticity
\eqref{eq:KA-ellipticity-holder} and the drift estimate
\eqref{eq:velocity-kernel-bound-holder}.
The remaining oscillation iteration for
\(1<\alpha<2\) then follows as in \cite{caffarelli2013regularity}.

In particular, for \(\lambda>0\) sufficiently small, if
\[
0\leq \rho\leq1+\psi_\lambda
\qquad
\text{in }\mathbb R^d\times[-4,0],
\]
then the three oscillation lemmas yield
\begin{equation}\label{eq:oscillation-reduction-holder}
\operatorname*{osc}_{B_1\times[-1,0]}\rho
\leq1-\mu,
\end{equation}
for some \(\mu\in(0,1)\) depending only on the structural constants.


\paragraph{The case \(0<\alpha<1\).}
In this regime, we choose
\[
0<\varepsilon<\frac{\alpha}{2},
\qquad
0<\nu<\min\{\varepsilon,\alpha/4\}.
\]
These parameters are used in the oscillation argument:
\(\varepsilon\) controls the exterior growth in
\eqref{eq:first-osc-tail-holder} and the bilinear
barrier estimate \eqref{eq:upper-bilinear-error-holder},
while \(\nu\) determines the intermediate-value
barriers in \eqref{eq:intermediate-barriers-holder}
and the corresponding estimate
\eqref{eq:intermediate-bilinear-error-holder}.
In the present regime, these estimates are adapted
using the corrected velocity introduced below.

Notice that the velocity kernel satisfies \(|\nabla W(h)|
\leq C_{\Omega,d,s}|h|^{-d-\alpha+1},
\) which is not integrable at infinity. Thus, as in
\cite{caffarelli2013regularity}, we remove the nonintegrable
far-field drift by working in a moving frame. Let
\(S(\tau)\) be defined by
\[
S'(\tau)=-\nabla L^{-1}\rho(S(\tau),\tau),
\qquad
S(0)=0,
\]
and set
\[
\widetilde \rho(x,\tau):=\rho(x+S(\tau),\tau).
\]
By translation invariance, we have
\[
\nabla L^{-1}\rho(x+S(\tau),\tau)
=
\nabla L^{-1}\widetilde \rho(x,\tau),
\qquad
\nabla L^{-1}\rho(S(\tau),\tau)
=
\nabla L^{-1}\widetilde \rho(0,\tau).
\]
Differentiating \(\widetilde \rho\) in time and using
the equation for \(\rho\), we obtain
\[
\partial_\tau\widetilde \rho
=
\nabla\cdot\left(
\widetilde \rho\nabla L^{-1}\widetilde \rho
\right)
-
\nabla L^{-1}\widetilde \rho(0,\tau)
\cdot\nabla\widetilde \rho.
\]
Since \(\nabla L^{-1}\widetilde \rho(0,\tau)\) is
independent of \(x\), this is equivalently
\begin{equation*}
\partial_\tau\widetilde \rho
=
\nabla\cdot\left(
\widetilde \rho
\left[
\nabla L^{-1}\widetilde \rho(x,\tau)
-
\nabla L^{-1}\widetilde \rho(0,\tau)
\right]
\right).
\end{equation*}
At the level of the velocity kernel, this change
of variables subtracts the drift at the base point.
Thus, the far-field contribution is governed by
the difference \(
\nabla W(x-y)-\nabla W(-y).
\) For \(x\) in a fixed bounded set and \(|y|\) large,
the mean-value theorem gives
\[
|\nabla W(x-y)-\nabla W(-y)|
\leq
|x|\sup_{\theta\in[0,1]}
|D^2W(-y+\theta x)|.
\]
Since
\(W(x)=|x|^{-s}\Omega(x/|x|)\), with
\(\Omega\in C^\infty(\mathbb S^{d-1})\), its
second derivatives satisfy
\[
|D^2W(z)|
\leq C_{\Omega,d,s}|z|^{-s-2},
\qquad z\neq0.
\]
Hence, for \(|y|\) sufficiently large compared
with the fixed spatial scale of \(x\),
\[
|\nabla W(x-y)-\nabla W(-y)|
\leq
C_{\Omega,d,s}|x||y|^{-s-2}
=
C_{\Omega,d,s}|x||y|^{-d-\alpha}.
\]
Since \(|y|^{-d-\alpha}\) is integrable at infinity,
the moving frame replaces the nonintegrable
far-field drift by an integrable tail. The
near-field estimates are unchanged: since
\(0<\alpha<1\), \(\nabla W\) is locally integrable,
and the singular terms are controlled by
\eqref{eq:velocity-kernel-bound-holder} and
the bilinear form \(\mathcal B_A\), as in the
previous case. The only modification occurs
in the far field, where the corrected kernel
has the integrable decay obtained above.
Consequently, in the transported variables,
the oscillation estimates are independent
of the global \(L^1\)-norm, which is used
only to control the coordinate shift.

Finally, we control the geometric distortion
produced by the moving frame invoking
the kernel representation of \(L^{-1}\). More preciesly, 
\begin{align*}
|S'(\tau)|
&=
\left|
\int_{\mathbb R^d}
\nabla W(S(\tau)-y)\rho(y,\tau)dy
\right|\\
&\leq
C_{\Omega,d,s}
\int_{|S(\tau)-y|\leq1}
|S(\tau)-y|^{-d-\alpha+1}
\rho(y,\tau)dy\\
&\quad+
C_{\Omega,d,s}
\int_{|S(\tau)-y|>1}
|S(\tau)-y|^{-d-\alpha+1}
\rho(y,\tau)dy\\
&\leq
C_{\Omega,d,s}
\left(
\|\rho(\cdot,\tau)\|_{L^\infty}
+
\|\rho(\cdot,\tau)\|_{L^1}
\right).
\end{align*}
where we used the fact that $|z|^{-d-\alpha+1}$ is integrable in the near field since $0<\alpha<1$, and that $|z|^{-d-\alpha+1}\leq1$ for $|z|>1$. Recall that the initial normalization gives \(0\leq \rho\leq1\) and \(\|\rho(\cdot,\tau)\|_{L^1}=1,\)
we obtain
\[
|S'(\tau)|\leq C_{\Omega,d,s},
\qquad \tau\in[-4,0].
\]
Hence, we have
\[
|S(\tau)-S(\sigma)|
\leq
C_{\Omega,d,s}|\tau-\sigma|,
\qquad \tau,\sigma\in[-4,0].
\]
Thus, the moving frame is uniformly Lipschitz
in time, and the geometric distortion of the
transported cylinders is controlled. After
this initial change of coordinates, the
successive rescalings and oscillation reductions
follow \cite[Section~12]{caffarelli2013regularity}.
Since the change of coordinates is Lipschitz,
the resulting H\"older estimate transfers back
to the original variables, yielding local
space-time H\"older continuity in the case
\(0<\alpha<1\).

\paragraph{The case \(\alpha=1\).}
In the critical regime \(s=d-1\), the velocity-kernel estimate \eqref{eq:velocity-kernel-bound-holder} reduces to
\begin{equation}\label{eq:critical-velocity-kernel-holder}
|\nabla W(h)|\leq C_{\Omega,d}|h|^{-d},
\qquad
h\neq0.
\end{equation}
Since this kernel is logarithmically nonintegrable at both the origin and infinity, the direct estimate for $1<\alpha<2$ and the single moving frame for $0<\alpha<1$ both fail to apply. Instead, we adapt the iterative partial-transport correction from \cite[Section~9]{caffarelli2016regularity}, highlighting only the modifications required by the anisotropic kernel.

Fix \(0<\varepsilon<1/2\) sufficiently small for the exterior barriers, and choose \(0<\nu<\min\{\varepsilon,1/4\}\) for the intermediate-value barriers in \eqref{eq:intermediate-barriers-holder}. We choose \(\varepsilon_0\in(0,1/8)\) sufficiently small for the near-field absorption below. Let \(\chi\in C^\infty(\mathbb R^d)\) be radial such that

\begin{equation}\label{eq:critical-cutoff-holder}
0\leq\chi\leq1,
\qquad
\chi\equiv0\ \text{in }B_{\varepsilon_0},
\qquad
\chi\equiv1\ \text{in }\mathbb R^d\setminus B_{2\varepsilon_0}, \qquad  \|\nabla\chi\|_{L^\infty}\leq C\varepsilon_0^{-1}
\end{equation}
where \(C\) independent of \(\varepsilon_0\). Let us split the velocity kernel as \(\nabla W=\nabla W(1-\chi)+\nabla W\chi\). Since the first part is supported in \(B_{2\varepsilon_0}\), it follows from \eqref{eq:velocity-kernel-bound-holder} that
\begin{equation}\label{eq:critical-near-kernel-holder}
|\nabla W(h)|(1-\chi(h))
\leq
C|h|K_A(h)\mathbf 1_{\{|h|\leq2\varepsilon_0\}}
\leq
C\varepsilon_0K_A(h).
\end{equation}
This replaces the isotropic radial estimate in the critical near-field regime and allows the corresponding contribution to be absorbed by the anisotropic Dirichlet form.

For the truncated far-field kernel, define
\begin{equation}\label{eq:critical-truncated-kernel-holder}
G_\chi(h):=\nabla W(h)\chi(h).
\end{equation}
Since \(W\) is homogeneous of degree \(-(d-1)\) and \(\Omega\in C^\infty(\mathbb S^{d-1})\), we have
\begin{equation}\label{eq:critical-anisotropic-derivatives-holder}
|D^2W(h)|\leq C_{\Omega,d}|h|^{-d-1},
\qquad
h\neq0.
\end{equation}
which implies
\begin{equation}\label{eq:critical-DGchi-holder}
|DG_\chi(h)|
\leq
C_{\Omega,d,\varepsilon_0}(1+|h|)^{-d-1},
\qquad
h\in\mathbb R^d.
\end{equation}

Indeed, by the homogeneity of \(W\) and the ellipticity of \(K_A\), we have \(|D^2W(h)|\leq C_{\Omega,d}|h|^{-d-1}\leq C K_A(h)\) for \(h\neq0\). Since
\[
DG_\chi=D^2W\chi+\nabla W\otimes\nabla\chi,
\]
and \(\chi\equiv0\) in \(B_{\varepsilon_0}\), while \(\nabla\chi\) is supported in \(B_{2\varepsilon_0}\setminus B_{\varepsilon_0}\), it follows that \(DG_\chi\) is bounded near the origin and decays as \(|h|^{-d-1}\) at infinity, yielding \eqref{eq:critical-DGchi-holder}.

For a nonnegative function \(v(\cdot,\tau)\in L^1(\mathbb R^d)\), define the truncated drift
\begin{equation}\label{eq:critical-truncated-drift-holder}
b_\chi[v](\tau)
:=
\int_{\mathbb R^d}\nabla W(-y)\chi(y)v(y,\tau)dy.
\end{equation}
The cutoff removes the singularity at the origin. Moreover, from \eqref{eq:critical-cutoff-holder} and \eqref{eq:critical-velocity-kernel-holder}, we have
\begin{equation*} 
|b_\chi[v](\tau)|
\leq
C_{\Omega,d,\varepsilon_0}
\|v(\cdot,\tau)\|_{L^1(\mathbb R^d)}.
\end{equation*}
Let \(S\) and \(\widetilde \rho\) be defined by
\begin{equation*}
\widetilde \rho(x,\tau):=\rho(x+S(\tau),\tau),
\qquad
S'(\tau)=-b_\chi[\widetilde \rho](\tau),
\qquad
S(0)=0.
\end{equation*}
By translation invariance, \(\widetilde \rho\) satisfies
\begin{equation*}
\partial_\tau\widetilde \rho
=
\nabla\cdot\left(\widetilde \rho\mathcal H_\chi[\widetilde \rho]\right),
\qquad
\mathcal H_\chi[v](x,\tau)
:=
\nabla L^{-1}v(x,\tau)-b_\chi[v](\tau),
\end{equation*}
where \(\mathcal H_\chi[v]\) denotes the corrected velocity field. Recalling the kernel representation of \(\nabla L^{-1}\) used in \eqref{eq:grad-Lminus-one-PV-holder}, together with \eqref{eq:critical-cutoff-holder}, \eqref{eq:critical-truncated-kernel-holder}, and \eqref{eq:critical-truncated-drift-holder}, we split \(\nabla W=\nabla W(1-\chi)+G_\chi\). Since \(W\) is even and \(\chi\) is radial, \(\nabla W(1-\chi)\) is odd, and its principal-value integral against the constant \(v(x,\tau)\) vanishes. Subtracting the truncated drift \(b_\chi[v]\), we obtain, initially for smooth integrable functions,
\begin{equation*}
\mathcal H_\chi[v]
=
\mathcal H_{\mathrm{near}}[v]
+
\mathcal H_{\mathrm{far}}[v],
\end{equation*}
where
\begin{align*}
\mathcal H_{\mathrm{near}}[v](x,\tau)
&:=
\pv\int_{\mathbb R^d}
\nabla W(x-y)(1-\chi(x-y))
\bigl(v(y,\tau)-v(x,\tau)\bigr)dy,
\\
\mathcal H_{\mathrm{far}}[v](x,\tau)
&:=
\int_{\mathbb R^d}
\bigl[G_\chi(x-y)-G_\chi(-y)\bigr]v(y,\tau)dy.
\end{align*}

For \(x\) in a fixed ball, the mean-value theorem and \eqref{eq:critical-DGchi-holder} imply
\begin{equation*}
|G_\chi(x-y)-G_\chi(-y)|
\leq
C|x|(1+|y|)^{-d-1}.
\end{equation*}
Thus, under the normalized exterior bound \(0\leq v(y,\tau)\leq1+(|y|^\varepsilon-2)_+\), we obtain
\[
|\mathcal H_{\mathrm{far}}[v](x,\tau)|
\leq
C|x|
\int_{\mathbb R^d}
(1+|y|)^{-d-1}(1+|y|^\varepsilon)dy
\leq
C|x|.
\]
In particular, the corrected far-field velocity is locally bounded independently of the global \(L^1\)-norm. From now on, we drop the tilde and denote the transported solution again by \(\rho\).

We now explain how the preceding estimates enter the oscillation lemmas. For the upper barriers \(\varphi_k\), set
\[
\rho_k^+=(\rho-\varphi_k)_+,
\qquad
\rho_k^-=(\rho-\varphi_k)_-,
\qquad
I_k(\tau):=\int_{\{\rho_k^+(\tau)>0\}}(\rho_{k-1}^+)^2dx.
\]
In the transported variables, the drift term corresponding to \eqref{eq:Q-upper-holder} is evaluated using the corrected velocity \(\mathcal H_\chi[\rho]\) and decomposed into its near- and far-field contributions. Using \eqref{eq:critical-near-kernel-holder}, the smallness of \(\varepsilon_0\) specified above, and consider \(\gamma>0\) fixed sufficiently small, we obtain
\begin{align}
|Q_{k,\mathrm{near}}|
&\leq
\gamma\mathcal B_A(\rho_k^+,\rho_k^+)
+
\gamma\mathcal B_A(\rho_k^+,\rho_k^-)
+
C_\gamma16^kI_k(\tau).
\label{eq:critical-near-Q-holder}
\end{align}
For the far-field contribution, the active set is contained in a fixed ball and \(|\nabla\varphi_k|/\varphi_k\leq C\) on \(\{\rho_k^+>0\}\). Hence the local bound above gives
\[
|\mathcal H_{\mathrm{far}}[\rho](x,\tau)|
\leq C
\qquad
\text{on }\{\rho_k^+(\tau)>0\}.
\]
Using \(\rho_k^+\leq \rho_{k-1}^+\) and the separation between consecutive barriers, we obtain
\begin{align}
|Q_{k,\mathrm{far}}|
&\leq
C\int_{\{\rho_k^+>0\}}\rho_k^+dx
\leq
C|\{\rho_k^+>0\}|^{1/2}I_k(\tau)^{1/2}
\leq
C4^kI_k(\tau)
\leq
C16^kI_k(\tau).
\label{eq:critical-far-Q-holder}
\end{align}

The barrier interaction is treated exactly as in \eqref{eq:upper-bilinear-error-holder}, using only \eqref{eq:KA-ellipticity-holder}, the Lipschitz bound for \(\varphi_k\), and \(K_A(h)\lesssim |h|^{-d-1}\). Together with \eqref{eq:upper-barrier-identity-holder}, \eqref{eq:critical-near-Q-holder}, and \eqref{eq:critical-far-Q-holder}, this yields \eqref{eq:upper-local-energy-holder}, with constants independent of the total mass. The lower-barrier and intermediate-value arguments follow from the same near-field/far-field decomposition.

After each rescaling, the cutoff in the transport correction changes. Following \cite[Sections~9.2--9.5]{caffarelli2016regularity}, we introduce a relative correction at each step. Let \(\mu\in(0,1)\) denote the oscillation-reduction constant in the critical regime, and set \(\vartheta:=1-\mu\). We choose \(B>1\) sufficiently large so that
\[
T:=\vartheta B>1,
\qquad
\vartheta B^\varepsilon>1,
\]
and the exterior barriers and transported cylinders remain admissible throughout the iteration, as in \cite[Sections~7 and~9]{caffarelli2016regularity}.

Starting from \(v_0:=\rho\), the scaling step, before the next transport correction, is given by
\[
v_k(x,\tau)
=
\vartheta v_{k+1}(Bx,T\tau).
\]
The underlying unconfined equation is invariant under this scaling, since \(\nabla W(Bh)=B^{-d}\nabla W(h)\) and \(T=\vartheta B\). The transport correction, however, inherits a rescaled cutoff, as shown below.

For a cutoff \(\Phi\), let \(b_\Phi\) be defined as in \eqref{eq:critical-truncated-drift-holder}. A change of variables and the critical homogeneity give
\begin{align}
b_\Phi[v_k](\tau)
&=
\vartheta
\int_{\mathbb R^d}
\nabla W(-y)\Phi(y)v_{k+1}(By,T\tau)dy
\nonumber\\
&=
\vartheta
\int_{\mathbb R^d}
\nabla W(-z)\Phi(z/B)v_{k+1}(z,T\tau)dz
=
\vartheta b_{\Phi(\cdot/B)}[v_{k+1}](T\tau).
\label{eq:critical-drift-scaling-holder}
\end{align}
Thus, the inherited correction at scale \(k+1\) uses the cutoff \(\chi(\cdot/B)\), whereas the natural correction uses \(\chi\). We therefore introduce
\begin{equation}\label{eq:critical-relative-cutoff-holder}
\chi_{\mathrm{rel}}(y)
:=
\chi(y)-\chi(y/B),
\qquad
\supp\chi_{\mathrm{rel}}
\subset
\overline B_{2B\varepsilon_0}\setminus B_{\varepsilon_0}.
\end{equation}
For \(k\geq0\), define the relative drift by
\begin{align}
b_{k+1}^{\mathrm{rel}}(\tau)
:=
b_\chi[v_{k+1}](\tau)
-
b_{\chi(\cdot/B)}[v_{k+1}](\tau)
=
\int_{\mathbb R^d}
\nabla W(-y)\chi_{\mathrm{rel}}(y)v_{k+1}(y,\tau)dy.
\label{eq:critical-relative-drift-holder}
\end{align}

As in \cite[Sections~7 and~9]{caffarelli2016regularity}, the oscillation iteration is performed under the normalized exterior bound
\[
0\leq v_j(y,\tau)
\leq 1+(|y|^\varepsilon-2)_+.
\]
The condition \(\vartheta B^\varepsilon>1\) controls the leading-order growth of the exterior barrier under rescaling, while its preservation in bounded regions is incorporated into the oscillation-reduction argument. Since \(\chi_{\mathrm{rel}}\) is supported in the annulus described in \eqref{eq:critical-relative-cutoff-holder}, we have
\[
v_{k+1}(y,\tau)
\leq 1+(2B\varepsilon_0)^\varepsilon
\qquad
\text{for }\varepsilon_0\leq|y|\leq2B\varepsilon_0.
\]
Thus, using \eqref{eq:critical-relative-drift-holder} and \eqref{eq:critical-velocity-kernel-holder}, we obtain
\begin{align}
|b_{k+1}^{\mathrm{rel}}(\tau)|
&\leq
C_{\Omega,d}
\bigl(1+(2B\varepsilon_0)^\varepsilon\bigr)
\int_{\varepsilon_0\leq|y|\leq2B\varepsilon_0}
|y|^{-d}dy
\nonumber\\
&=
C_{\Omega,d}
\bigl(1+(2B\varepsilon_0)^\varepsilon\bigr)
|\mathbb S^{d-1}|\log(2B)
=:C_B.
\label{eq:critical-relative-drift-bound-holder}
\end{align}
Since \(B\), \(\varepsilon_0\), and \(\varepsilon\) are fixed, this bound is uniform in \(k\) and independent of the masses of the rescaled solutions.

Combining \eqref{eq:critical-drift-scaling-holder} and \eqref{eq:critical-relative-drift-holder}, we obtain
\[
b_\chi[v_k](\tau)
=
\vartheta\left(
b_\chi[v_{k+1}](T\tau)
-
b_{k+1}^{\mathrm{rel}}(T\tau)
\right).
\]
Iterating this identity yields
\[
b_\chi[v_0](\tau)
+
\sum_{j=1}^{k-1}
\vartheta^j b_j^{\mathrm{rel}}(T^j\tau)
=
\vartheta^{k-1}
b_\chi[v_{k-1}](T^{k-1}\tau).
\]

Following \cite[Sections~9.3--9.5]{caffarelli2016regularity}, the successive relative corrections contribute geometric factors \(\vartheta^j\) when expressed in the initial coordinates. For every \(k\geq2\) and \(\tau\in[-4T^{-(k-1)},0]\), we have \(T^j\tau\in[-4,0]\) for \(1\leq j\leq k-1\). Thus, it follows that
\begin{align*}
\left|
b_\chi[v_0](\tau)
+
\sum_{j=1}^{k-1}
\vartheta^j b_j^{\mathrm{rel}}(T^j\tau)
\right|
&\leq
C\sup_{\sigma\in[-4,0]}
\|v_0(\cdot,\sigma)\|_{L^1(\mathbb R^d)}
+
C_B\sum_{j=1}^{k-1}\vartheta^j\\
&\leq
C+\frac{\vartheta C_B}{1-\vartheta},
\end{align*}
where we used that the initial normalization and the moving-frame translation preserve the unit mass of \(v_0\). Thus, the initial drift is controlled by the mass of \(v_0\), while the subsequent relative corrections are uniformly bounded, independently of the masses of the rescaled solutions. Since their weighted contributions form a convergent geometric series, the accumulated transport correction remains uniformly bounded. The remaining oscillation and geometric arguments follow \cite[Sections~9--10]{caffarelli2016regularity}, yielding local space-time H\"older continuity in the critical case \(\alpha=1\).

Combining the three regimes, we obtain local space-time H\"older continuity, with a common exponent \(\eta\in(0,1)\), after decreasing it if necessary. Finally, undoing the initial translation and scaling gives the stated estimate for every \(T_1<t_0<t_1<T_2\).

\section*{Acknowledgements} JAC was supported by the Advanced Grant Nonlocal-CPD (Nonlocal PDEs for Complex Particle Dynamics: Phase Transitions, Patterns and Synchronization) of the European Research Council Executive Agency (ERC) under the European Union Horizon 2020 research and innovation programme (grant agreement No. 883363) and  partially supported by the EPSRC EP/V051121/1. JAC was partially supported by the ``Maria de Maeztu'' Excellence Unit IMAG, reference CEX2020-001105-M, funded by MCIN/AEI/10.13039/501100011033/. DPG was partially supported by the Zuckerman STEM Leadership Program and by the FWF ESPRIT project ESP3828725. RMMG was supported by the FAPESP project 2025/08196-9 and 2024/00744-4.

\section*{Statement on generative AI}

The authors used ChatGPT-6 Astra to assist the exploration of implications among the assumptions \eqref{hw}, \eqref{hl} and \eqref{ha}. In particular, the proof of Lemma \ref{lem:HL-implies-HW} was done with its assistance. No other place in this work has used AI tools except for final editing.

\bibliographystyle{alpha}
\bibliography{bibliography.bib}

\end{document}